\documentclass[12pt,a4paper]{article}
\usepackage{color}
\usepackage{amsmath, amssymb, theorem, latexsym}
\usepackage{graphicx}
\usepackage{epsf}
\usepackage{fullpage}
\usepackage[pdftex]{hyperref}
\usepackage{multicol}
\newtheorem{theorem}{Theorem}[section]
\newtheorem{lemma}[theorem]{Lemma}
\newtheorem{proposition}[theorem]{Proposition}

\newtheorem{corollary}[theorem]{Corollary}

\newtheorem{properties}[theorem]{Properties}

\numberwithin{equation}{section}
\numberwithin{figure}{section}
\numberwithin{table}{section}

\allowdisplaybreaks
\newcommand{\clr}{\color{red}}

\newcommand{\noib}{\noindent $\bullet $~}

\newcommand{\ie}{i.e.}

\newcommand{\resp}{\emph{resp. }}
\newcommand{\etc}{\emph{etc\,}}

\newcommand{\cB}{\mathcal{B}}

\newcommand{\cD}{\mathcal{D}}
\newcommand{\cE}{\mathcal{E}}
\newcommand{\cF}{\mathcal{F}}

\newcommand{\cH}{\mathcal{H}}

\newcommand{\cL}{\mathcal{L}}

\newcommand{\cS}{\mathcal{S}}
\newcommand{\cT}{\mathcal{T}}

\newcommand{\E}{\mathbb{E}}

\newcommand{\Nb}{\mathbb{N}^{\bullet}}

\newcommand{\R}{\mathbb{R}}

\newcommand{\T}{\mathbb{T}}
\newcommand{\Z}{\mathbb{Z}}

\title{ Courant-sharp eigenvalues for the equilateral torus, and for the equilateral triangle}

\author{P. B\'erard \\ Institut Fourier, Universit\'{e} de Grenoble and CNRS, B.P.74,
\\ F 38402 Saint Martin d'H\`{e}res Cedex, France.\\
and \\
B. Helffer \\
Laboratoire de Math\'ematiques, Univ. Paris-Sud 11 and CNRS,\\
F 91405 Orsay Cedex, France, and\\
Laboratoire  de Math\'{e}matiques Jean Leray, Universit\'{e} de Nantes.}

\date{June 25, 2015}

\begin{document}
\maketitle

\medskip
{\large{To the memory of Louis Boutet de Monvel}}
\bigskip

\begin{abstract}
We address the question of determining the eigenvalues $\lambda_n$ (listed in nondecreasing order, with multiplicities) for which Courant's nodal domain theorem is sharp \ie, for which there exists an associated eigenfunction with $n$ nodal domains (Courant-sharp eigenvalues). Following ideas going back to Pleijel (1956), we prove that the only Courant-sharp eigenvalues of the flat equilateral torus are the  first and second, and that the only Courant-sharp Dirichlet eigenvalues of the equilateral triangle are the  first, second, and fourth eigenvalues.  In the last section we sketch similar results for the right-angled isosceles triangle and for the hemiequilateral triangle.
\end{abstract}


\section{Introduction}\label{S-intro}

Let $(M,g)$ be a compact Riemannian manifold, and $\Delta$ the nonpositive Laplace-Beltrami operator. Consider the eigenvalue problem $-\Delta u = \lambda u$, with Dirichlet boundary condition in case $M$ has a boundary. Write the eigenvalues in nondecreasing order, with multiplicities,
$$
\lambda_1(M) < \lambda_2(M) \le \lambda_3(M) \le \cdots \,.
$$

Courant's theorem \cite{Cou} states that an eigenfunction $u$, associated with the $k$-th eigenvalue $\lambda_k$, has at most $k$ nodal domains. We say that $\lambda_k$ is \emph{Courant-sharp} if there exists an eigenfunction $u$, associated with $\lambda_k$, with exactly $k$ nodal domains. Pleijel \cite{Pl} (see also \cite{BeHe1}) has shown that the only Dirichlet eigenvalues of a square which are Courant-sharp are the first, second and fourth eigenvalues. L\'{e}na \cite{Len} recently proved that the only Courant-sharp eigenvalues of the square flat torus are the first and second eigenvalues.

In this paper,  following \cite{Pl} and \cite{Len}, we determine the  Courant-sharp eigenvalues for the equilateral torus and the equilateral triangle (with Dirichlet boundary condition).

\begin{theorem}\label{PETO-1}
The only Courant-sharp eigenvalues of the equilateral torus are the first and second eigenvalues.
\end{theorem}

\begin{theorem}\label{PETR-2}
The only Courant-sharp Dirichlet eigenvalues of the equilateral triangle are the first, second, and fourth eigenvalues.
\end{theorem}

\medskip

The Dirichlet eigenvalues and eigenfunctions of the equilateral triangle were originally studied by Lam\'{e} \cite{Lam}, and later in the papers \cite{Pin1, Ber, Pin2}.  We refer to \cite{McC} for a detailed account, and to \cite[\S 6]{KeRu} for another approach.

In \cite{Ber}, the author determined the Dirichlet and Neumann eigenvalues and eigenfunctions for the fundamental domains of crystallographic groups, namely the fundamental domains (or \emph{alcoves})  of the action of the affine Weyl group $W_a$ generated by the reflections associated with a root system \cite{Bou}. In the two-dimensional case, the domains are the rectangle, the equilateral triangle, the right isosceles triangle, and the  hemiequilateral triangle, \ie, with angles $(30,60,90)$ degrees. The same ideas \cite{BerBes} can be applied to determine the Dirichlet and Neumann eigenvalues of the spherical domains obtained as intersections of the round sphere with the fundamental domain (or \emph{chamber}) of the Weyl group $W$ associated with a root system. \medskip

In Section~\ref{S-GF}, we describe the proper geometric framework to determine the eigenvalues of the equilateral triangle, following the ideas of \cite{Ber}, but without references to root systems. For the convenience of the reader, we use the same notation as in \cite{Bou, Ber}. \medskip

In Section~\ref{S-ETO},  following \cite{Len}, we study  the Courant-sharp eigenvalues of the equilateral torus, and prove Theorem~\ref{PETO-1}.\medskip

Starting from Section~\ref{S-ETR}, we describe the spectrum of the equilateral triangle, and the symmetries of the eigenfunctions. In Section \ref{S5}, following Pleijel's approach  \cite{Pl}, we reduce the question of Courant-sharp Dirichlet eigenvalues of the equilateral triangle to the analysis of only two eigenspaces associated with the eigenvalues $\lambda_5$ and $\lambda_7$.  The number of nodal domains of eigenfunctions in these spaces  is studied in Sections~\ref{SS-E13} and \ref{SS-E23} respectively, after a presentation of the strategy in Section~\ref{sProl}.

In Section~\ref{S-OT}, we sketch similar results for the right-angled isosceles triangle and for the hemiequilateral triangle. In both cases, the only Courant-sharp eigenvalues are the first and second.

The authors would like to thank Virginie Bonnaillie-No\"{e}l for communicating her numerical computations of eigenvalues and nodal sets for the triangular domains.

\section{The geometric framework}\label{S-GF}

We call $\E^2$ the vector (or affine) space $\R^2$ with the canonical basis $\{e_1 = (1,0), e_2=(0,1)\}$, and with the standard  Euclidean inner product $\langle \cdot,\cdot \rangle$ and norm $|\cdot|$.

Define the vectors
\begin{equation}\label{GF-2}
\alpha_1 = \left( 1, -\frac{1}{\sqrt{3}} \right),\,
\alpha_2 = \left( 0, \frac{2}{\sqrt{3}} \right),\, \text{~and~}
\alpha_3 = \left( 1, \frac{1}{\sqrt{3}}\right).
\end{equation}

The vectors $\alpha_1, \alpha_2$ span $\E^2$. Furthermore,
$|\alpha_i| = \frac{2}{\sqrt{3}}$, $\measuredangle(\alpha_1,\alpha_2) = \frac{2\pi}{3}$, and \goodbreak $\measuredangle(\alpha_1,\alpha_3) =\measuredangle(\alpha_3,\alpha_2) = \frac{\pi}{3}$  (here $\measuredangle$ denotes the angle between two vectors). \medskip

Introduce the vectors
$$
\alpha_{i}^{\vee} = \frac{2}{\langle \alpha_i,\alpha_i \rangle} \alpha_i\,, \text{for~} i=1, 2, 3\,.
$$
Then,
\begin{equation}\label{GF-4}
\alpha_{1}^{\vee} = \left( \frac{3}{2}, -\frac{\sqrt{3}}{2}\right),\,
\alpha_{2}^{\vee} = \left( 0, \sqrt{3}\right),\, \text{~and~}
\alpha_{3}^{\vee} = \left( \frac{3}{2}, \frac{\sqrt{3}}{2}\right).
\end{equation}

For $i \in \{1,2,3\}$, introduce the lines
\begin{equation}\label{GF-6}
L_i = \left\{ x \in \E^2 ~|~ \langle x,\alpha_i \rangle = 0\right\},
\end{equation}
and the orthogonal symmetries with respect to these lines,
\begin{equation}\label{GF-8}
s_i(x) = x - \frac{2\langle x,\alpha_i \rangle}{\langle \alpha_i,\alpha_i \rangle}  \alpha_i \,.
\end{equation}

Let $W$ be the group  generated by the symmetries $\{s_1, s_2, s_3\}$. \medskip

More generally, for $i \in \{1,2,3\}$ and $k \in \Z$, consider the lines
\begin{equation}\label{GF-6a}
L_{i,k} = \left\{ x \in \E^2 ~|~ \langle x,\alpha_i \rangle = k\right\},
\end{equation}
and the orthogonal symmetries with respect to these lines,
\begin{equation}\label{GF-8a}
s_{i,k}(x) = x - \frac{2\langle x,\alpha_i \rangle}{\langle \alpha_i,\alpha_i \rangle} \alpha_i + k \alpha_{i}^{\vee} = s_i(x)+k\alpha_{i}^{\vee}\,.
\end{equation}

Let $W_a$ be the group generated by the symmetries $\{s_{i,k} ~|~ i = 1, 2, 3\, \text{~and~}k \in \Z\}$.

These groups are called respectively the \emph{Weyl group} and the \emph{affine Weyl group} \cite{Bou}. \medskip

Call $\Gamma$ the group generated by $\alpha_{1}^{\vee}$ and $\alpha_{2}^{\vee}$ (see the notation $Q(R^{\vee})$ in \cite{Bou,Ber}) \ie,
\begin{equation}\label{GF-10}
\Gamma = \Z \, \alpha_{1}^{\vee} \bigoplus \Z \, \alpha_{2}^{\vee}\,.
\end{equation}

\begin{figure}
  \centering
  \includegraphics[width=8cm]{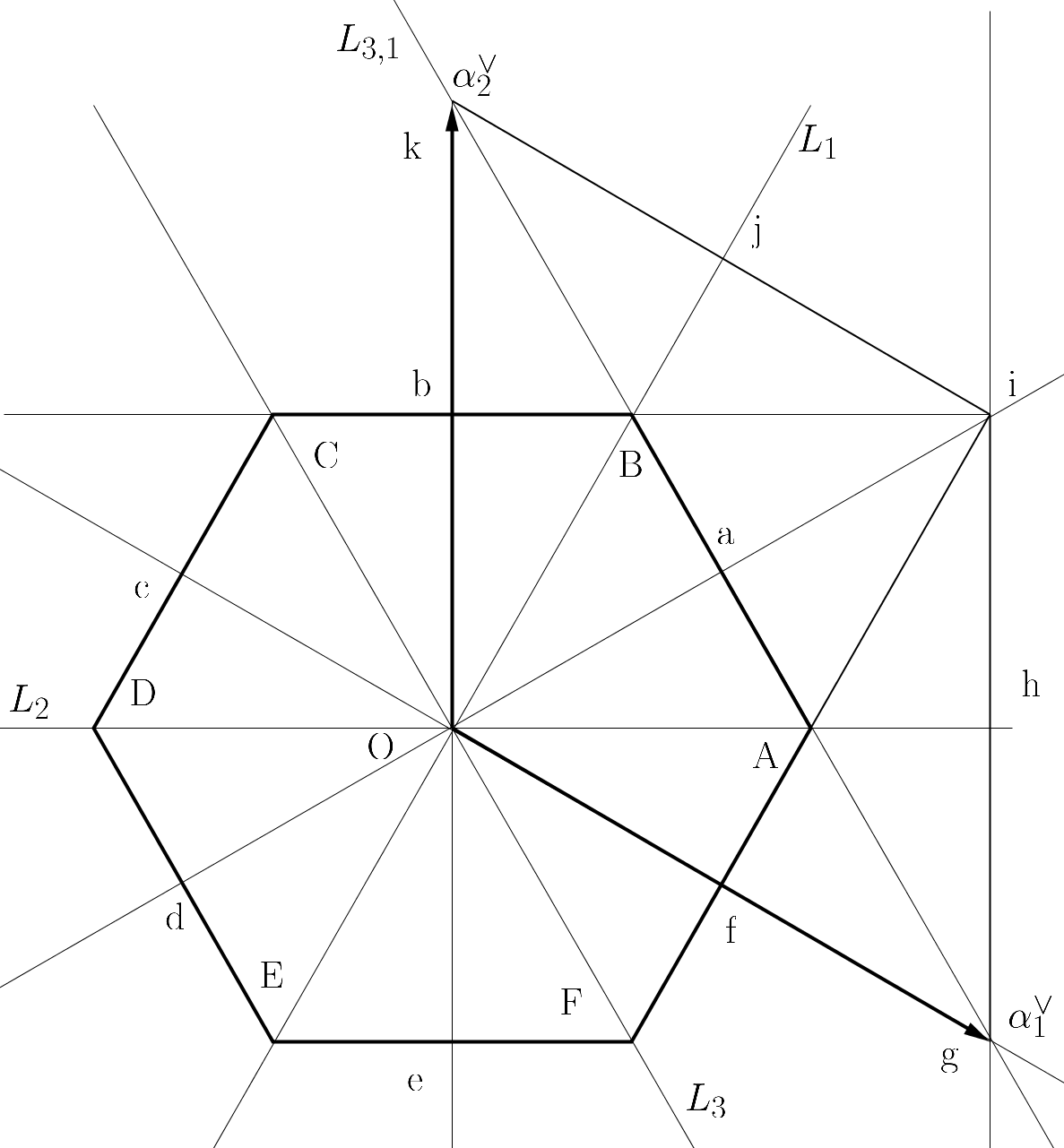}
  \caption{The geometric framework}\label{FGF-1}
\end{figure}

The following properties are easy to establish (see Figure~\ref{FGF-1}).

\begin{proposition}\label{GF-P1}
With the above notation, the following properties hold.
\vspace{-3mm}
\begin{enumerate}
  \item The group $W$ has six elements; it is isomorphic to the symmetric group in three letters $\cS_3$.
  \item The group $W$ acts simply transitively on the connected components of $\E^2 \setminus L_1 \cup L_2 \cup L_3$ (the Weyl chambers).
  \item The group $W_a$ is the semi-direct product $\Gamma \rtimes W$. A fundamental domain for the action of $W_a$ on $\E^2$ is the closed equilateral triangle $\cT$, with vertices $O=(0,0)$, $A=(1,0)$, and $B=(\frac{1}{2},\frac{\sqrt{3}}{2})$.
  \item The group $W$ acts simply transitively on the equilateral triangles which tile the regular hexagon $[A,B,\ldots,E,F]$.
  \item The closed regular hexagon $[A,B,\ldots,E,F]$ is a fundamental domain for the action of the lattice $\Gamma$ on $\E^2$.
\end{enumerate}
\end{proposition}%

\begin{figure}[htb]
  \centering
  \includegraphics[width=8cm]{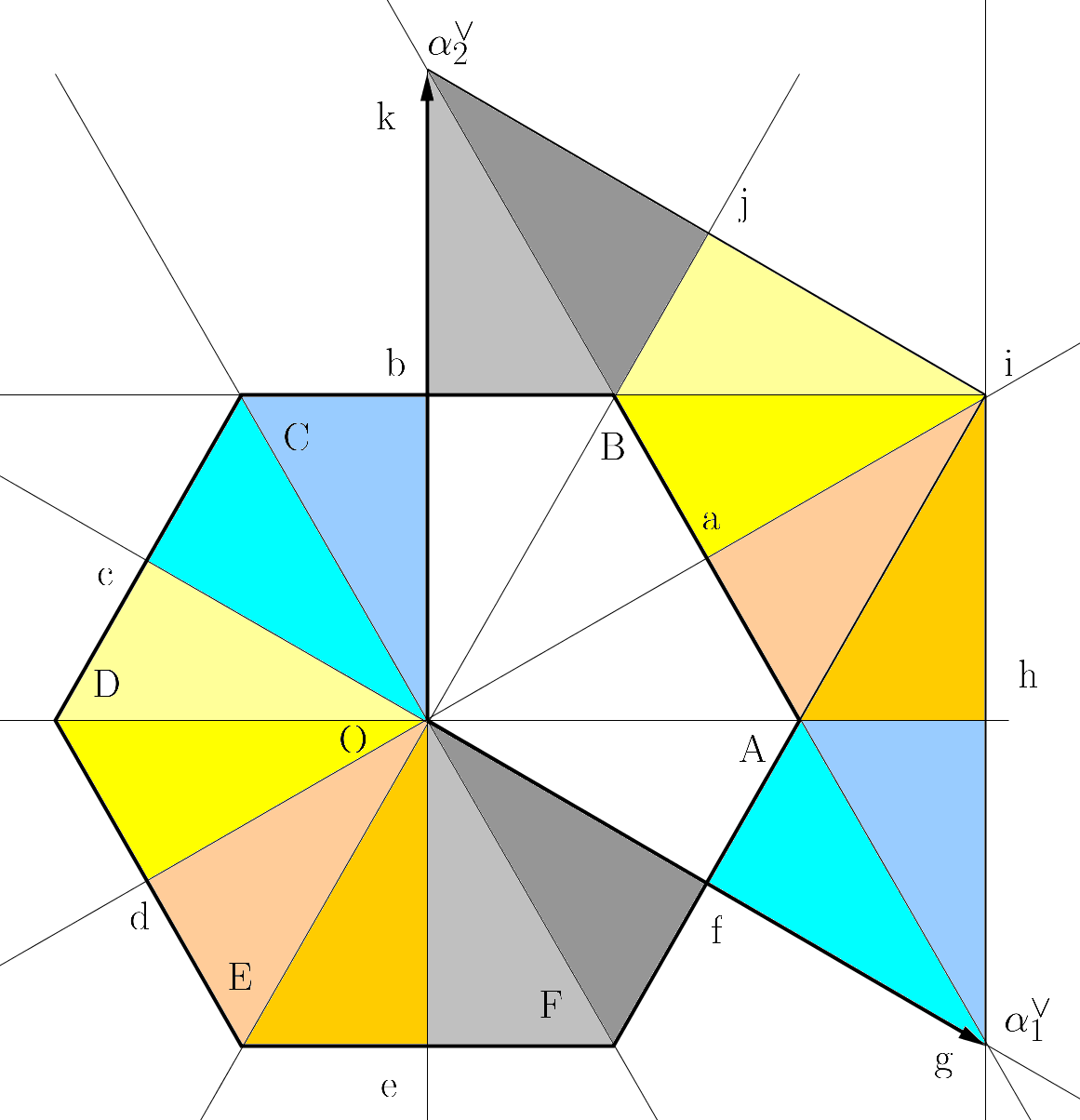}
  \caption{The hexagon is a fundamental domain for $\Gamma$}\label{FGF-2}
\end{figure}

\textbf{Proof}. \\
(1)  Clearly, $s_1^2 = 1$ (where 1 is the identity), $s_2^2 =1$, $s_3^2=1$, and $s_3 = s_1\circ s_2\circ s_1 = s_2\circ s_1\circ s_2\,$. Furthermore, $s_1\circ s_2$, $s_3\circ s_1$ and $s_2\circ s_3$ are equal to the rotation with center $O$ and angle $\frac{2\pi}{3}$; $s_2\circ s_1$, $s_1\circ s_3$ and $s_3\circ s_2$ are equal to the rotation with center $O$ and angle $-\frac{2\pi}{3}$.  It follows that
\begin{equation*}
W = \left\{ 1, s_1, s_2, s_3, s_1 \circ s_2, s_2 \circ s_1 \right\}.
\end{equation*}
To see that $W$ is isomorphic to $\cS_3$, look at its action on the lines $\R\,\alpha_i\,$, $i=1, 2, 3\,$. \\
(2) Clear.\\
(3) Write the action of $(\gamma,\sigma) \in \Gamma\times W$ on $\E^2$, as
$$
x \mapsto (\gamma,\sigma)\cdot x = \sigma(x) + \gamma\,.
$$
Then, $(\gamma,\sigma).(\gamma',\sigma') = (\gamma+\sigma(\gamma'), \sigma\circ\sigma')$. Clearly, $s_{i,k}$ can be written as $(k\alpha_{i}^{\vee},s_i)$ and hence, $W_a \subset \Gamma \rtimes W$.
To prove the reverse inclusion, we first remark that $s_i(\alpha_{i}^{\vee}) = - \alpha_{i}^{\vee}$, and that $s_1(\alpha_2^{\vee}) = s_2(\alpha_1^{\vee}) = \alpha_1^{\vee}+\alpha_2^{\vee}$. Note also that for $m, n \in \Z$, we have $s_{i,m}\circ s_{i,n} = \big( (m-n)\alpha_{i}^{\vee}, 1\big)$. One can then write
$$
(m\alpha_{1}^{\vee}+n\alpha_{2}^{\vee},s_1) = s_{2,n}\circ s_2 \circ s_{1,m}\,,
$$
and
$$
(m\alpha_{1}^{\vee}+n\alpha_{2}^{\vee},s_2\circ s_1) = s_2\circ s_{2,m}\circ s_{2,n}\circ s_{1,m}\,,
$$
and similar identities to conclude. The fact that the equilateral triangle with vertices  $O=(0,0)$, $A=(1,0)$, and $B = (\frac{1}{2},\frac{\sqrt{3}}{2})$ is a fundamental domain for $W_a$ follows from the fact that the sides of the triangle are supported by the lines $L_1, L_2$ and $L_{3,1}$.\\
(4) Clear.\\
(5) This assertion is illustrated by Figure~\ref{FGF-2}.  More precisely, we have the following correspondences of triangles which move the hexagon onto the fundamental domain $\{x\alpha_{1}^{\vee}+y\alpha_{1}^{\vee}~|~ 0 \le x,y \le 1\}$ for the action of $\Gamma$ on $\E^2$.

\begin{center}
\begin{tabular}{|c|c|c|}
  \hline
  Triangle & sent to & by \\\hline
  $[0,b,C]$ & $[g,h,A]$ & $\alpha_{1}^{\vee}$ \\
  $[0,C,c]$ & $[g,A,f]$ & $\alpha_{1}^{\vee}$ \\
  $[0,c,D]$ & $[i,j,B]$ & $\alpha_{3}^{\vee}$ \\
  $[0,D,d]$ & $[i,B,a]$ & $\alpha_{3}^{\vee}$ \\
  $[0,d,E]$ & $[i,a,A]$ & $\alpha_{3}^{\vee}$ \\
  $[0,E,e]$ & $[i,A,h]$ & $\alpha_{3}^{\vee}$ \\
  $[0,e,F]$ & $[k,b,B]$ & $\alpha_{2}^{\vee}$ \\
  $[0,F,f]$ & $[k,B,j]$ & $\alpha_{2}^{\vee}$ \\
  \hline
\end{tabular}
\end{center}
\hfill $\square$

\section{The equilateral torus}\label{S-ETO}

\subsection{Preliminaries}

The \emph{equilateral torus} is the flat torus $\T = \E^2/\Gamma$, where $\Gamma$  is the lattice $\Z \, \alpha_{1}^{\vee} \bigoplus \Z \, \alpha_{2}^{\vee}$.

Let $\Gamma^{*}$ be the dual lattice (see the notation $P(R)$ in \cite{Bou, Ber})
\begin{equation}\label{ETO-2}
\Gamma^{*} = \left\lbrace x \in \E^2 ~|~ \langle x,\gamma \rangle \in \Z, ~ \forall \gamma \in \Gamma \right\rbrace \,.
\end{equation}
The lattice $\Gamma^{*}$ admits the basis $\{\varpi_1,\varpi_2\}$, where
\begin{equation}\label{ETO-4}
\varpi_1 = (\frac{2}{3},0),\, \varpi_2 = (\frac{1}{3}, \frac{1}{\sqrt{3}}),
\end{equation}
which is dual to the basis $\{ \alpha_{1}^{\vee},\alpha_{2}^{\vee} \}$ of $\Gamma$.

Up to normalization, a complete set of eigenfunctions of the torus $\T$ is given by
\begin{equation}\label{ETO-6}
\phi_{p}(x) = \exp(2i\pi \langle x,p \rangle),
\end{equation}
where $p$ ranges over $\Gamma^{*}$,  and $x \in \E^2$\,. More precisely, for
$p \in \Gamma^{*}$, the function $\phi_{p}$ satisfies
\begin{equation*}
-\Delta \phi_{p} = 4\pi^2 |p|^2 \phi_{p}.
\end{equation*}

Writing $p = m \varpi_1 + n \varpi_2$, we find that the eigenvalues of the equilateral torus $\T$ are the numbers
\begin{equation}\label{ETO-8}
\hat{\lambda}(m,n) = \frac{16\pi^2}{9} (m^2+mn+n^2), \text{~for~} m,n \in \Z\,.
\end{equation}

The multiplicity of the eigenvalue $\hat{\lambda}(m,n)$ is the number
$$\# \left\{ (i,j) \in \Z^2 ~|~ i^2+ij+j^2 = m^2+mn+n^2 \right\}.$$

As usual, the counting function $N_{\T}(\lambda)$ of the torus $\T$ is defined by
\begin{equation}\label{ETO-10}
\begin{array}{lll}
N_{\T}(\lambda) & = &\# \left\{n \ge 1 ~|~ \lambda_n(\T) < \lambda \right\},\\[5pt]
& = &\# \left\{ (m,n) \in \Z^2 ~|~ \frac{16\pi^2}{9}(m^2+mn+n^2) < \lambda\right\}\,.
\end{array}
\end{equation}

\begin{lemma}\label{LETO-1}
The counting function of the equilateral torus satisfies
\begin{equation}\label{ETO-12}
N_{\T}(\lambda) \ge \frac{3\sqrt{3}}{2} \frac{\lambda}{4\pi} - \frac{9}{2\pi} \sqrt{\lambda} + 1\,.
\end{equation}
\end{lemma}

\textbf{Proof}. The lemma follows from easy counting arguments. Define the sets\footnote{The set $C$ is an open Weyl chamber \cite{Bou}.}
\begin{equation}\label{ETO-14}
\begin{array}{l}
C = \left\lbrace x\varpi_1+y\varpi_2 ~|~ x,y > 0\right\rbrace ,\\[5pt]
C_1 = \left\lbrace x\varpi_1+y\varpi_2 ~|~ x,y > 1\right\rbrace ,\\[5pt]
B(r) = \left\lbrace x \in \E^2 ~|~ |x| < r \right\rbrace .
\end{array}
\end{equation}

For $(i,j) \in \Z^2$, define the closed lozenge $\cL_{i,j}$,
\begin{equation}\label{ETO-16}
\cL_{i,j} = \left\lbrace x\varpi_1 + y\varpi_2 ~|~ i \le x \le i+1,\, j \le y \le j+1 \right\rbrace .
\end{equation}

Define the sets
\begin{equation}\label{ETO-18}
\begin{array}{l}
\cL(r) = \left\lbrace (i,j) \in \Z^2 ~|~ i\varpi_1+j\varpi_2 \in B(r) \right\rbrace\,,\\[5pt]
\cL_2(r) = \left\lbrace (i,j) \in \Z^2 ~|~ i\varpi_1+j\varpi_2 \in C\cap B(r) \right\rbrace\,,\\[5pt]
\cL_1(r) = \left\lbrace i\in \Z ~|~ i \ge 1 \text{~and~} i\varpi_1 \in B(r)\right\rbrace\,.\\
\end{array}
\end{equation}

Denote by $A(\Omega)$ the area of the set $\Omega$ contained in $\E^2$ or $\T$. Define $h_{\varpi}$ to be the height of the equilateral triangle $(0,\varpi_1, \varpi_2)$, \ie, $h_{\varpi} = \frac{1}{\sqrt{3}}$, and define $A_{\varpi}$ to be the area of the lozenges $\cL_{i,j}$.

\begin{figure}
  \centering
  \includegraphics[width=10cm]{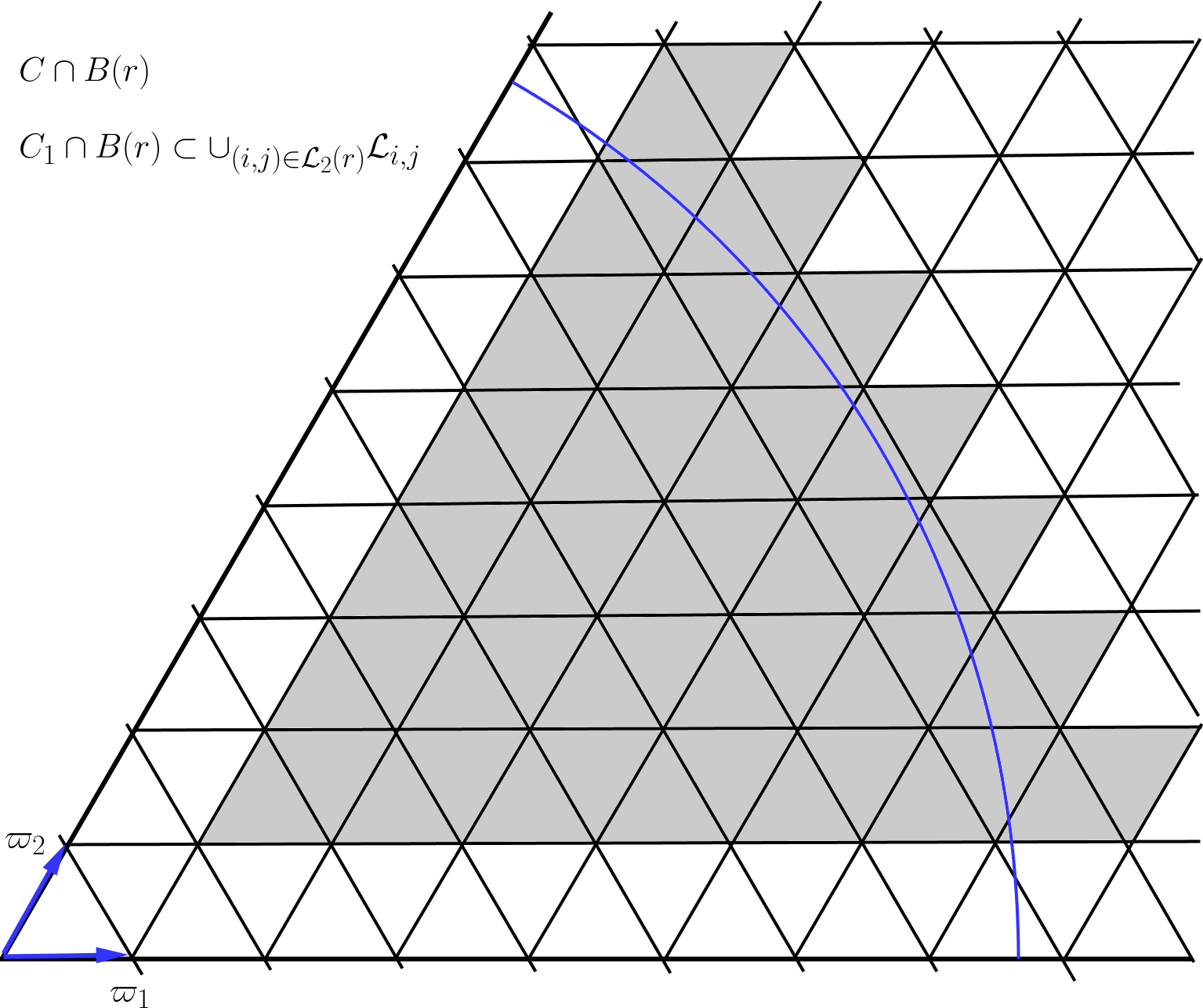}
  \caption{Illustration of \eqref{ETO-22}}\label{FETO-1}
\end{figure}

For symmetry reasons, we have the relation
\begin{equation}\label{ETO-20}
\# \cL(r) = 6 \, \# \cL_2(r) + 6 \, \# \cL_1(r) + 1 \,.
\end{equation}
Furthermore, see Figure~\ref{FETO-1}, it is easy to check, that
\begin{equation}\label{ETO-22}
C_1 \cap B(r) \subset \bigcup_{(i,j) \in \cL_2(r)} \cL_{i,j}.
\end{equation}
It follows that $A_{\varpi}\, \# \cL_2(r) \ge A(C_1 \cap B(r))$, and hence that $$A_{\varpi}\, \# \cL_2(r) \ge A(C\cap B(r)) - 2rh_{\varpi} + A_{\varpi}.$$

Using the fact that $h_{\varpi} = \frac{1}{\sqrt{3}}$ and $A_{\varpi} = \frac{2}{3\sqrt{3}}$, we obtain that
\begin{equation}\label{ETO-24}
\# \cL_2(r) \ge \frac{\sqrt{3}\pi}{4}r^2- 3r +1\,.
\end{equation}

Similarly, since $|\varpi_1| = \frac{2}{3}\,$, we have that
\begin{equation}\label{ET0-26}
\# \cL_1(r) = \left[ \frac{r}{|\varpi_1|} \right] \ge \frac{3}{2}r-1\,.
\end{equation}

Finally, we obtain
\begin{equation}\label{ETO-28}
\# \cL(r) \ge \frac{3\sqrt{3}\pi}{2} r^2 -9 r +1\,.
\end{equation}

Since $N_{\T}(\lambda) = \# \cL(\frac{\sqrt{\lambda}}{2\pi})$, we obtain the estimate
\begin{equation}\label{ETO-30}
N_{\T}(\lambda) \ge \frac{3\sqrt{3}}{2} \frac{\lambda}{4\pi} - \frac{9}{2\pi} \sqrt{\lambda} + 1\,.
\end{equation}
\hfill $\square$

\textbf{Remark}. Notice that the area of the equilateral torus is $$ A(\T) = \frac{3\sqrt{3}}{2}\,,$$
 so that the above lower bound is asymptotically sharp (Weyl's asymptotic law).

Denote by $a(\T)$ the length of the shortest closed geodesic on the flat torus $\T$ \ie,
\begin{equation}\label{lengthtore}
a(\T) = \frac{\sqrt{3}}{2}|\alpha_{1}^{\vee}| = \frac{3}{2}\,.
\end{equation}
 According to \cite[\S 7]{HoHuMo}, we  have the following isoperimetric inequality  in the particular case of the equilateral torus.

 \begin{lemma}
 If a domain $\Omega \subset \T$ has area $A(\Omega) \le \frac{(a(\T))^2}{\pi}$, then the length $\ell(\partial \Omega)$ of its boundary satisfies $$\ell^2(\partial \Omega) \ge 4\pi A(\Omega)\,.$$
 \end{lemma}

Together with \eqref{lengthtore}, this implies the

\begin{proposition}\label{LETO-2}
Let $\Omega \subset \T$ be a domain whose area satisfies $A(\Omega) \le \frac{9}{4\pi}$. Then, the first Dirichlet eigenvalue of $\Omega$ satisfies the inequality
\begin{equation}\label{ETO-32}
\lambda(\Omega) \ge \dfrac{\pi j^2_{0,1}}{A(\Omega)}\,,
\end{equation}
where $j_{0,1}$ is the first positive zero of the Bessel function of order $0$.\end{proposition}
We recall that $j_{0,1} \sim 2.4048255577$.

\textbf{Proof}. Apply the proof of the classical Faber-Krahn inequality given for example in \cite{FK}. \hfill $\square$

\subsection{Courant-sharp eigenvalues of the equilateral torus}

If $\lambda_n = \lambda_n(\T)$ is Courant-sharp, then $\lambda_{n-1}(\T) < \lambda_n(\T)$ and hence $N_{\T}(\lambda_n)=n-1$. In view of Lemma~\ref{LETO-1}, we obtain
$$
\frac{3\sqrt{3}}{8\pi}\lambda_n - \frac{9}{2\pi}\sqrt{\lambda_n} -(n-2) \le 0\,.
$$
It follows that
\begin{equation}\label{ETO-34}
\lambda_n(\T) \text{~Courant-sharp~} \Rightarrow \lambda_n(\T) \le 12 \left( 1 + \sqrt{1+\frac{2\pi}{9\sqrt{3}}(n-2)} \, \right)^2\,.
\end{equation}

If $\lambda_n$ is Courant-sharp, there exists an eigenfunction $u$ associated with $\lambda_n$ with exactly $n$ nodal domains. One of them, call it $\Omega$, satisfies $$A(\Omega) \le \frac{A(\T)}{n} = \frac{3\sqrt{3}}{2n}\,.$$
 If $n \ge 4$, $\Omega$ satisfies the assumption of Proposition~\ref{LETO-2} and hence,
\begin{equation*}
\lambda_n(\T) = \lambda(\Omega) \ge \frac{\pi j^2_{0,1}}{A(\Omega)}\,.
\end{equation*}
It follows that
\begin{equation}\label{ETO-36}
n \ge 4 \text{~and~} \lambda_n(\T) \text{~Courant-sharp~}  \Rightarrow
\lambda_n(\T) \ge \frac{2\pi j^2_{0,1}}{3\sqrt{3}}\, n\,.
\end{equation}

Comparing \eqref{ETO-34} and \eqref{ETO-36}, we see that if the eigenvalue  $\lambda_n(\T), \,n\ge 4\,$, is Courant-sharp, then $n\le 63\,$.  {\bf Hence it remains to examine condition \eqref{ETO-36} for the first $63$ eigenvalues.}

The following table gives the first 85 normalized eigenvalues $\bar{\lambda}_k$ such that $$\lambda_k = \frac{16\pi^2}{9} \bar{\lambda}_k\,.$$
The condition \eqref{ETO-36} becomes
\begin{equation}\label{ETO-36n}
n \ge 4 \text{~and~} \lambda_n(\T) \text{~Courant-sharp~}  \Rightarrow
\frac{\bar{\lambda}_n(\T)}{n} \ge \frac{\sqrt{3}j^2_{0,1}}{8\pi} \sim 0.3985546913 \,.
\end{equation}

The first column in the table displays the normalized eigenvalue $\bar{\lambda}$;  the second column the least integer $k$ such that $\bar{\lambda}_k = \bar{\lambda}\,$; the third column the largest integer $k$  such that $\bar{\lambda}_k = \bar{\lambda}$; the fourth column the multiplicity of $\bar{\lambda}$. The last column displays the ratio $\bar{\lambda}_k/k$ which should be larger than
$0.3985546913$ provided that $n\ge 4$ and $\lambda_n$ is Courant-sharp. The table proves Theorem~\ref{PETO-1}.

\begin{center}
\begin{tabular}{|c|c|c|c|c|}
  \hline
  eigenvalue& minimum index& maximum index& multiplicity& ratio \\ \hline
  0& 1& 1& 1 & \\ \hline
  1& 2& 7& 6 & \\ \hline
  3& 8& 13& 6 & 0.3750000000 \\ \hline
  4& 14& 19& 6 & 0.2857142857 \\ \hline
  7& 20& 31& 12 & 0.3500000000 \\ \hline
  9& 32& 37& 6 & 0.2812500000 \\ \hline
  12& 38& 43& 6 & 0.3157894737 \\ \hline
  13& 44& 55& 12 & 0.2954545455 \\ \hline
  16& 56& 61& 6 & 0.2857142857 \\ \hline
  19& 62& 73& 12 & 0.3064516129 \\ \hline
  21& 74& 85& 12 & 0.2837837838 \\
  \hline
\end{tabular}
\end{center}

\textbf{Remark}. The ratio $\bar{\lambda}_n(\T)/n$ is meaningful to study Courant-sharpness for $n\ge 4$ only. This is why this information is not calculated in the first two lines.

\section{ The equilateral triangle: spectrum and action of symmetries}\label{S-ETR}

In this section, we start the analysis of the case of the equilateral triangle by recalling its spectrum and exploring the action of symmetries in each eigenspace. We keep the notation of the previous sections.

\subsection{Eigenvalues and eigenfunctions}

Recall the following result \cite[Proposition~9]{Ber}.

\begin{proposition}\label{PETR-1}
Up to normalization, a complete set of eigenfunctions of the Dirichlet Laplacian in the equilateral triangle $\cT = \{0,A,B\}$,  with sides of length $1$, is given by the functions
\begin{equation}\label{ETR-2}
\begin{array}{lll}
\Phi_{p} &=& \sum_{w \in W} \epsilon(w) \phi_{w(p)} \,,\\[5pt]
\Phi_{p}(x) &=& \sum_{w \in W} \epsilon(w) \exp\left( 2i\pi \langle x,w(p)\rangle\right) \,,
\end{array}
\end{equation}
where $p$ ranges over the set $C \cap \Gamma^{*}$, and where $\epsilon(w)$ is the determinant of $w$. The associated eigenvalues are the numbers $4\pi^2|p|^2$ for $p \in C \cap \Gamma^{*}$. The multiplicity of the eigenvalue $4\pi^2|p|^2$ is given by $\# \{q \in C \cap \Gamma^{*} ~|~ |q|=|p|\}$.
\end{proposition}%

\textbf{Remark}. Notice that
$$
C \cap \Gamma^{*} = \left\lbrace m\varpi_1 + n\varpi_2 ~|~ m,n \in \Nb \right\rbrace \,,
$$
so that the eigenvalues of the equilateral triangle $\cT$, with sides of length $1$, are the numbers %
$$
\frac{16\pi^2}{9}(m^2+mn+n^2)\, \text{~for~} m, n \in \Nb\,.
$$

\textbf{Idea of the proof}.\\
We follow   \cite{Ber}. Given a Dirichlet eigenfunction $\phi$ of the triangle, we extend it to a function $\psi$ on $\E^2$ using the symmetries $s_{i,k}$, in such a way that $\psi\circ w = \epsilon(w) \psi$ for any $w$ in the group $W_a$. This is possible because $\cT$ is a fundamental domain for $W_a$. The function $\psi$ turns out to be smooth. Because $W_a = \Gamma \rtimes W$, the function $\psi$ is $\Gamma$-periodic, and hence defines an eigenfunction $\Phi$ on the torus $\E^2/\Gamma$ which satisfies $\Phi \circ w = \epsilon(w) \Phi$ for all $w \in W$. Conversely, any eigenfunction $\Phi$ on the torus, which satisfies this condition, gives a Dirichlet eigenfunction of the triangle. It remains to identify the eigenfunctions of the torus which satisfy the condition. This is done by making use of \cite[\S VI.3, Proposition~1]{Bou}.
\hfill $\square$ \medskip

The Dirichlet eigenfunctions of $\cT$ look a little simpler in the following parametrization $\cF$ of $\E^2$,
\begin{equation}\label{ETR-4}
\begin{array}{lll}
\cF &:& \R^2 \to \E^2 \,,\\[5pt]
\cF &:& (s,t) \mapsto s \alpha_{1}^{\vee} + t \alpha_{2}^{\vee} \,.
\end{array}%
\end{equation}

Given a point $p\in \E^2$, we will denote by $(x,y)$ its coordinates with respect to the standard basis $\{e_1,e_2\}$, and by $(s,t)_{\cF}$ its coordinates in the parametrization $\cF$. More precisely, we will write,
\begin{equation}\label{ETR-4a}
\left\{
\begin{array}{ll}
p = (x,y), & \text{~for~} p = xe_1+ye_2,\\
p = (s,t)_{\cF} & \text{~for~} p = s\alpha_1^{\vee} + t \alpha_2^{\vee}.
\end{array}
\right.
\end{equation}

\begin{figure}
  \centering
  \includegraphics[width=6cm]{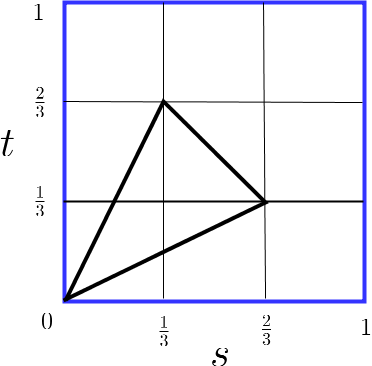}
  \caption{Parametrization $\cF$}\label{FETR-1}
\end{figure}

In the parametrization $\cF$, a fundamental domain for the action of $\Gamma$ on $\R^2$ is the square $\{0\le s \le 1\}\times\{0\le t \le 1\}$; a fundamental domain for the action of $W_a$ on $\R^2$ is the triangle with vertices $(0,0)_{\cF}$, $(\frac{2}{3},\frac{1}{3})_{\cF}$ and $(\frac{1}{3},\frac{2}{3})_{\cF}$, see Figure~\ref{FETR-1}.

\medskip Notice that the parametrization $\cF$ is not orthogonal, and that the Laplacian is given in this parametrization by
\begin{equation}\label{ETR-6}
\Delta_{\cF} = \frac{4}{9}\left( \partial^2_{ss} + \partial^2_{st} + \partial^2_{tt} \right)\,.
\end{equation}
 Recall that the function $\phi_{p}$ is defined by $\phi_{p}(x) = \exp\left( 2i\pi \langle x,p \rangle\right)$, for $x \in \E^2$ and $p \in \Gamma$. In the parametrization $\cF$, writing $x=s\alpha_{1}^{\vee}+t\alpha_{2}^{\vee}$ and $p = m\varpi_1+n\varpi_2$, the function  $\phi_p$ will be written as
\begin{equation}\label{ETR-6a}
\phi_{m,n}(s,t) = \exp\left( 2i\pi (ms+nt) \right).
\end{equation}

To write the eigenfunctions of the Dirichlet Laplacian in the equilateral triangle, we have to compute the scalar products $\langle x, w(p) \rangle$ or equivalently $\langle w(x), p \rangle$ for $x = s\alpha_{1}^{\vee}+t\alpha_{2}^{\vee}$ and $p = m\varpi_1+n\varpi_2$. Table~\ref{TETR-1} displays the result.\medskip

\begin{table}[h!t]\label{TETR-1}
\begin{center}
\begin{tabular}{|c|r|c|}
  \hline
  $w$ & $\det(w)$ & $\langle w(x),p \rangle$ \\[5pt]\hline
  $1$ & $1$ & $ms+nt$ \\[5pt]\hline
  $s_1$ & $-1$ & $-ms+(m+n)t$ \\[5pt]\hline
  $s_2$ & $-1$ & $(m+n)s-nt$ \\[5pt]\hline
  $s_3$ & $-1$ & $-ns-mt$ \\[5pt]\hline
  $s_1\circ s_2$ & $1$ & $ns-(m+n)t$ \\[5pt]\hline
  $s_2\circ s_1$ & $1$ & $-(m+n)s+mt$ \\
  \hline
\end{tabular}
\caption{Action of $W$}
\end{center}%
\end{table}%

 Define the functions
\begin{equation}\label{ETR-8E}
\begin{split}
E_{m,n}(s,t) =
& \exp\left(2i\pi\big( ms+nt\big)\right) - \exp\left(2i\pi\big( -ms+(m+n)t\big)\right) \\
& - \exp\left(2i\pi\big( (m+n)s-nt\big)\right) - \exp\left(2i\pi\big( -ns-mt\big)\right) \\
& + \exp\left(2i\pi\big( ns-(m+n)t\big)\right) + \exp\left(2i\pi\big( -(m+n)s+mt\big)\right) \,,
\end{split}
\end{equation}

 Using Table~\ref{TETR-1}, we see that
\begin{equation}\label{symm}
E_{n,m}(s,t) = - \overline{E_{m,n}(s,t)}\,.
\end{equation}
Looking at the pairs $[m,n]$ and $[n,m]$ simultaneously, and making use of real eigenfunctions instead of complex ones, we find that eigenfunctions associated with the pairs of positive integers $[m,n]$ and $[n,m]$, are $C_{m,n}(s,t)$ and $S_{m,n}(s,t)$, given by the following formulas,
\begin{equation}\label{ETR-8a}
\begin{split}
C_{m,n}(s,t) =
& \cos2\pi\big( ms+nt\big) - \cos2\pi\big( -ms+(m+n)t\big)\\
& - \cos2\pi\big( (m+n)s-nt\big)- \cos2\pi\big( -ns-mt\big)\\
& + \cos2\pi\big( ns-(m+n)t\big) + \cos2\pi\big( -(m+n)s+mt\big)\,,
\end{split}
\end{equation}
and
\begin{equation}\label{ETR-8b}
\begin{split}
S_{m,n}(s,t) =
& \sin2\pi\big( ms+nt\big) - \sin2\pi\big( -ms+(m+n)t\big)\\
& - \sin2\pi\big( (m+n)s-nt\big) - \sin2\pi\big( -ns-mt\big)\\
& + \sin2\pi\big( ns-(m+n)t\big) + \sin2\pi\big( -(m+n)s+mt\big)\,.
\end{split}
\end{equation}

\textbf{Remarks}. Notice that by \eqref{symm} $C_{m,n} = - C_{n,m}$ and $S_{m,n}=S_{n,m}$. The line $\{s=t\}$ corresponds to the median of the equilateral triangle $\cT$ issued from $O$. This median divides $\cT$ into two congruent triangles with angles $\{30,60,90\}$ degrees. The eigenfunctions $C_{m,n}$ correspond to the Dirichlet eigenfunctions of these triangles. When $m=n$, we have $C_{m,m} \equiv 0$ and
\begin{equation*}
S_{m,m}(s,t) =2 \left\{ \sin2\pi m(s+t)\ - \sin2\pi m (2t-s) - \sin2\pi m (2s-t)\right\}\,,
\end{equation*}
\ie,
\begin{equation}\label{ETR-8S}
 S_{m,m}(s,t) = S_{1,1}(ms,mt)\,.
\end{equation}

\subsection{Symmetries}\label{S-sym}

Call $F_C$ the centroid (or center of mass) of the equilateral triangle $\cT = \{0,A,B\}$.  With the conventions \eqref{ETR-4a}, $F_C=(\frac{1}{2},\frac{\sqrt{3}}{6})=(\frac{1}{3},\frac{1}{3})_{\cF}$. The isometry group $G_{\cT}$ of the triangle $\cT$ has six elements: the three orthogonal symmetries with respect to the medians of the triangle (which fix one vertex and exchange the two other vertices), the rotations $\rho_{\pm}$, with center $F_C$ and angles $\pm\frac{2\pi}{3}$ (which permute the vertices), and the identity. The group $G_{\cT}$ fixes the centroid $F_C$.

 In the parametrization $\cF$, the symmetries are given by
\begin{equation}\label{sym-2}
\begin{array}{lll}
\sigma_1 (s,t) & = & (t,s)\,,\\[5pt]
\sigma_2 (s,t) & = & (-s+\frac{2}{3},t-s+\frac{1}{3})\,,\\[5pt]
\sigma_3 (s,t) & = & (s-t+\frac{1}{3},-t+\frac{2}{3})\,.
\end{array}
\end{equation}
\ie, respectively, the symmetries with respect to the median issued from $0$, to the median issued from $B$, and to the median issued from $A$.

As a matter of fact, the group $G_{\cT}$ is generated by the symmetries $\sigma_1$ and $\sigma_2$. The rotations are
$$ \rho_{+} = \sigma_2\circ \sigma_1 : (s,t) \mapsto (-t+\frac{2}{3},s-t+\frac{1}{3}),$$
with angle $\frac{2\pi}{3}$, and
$$ \rho_{-} = \sigma_1\circ \sigma_2 : (s,t) \mapsto (t-s+\frac{1}{3},-s+\frac{2}{3}),$$
with angle $ -\frac{2\pi}{3}$. \\
Furthermore, $\sigma_3 = \sigma_1 \circ \sigma_2\circ \sigma_1 = \sigma_2 \circ \sigma_1\circ \sigma_2\,$.\medskip

The action of the  symmetries $\sigma_1, \sigma_2$ on the eigenfunctions $C_{m,n}$ and $S_{m,n}$ is given by
\begin{equation}\label{sym-4}
\begin{array}{lll}
\sigma_1^{*}C_{m,n} &=& - C_{m,n}\,,\\[5pt]
\sigma_1^{*}S_{m,n} &=& S_{m,n}\,,\\[5pt]
\sigma_2^{*}C_{m,n} &=& - \cos\alpha_{m,n} \, C_{m,n} - \sin\alpha_{m,n} \, S_{m,n}\,,\\[5pt]
\sigma_2^{*} S_{m,n}  &=& - \sin\alpha_{m,n} \, C_{m,n} + \cos\alpha_{m,n} \, S_{m,n}\,,
\end{array}
\end{equation}
where $\alpha_{m,n} = \frac{2\pi(2m+n)}{3}\,$.

We work in the parametrization $\cF$, and use the convention \eqref{ETR-4a}.  For the pairs $(m,n)$ and $(n,m)$, we define the family of eigenfunctions
\begin{equation}\label{sym-6}
\Psi_{m,n}^{\theta}(s,t) = \cos \theta\,  C_{m,n}(s,t) + \sin \theta \, S_{m,n} (s,t)\,,
\end{equation}
  associated with the eigenvalue $\frac{16 \pi^2}{9} (m^2 + mn + n^2)$.

The action of the symmetries on the family  $\Psi_{m,n}^{\theta}$ defined by \eqref{sym-6} is given by
\begin{equation}\label{sym-8}
 \begin{array}{lll}
\sigma_1^{*}\Psi_{m,n}^{\theta} &=& \Psi_{m,n}^{\pi-\theta}\,,\\[5pt]
\sigma_2^{*}\Psi_{m,n}^{\theta} &=& \Psi_{m,n}^{\pi+\alpha_{m,n}-\theta}\,,\\[5pt]
\sigma_3^{*}\Psi_{m,n}^{\theta} &=& \Psi_{m,n}^{\pi-\alpha_{m,n}-\theta}\,.\\[5pt]
\end{array}
\end{equation}


\section{Pleijel's approach for Courant-sharp eigenvalues}\label{S5}

In order to investigate the Courant-sharp Dirichlet eigenvalues of the equilateral triangle $\cT$, we use the same methods as  in \cite{Pl} and Section~\ref{S-ETO}.

Notice that the counting function of the Dirichlet eigenvalues of the equilateral triangle is given by
\begin{equation}\label{ETR-10}
N_{\cT}(\lambda) = \cL_{2}(\frac{\sqrt{\lambda}}{2\pi}) = \# \left\lbrace
(k,\ell) \in \Nb \times\Nb ~|~ |k\varpi_1+\ell\varpi_2|^2 < \frac{\lambda}{4\pi^2}\right\rbrace \,.
\end{equation}

Using \eqref{ETO-24}, we conclude that
\begin{equation}\label{ETR-12}
N_{\cT}(\lambda) \ge \frac{\sqrt{3}}{4} \frac{\lambda}{4\pi} - \frac{3}{2\pi} \sqrt{\lambda} + 1\,.
\end{equation}
Notice that this lower bound is asymptotically sharp (Weyl's asymptotic law).

Assuming  that $\lambda_n(\cT)$ is Courant-sharp, we have
$\lambda_{n-1}(\cT) < \lambda_{n}(\cT)$, and hence\goodbreak $N_{\cT}(\lambda_n(\cT)) = n-1$. It follows that
\begin{equation}\label{ETR-14}
\lambda_n(\cT) \text{~Courant-sharp~} \Rightarrow \lambda_n(\cT) \le 48 \left( 1 + \sqrt{1 + \frac{\pi}{3\sqrt{3}}(n-2)} \right)^2\,.
\end{equation}

On the other-hand, if $\lambda_n(\cT)$ is Courant-sharp, there exists an eigenfunction $u$ with exactly $n$ nodal domains $\Omega_1, \ldots \Omega_n$, for which we can write
\begin{equation*}
\lambda_n(\cT) = \lambda(\Omega_i) \ge \frac{\pi j^2_{0,1}}{A(\Omega_i)}\,,
\end{equation*}
where we have used the Faber-Krahn inequality in $\E^2$.\\
 Summing up in $i$, it follows that
\begin{equation}\label{ETR-18}
\lambda_n(\cT) \text{~Courant-sharp~} \Rightarrow \frac{\lambda_n(\cT)}{n} \ge \frac{\pi j^2_{0,1}}{A(\cT)} = \frac{4\pi j^2_{0,1}}{\sqrt{3}}\,.
\end{equation}

Combining \eqref{ETR-14} and \eqref{ETR-18}, we find that if $\lambda_n(\cT)$ is Courant-sharp, then $n \le  40 \,$. It follows that to determine the Courant-sharp eigenvalues, it suffices to look at the first $40$ eigenvalues of the equilateral triangle. Using \eqref{ETR-18} again, we compute the ratios $\frac{\bar{\lambda}_n(\cT)}{n}$ for the first $40$ normalized eigenvalues ($\bar{\lambda}_n(\cT) = \frac{9}{16\pi^2}\lambda_n(\cT)$),  and compare them with the value $$\frac{3\sqrt{3}j^2_{0,1}}{4\pi} \sim 2.391328148\,.$$
 This is given by the following table in which the first column gives the normalized eigenvalue  $\bar{\lambda}$; the second column the smallest index $i$ such that $\bar{\lambda}_i = \bar{\lambda}$; the third column the largest index $j$ such that $\bar{\lambda}_j = \bar{\lambda}$; the fourth column the multiplicity of $\bar{\lambda}$; and the last one the ratio normalized eigenvalue/least index (which is the relevant information for checking the Courant-sharp property).

\begin{table}[ht]
\begin{center}
\begin{tabular}{|c|c|c|c|c|}
  \hline
&&&&\\
  $\bar{\lambda}$& $\bar{\lambda}_i=\bar{\lambda}$& $\bar{\lambda}_j=\bar{\lambda}$& $\mathrm{mult}(\bar{\lambda})$& $\frac{\bar{\lambda}_i(\cT)}{i}$ \\[5pt] \hline
  3& 1& 1& 1 & \textbf{3} \\ \hline
  7& 2& 3& 2 & \textbf{3.5}\\ \hline
  12& 4& 4& 1 & \textbf{3} \\ \hline
  13& 5& 6& 2 & \textbf{2.6000000} \\ \hline
  19& 7& 8& 2 & \textbf{2.7142857} \\ \hline
  21& 9& 10& 2 & 2.333333333 \\ \hline
  27& 11& 11& 1 & \textbf{2.45454545} \\ \hline
  28& 12& 13& 2 & 2.333333333 \\ \hline
  31& 14& 15& 2 & 2.214285714 \\ \hline
  37& 16& 17& 2 & 2.312500000 \\ \hline
  39& 18& 19& 2 & 2.166666667 \\ \hline
  43& 20& 21& 2 & 2.150000000 \\ \hline
  48& 22& 22& 1 & 2.181818182 \\ \hline
  49& 23& 24& 2 & 2.130434783 \\ \hline
  52& 25& 26& 2 & 2.080000000 \\ \hline
  57& 27& 28& 2 & 2.111111111 \\ \hline
  61& 29& 30& 2 & 2.103448276 \\ \hline
  63& 31& 32& 2 & 2.032258065 \\ \hline
  67& 33& 34& 2 & 2.030303030 \\ \hline
  73& 35& 36& 2 & 2.085714286 \\ \hline
  75& 37& 37& 1 & 2.027027027 \\ \hline
  76& 38& 39& 2 & 2. \\ \hline
  79& 40& 41& 2 & 1.975000000 \\ \hline
\end{tabular}
\end{center}
\caption{Courant-sharp eigenvalues satisfy $\frac{\bar{\lambda}_n(\cT)}{n}\ge 2.391328148$}\label{ETR-TCS}
\end{table}

\begin{lemma}\label{ETR-L1}
The only possible Courant-sharp eigenvalues of the equilateral triangle are the $\lambda_k(\cT)$, for $k \in \{1,2,4,5,7,11\}$.
\end{lemma}%

\begin{figure}[ht]
  \centering
  \includegraphics[width=10cm]{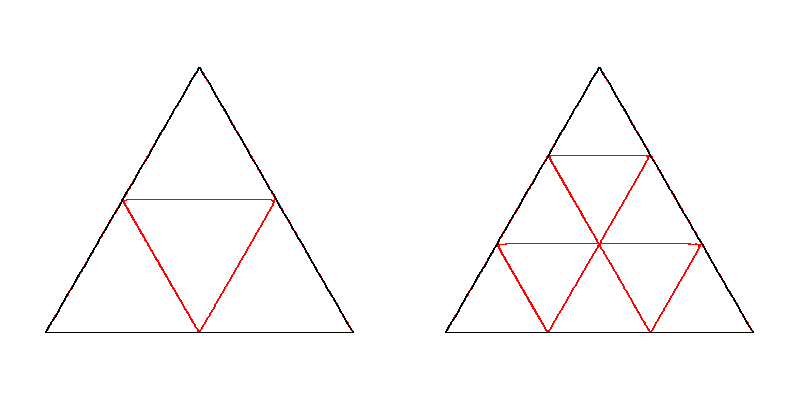}
  \caption{Nodal sets corresponding to $\lambda_4$ and $\lambda_{11}$}\label{ETF-F3}
\end{figure}

Clearly, the eigenvalues $\lambda_1$ and $\lambda_2$ are Courant-sharp. The eigenvalues $\lambda_4 = \hat{\lambda}(2,2)$ and $\lambda_{11} = \hat{\lambda}(3,3)$ are simple. It is easy to see that the number of nodal domains is $4$, \resp $9$, for these eigenvalues, see \eqref{ETR-8S} and Figure~\ref{ETF-F3}. Hence $\lambda_4$ is Courant sharp, and $\lambda_{11}$ is not Courant sharp.

\begin{figure}[ht]
  \centering
  \includegraphics[width=10cm]{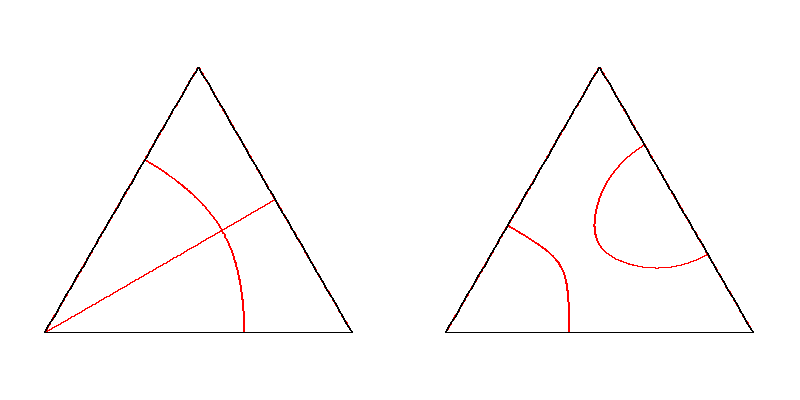}
  \caption{Nodal sets  of the eigenfunctions $C_{1,3}$ and $S_{1,3}$ corresponding to $\lambda_5$}\label{ETR-F4}
\end{figure}

\begin{figure}[ht]
  \centering
  \includegraphics[width=10cm]{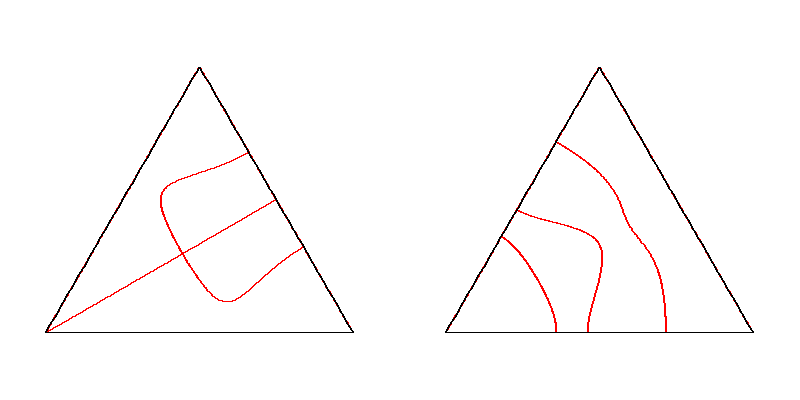}
  \caption{Nodal sets  of  the eigenfunctions $C_{2,3}$ and $S_{2,3}$ corresponding to $\lambda_7$}\label{ETR-F5}
\end{figure}

\medskip In order to determine the Courant-sharp eigenvalues, it therefore remains to consider the eigenvalues $\lambda_5 = \hat{\lambda}(1,3) = \hat{\lambda}(3,1)$ and $\lambda_{7} = \hat{\lambda}(2,3) = \hat{\lambda}(3,2)$ which have multiplicity $2\,$. We study these eigenvalues in Sections~\ref{SS-E13} and \ref{SS-E23} respectively, see Figures~\ref{ETR-F4} and \ref{ETR-F5}.

\begin{figure}[ht]
  \begin{minipage}[b]{0.45\linewidth}
    \centering
    \includegraphics[width=1\linewidth]{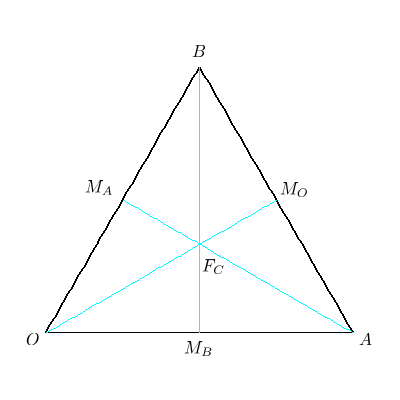}
    \par\vspace{-25pt}
  \end{minipage}%
  \begin{minipage}[b]{0.45\linewidth}
    \centering
\begin{tabular}[b]{|c|c|c|}
\hline
Point & $\E^2$ & $\cF$\\[5pt] \hline
Vertex $O$ & $(0,0)$ & $(0,0)_{\cF}$\\[5pt] \hline
Vertex $A$ & $(1,0)$ & $(\frac{2}{3},\frac{1}{3})_{\cF}$\\[5pt] \hline
Vertex $B$ & $(\frac{1}{2},\frac{\sqrt{3}}{2})$ & $(\frac{1}{3},\frac{2}{3})_{\cF}$\\[5pt] \hline
Centroid $F_C$ & $(\frac{1}{2},\frac{\sqrt{3}}{6})$ & $(\frac{1}{3},\frac{1}{3})_{\cF}$\\[5pt] \hline
Mid-point $M_O$ & $(\frac{3}{2},\frac{\sqrt{3}}{4})$ & $(\frac{1}{2},\frac{1}{2})_{\cF}$\\[5pt] \hline
Mid-point $M_A$ & $(\frac{1}{2},\frac{\sqrt{3}}{4})$ & $(\frac{1}{6},\frac{1}{3})_{\cF}$\\[5pt] \hline
Mid-point $M_B$ & $(\frac{1}{2},0)$ & $(\frac{1}{3},\frac{1}{6})_{\cF}$\\[5pt] \hline
\end{tabular}
    \par\vspace{0pt}
    \end{minipage}
\caption{Triangle $\cT$} \label{RC-F1}
\end{figure}

\medskip We will work in the triangle $\cT$, with vertices $\{O,A,B\}$. We denote by $F_C$ the centroid of the triangle. The median issued from the vertex $O$ is denoted by $[OM]$, the mid-point of the side $BA$ by $M_O\,$. We use similar notation for the other vertices, see Figure~\ref{RC-F1}\,.

\textbf{Remark}. We denote by $N(\Psi_{m,n}^{\theta})$ the nodal set of the eigenfunction $\Psi_{m,n}^{\theta}$ \ie, the closure of the set of zeros of the eigenfunction in the interior of the triangle. The eigenfunctions also vanish on the edges of the triangles, and we will analyze separately the points in the open edges and the vertices. For a summary of the general properties of nodal sets of eigenfunctions, we refer to \cite[Section~5]{BeHe1}. Observe that the eigenfunctions $\Psi_{m,n}^{\theta}$ are defined over the whole plane, so that even at the boundary of the triangle, we can use the local structure of the nodal set. \medskip

The remaining part of this paper is devoted to the proof of Theorem~\ref{PETR-2}, using the following strategy.

\section{ Playing on the checkerboard, a strategy to\\ determine nodal patterns}\label{sProl}
After the reduction \`a la Pleijel, we now have to  determine the nodal pattern of a general eigenfunction in two $2$-dimensional eigenspaces of the Laplacian in $\cT$, the interior of the equilateral triangle.  Recall that the nodal set $N(\Psi)$ of a Dirichlet eigenfunction $\Psi$ is the closure of the set $\{x \in \cT ~|~ \Psi(x)=0\}$ in $\overline{\cT}$. The function $\Psi$ actually extends smoothly to the whole plane, and the nodal set of $\Psi$ consists of finitely many regular arcs which may intersect, or hit the boundary of $\cT$, with the equiangular property, \cite[Section~2.1]{HOMN}.

\begin{figure}
  \centering
  \includegraphics[width=14cm]{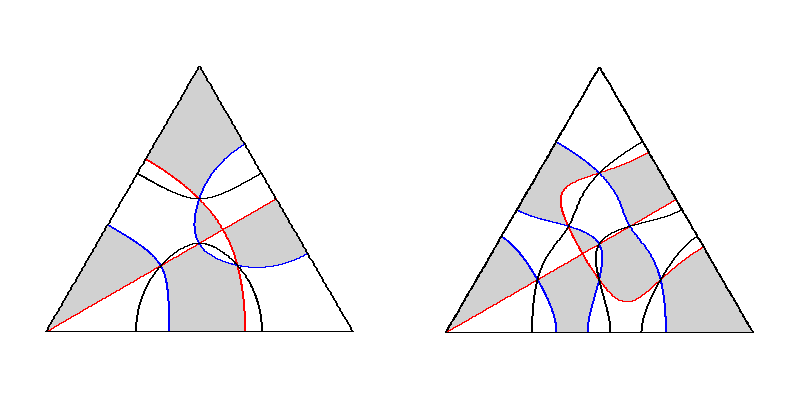}\\
  \caption{Checkerboards for $\lambda_5$ and $\lambda_7$, and nodal set $N(\Psi^{\frac{\pi}{6}})$}\label{I-F1}
\end{figure}

Let $\lambda$ be either $\lambda_5$ or $\lambda_7$, and let $\cE$ be the associated ($2$-dimensional) eigenspace. Determining whether $\lambda$ is Courant-sharp amounts to describing the possible nodal sets of the family $\Psi^{\theta} = \cos\theta \, C + \sin\theta \, S$, where $\{C,S\}$ is a basis of $\cE$, and $\theta \in [0,2\pi]$. For this purpose, we use the following ideas.
\vspace{-3mm}
\begin{enumerate}
  \item Using the natural \emph{symmetries} of the triangle, we restrict the analysis to $\theta \in [0,\frac{\pi}{6}]$.
  \item The set of \emph{fixed points} $N(C) \cap N(S)$ \ie, the set of common zeros of the functions $\Psi^{\theta}$, plays an important role in the analysis of the family of nodal sets $N(\Psi)$. When the variables are separated, as in the case of the square membrane, it is easy to determine the fixed points, \cite{St,BeHe1}. In the case of the equilateral triangle, we shall only determine the fixed points located on the medians of $\cT$.
  \item If $\theta \in ]0,\frac{\pi}{6}]$, the nodal set $N(\Psi^{\theta})$ is contained in the set   $\{C\, S < 0\} \cup \left( N(C) \cap N(S)\right)$, \cite{St,BeHe1}. This \emph{checkerboard argument} is illustrated in Figure~\ref{I-F1} which displays the checkerboards for $\lambda_5$ and $\lambda_7$, with   $C = \Psi^{0}$ and $S = \Psi^{\frac{\pi}{2}}$, as well as the nodal set of $\Psi^{\frac{\pi}{6}}$ (Maple simulations). For the equilateral triangle, the variables are not separated, and we shall not use this argument directly, but rather \emph{separation lemmas} involving the medians and other natural lines.
  \item Critical zeros \ie, points at which both the eigenfunction and its first  partial derivatives vanish, play a central role in our analysis. They are easy to determine when the variables are separated, \cite{BeHe1}. In the case of the equilateral triangle, the determination of critical zeros is more involved. We first analytically determine the critical zeros on the boundary \ie, the points at which the nodal sets hits $\partial \cT$. Later on in the proof, we show that $\Psi^{\theta}$ does not have interior critical zeros when $\theta \in ]0,\frac{\pi}{6}]$. For this purpose, we use the following \emph{energy argument}.
  \item \emph{Energy argument}. The medians divide the triangle $\cT$ into six isometric  hemiequilateral triangles, whose first Dirichlet eigenvalue is easily seen to be strictly larger than $\lambda_7$. This argument is needed to discard possible interior critical zeros, and also simply closed nodal arcs.
\end{enumerate} \vspace{-3mm}
For a given $\theta \in ]0,\frac{\pi}{6}]$, the above information provides some points in $N(\Psi^{\theta})$. It remains to determine how the nodal arcs can join these points, without crossing the established barriers.

\section{The eigenvalue $\lambda_5(\cT)$ and its eigenspace}\label{SS-E13}

We call $\cE_{5}$ the $2$-dimensional eigenspace associated with the eigenvalue $\lambda_5(\cT)$, \ie, with the pairs $[1,3]$ and $[3,1]$. Recall the eigenfunctions,
\begin{equation}\label{S13-0}
\begin{array}{ll}
C_{1,3}(s,t) = & \cos 2\pi(s+3t) - \cos 2\pi(-s+4t) - \cos 2\pi(4s-3t)\\
& - \cos 2\pi(-3s-t) + \cos 2\pi(3s-4t) + \cos 2\pi(-4s+t)\,,\\[5pt]
S_{1,3}(s,t) = & \sin2\pi(s+3t) - \sin2\pi(-s+4t) - \sin 2\pi(4s-3t)\\
& - \sin 2\pi(-3s-t) + \sin 2\pi(3s-4t) + \sin 2\pi(-4s+t)\,.
\end{array}
\end{equation}

In this section, we study   more carefully the nodal sets of the family $\Psi_{1,3}^{\theta} = \cos\theta \, C_{1,3} + \sin\theta \, S_{1,3}\,,$ following the strategy sketched in Section~\ref{sProl}. The first five subsections contain preliminary results. In Subsection~\ref{S13-NCS}, we determine rigorously the nodal sets of the functions $C_{1,3}$ and $S_{1,3}$. In Subsection~\ref{S13-NSP}, we determine the nodal sets of the function $\Psi^{\theta}_{1,3}$ when $\theta \in ]0,\frac{\pi}{6}]$. As a consequence, we obtain that the number of nodal domains is at most $4$, so that the eigenvalue $\lambda_5(\cT)$ is not Courant-sharp.

\subsection{Symmetries}\label{S13-sym}

Taking Subsection~\ref{S-sym} into account, we find that
\begin{equation}\label{E13-sym2}
\begin{array}{lll}
\sigma_1^{*}\Psi_{1,3}^{\theta} &=& \Psi_{1,3}^{\pi-\theta}\,,\\[5pt]
\sigma_2^{*}\Psi_{1,3}^{\theta} &=& \Psi_{1,3}^{\frac{\pi}{3}-\theta}\,,\\[5pt]
\sigma_3^{*}\Psi_{1,3}^{\theta} &
=& \Psi_{1,3}^{\frac{5\pi}{3}-\theta} \,,\\[5pt]
(\sigma_2 \circ \sigma_1)^{*}\Psi_{1,3}^{\theta} &=& \Psi_{1,3}^{\frac{2\pi}{3}+\theta}\,,\\[5pt]
(\sigma_1 \circ \sigma_2)^{*}\Psi_{1,3}^{\theta} &=& \Psi_{1,3}^{\frac{4\pi}{3}+\theta}\,.
\end{array}
\end{equation}

It follows that, up to multiplication by a scalar, and for $V \in \{O,A,B\}\,$, the eigenspace $\cE_5$ contains a unique eigenfunction $S_V$, \resp a unique eigenfunction $C_V$, which is invariant, \resp anti-invariant, under the symmetry $\sigma_{i(V)}$ with respect to the median $[VM]$ issued from the vertex $V$, where $i(O)=1, i(A)=3\,$, and $i(B)=2\,$. More precisely, we  define
\begin{equation}\label{E13-sym4}
\begin{array}{l}
S_O = \Psi_{1,3}^{\frac{\pi}{2}} \text{~and~} C_O = \Psi_{1,3}^{0}\,,\\[5pt]
S_A = \Psi_{1,3}^{\frac{11\pi}{6}} \text{~and~} C_A = \Psi_{1,3}^{\frac{4\pi}{3}}\,,\\[5pt]
S_B = \Psi_{1,3}^{\frac{7\pi}{6}} \text{~and~} C_B = \Psi_{1,3}^{\frac{2\pi}{3}}\,.
\end{array}
\end{equation}

Note that $C_O$, \resp $S_O$, are the functions $C_{1,3}\,$, \resp $S_{1,3}\,$, and that,
\begin{equation}\label{E13-sym4a}
\begin{array}{l}
S_A = (\sigma_1 \circ \sigma_2)^{*} S_O \text{~and~} C_A = (\sigma_1 \circ \sigma_2)^{*} C_O\,,\\[5pt]
S_B = (\sigma_2 \circ \sigma_1)^{*} S_O \text{~and~} C_B = (\sigma_2 \circ \sigma_1)^{*} C_O\,.
\end{array}%
\end{equation}

The eigenfunctions $S_V\,$, \resp $C_V\,$, are permuted under the action of the rotations   $\rho_{+}=\sigma_2 \circ \sigma_1$ and $\rho_{-}=\sigma_1 \circ \sigma_2\,$. The eigenspace $\cE_5$ does not contain any non trivial rotation invariant eigenfunction.\medskip

Since $\Psi_{1,3}^{\theta+\pi} = - \Psi_{1,3}^{\theta}\,$, it follows from \eqref{E13-sym2} that, up to the symmetries $\sigma_i$, the nodal sets of the family $\Psi_{1,3}^{\theta}\,$, $\theta \in [0,2\pi]$ are determined by the nodal sets of the sub-family $\theta \in [0,\frac{\pi}{6}]\,$.\medskip

\textbf{From now on, we assume that $\theta \in [0,\frac{\pi}{6}]$.}

\subsection{Behaviour at the vertices}\label{S13-vert}

The vertices of the equilateral triangle $\cT$ belong to the nodal set $N(\Psi_{1,3}^{\theta})$ for all $\theta$. For geometric reasons, the order of vanishing at a vertex is at least $3$. More precisely,

\begin{properties}\label{S13-vert-P1}
Behaviour of $\Psi_{1,3}^{\theta}$ at the vertices.\vspace{-4mm}
\begin{enumerate}
  \item  The function $\Psi_{1,3}^{\theta}$ vanishes at order $6$ at $O$ if and only if $\theta \equiv 0 \pmod{\pi}$; otherwise it vanishes at order $3\,$.
  \item  The function $\Psi_{1,3}^{\theta}$ vanishes at order $6$ at $A$ if and only if $\theta \equiv \frac{\pi}{3} \pmod{\pi}$; otherwise it vanishes at order $3\,$.
  \item  The function $\Psi_{1,3}^{\theta}$ vanishes at order $6$ at $B$ if and only if $\theta \equiv \frac{2\pi}{3} \pmod{\pi}$; otherwise it vanishes at order $3\,$.
\end{enumerate}\vspace{-4mm}
In other words, up to multiplication by a scalar, the only eigenfunction $\Psi_{1,3}^{\theta}$ which vanishes at higher order at a vertex $V \in \{O,A,B\}$ is $C_V$, the anti-invariant eigenfunction with respect to the median issued from the vertex $V$. In particular, when $\theta \in ]0,\frac{\pi}{6}]$, the three vertices are critical zeros of order three for the eigenfunction $\Psi_{1,3}^{\theta}$, and no interior nodal curve of such an eigenfunction can arrive at a vertex.
\end{properties}

\textbf{Proof}. Taking \eqref{E13-sym4} into account, it suffices to compute the Taylor expansions of $\Psi_{1,3}^{\theta}$ at the  point $(0,0)$. \hfill $\square$

\subsection{Fixed points on the medians}\label{S13-FP}

Since the median $[OM]$ is contained in the nodal set $N(C_{1,3})$, the intersection points of $[OM]$ with  $N(S_{1,3})$ are fixed points of the family $N(\Phi_{1,3}^{\theta})$, \ie, common zeros of the functions $\Psi_{1,3}^{\theta}$. If we parametrize $[OM]$ by $u \mapsto (u,u)$ with $u \in [0,\frac{1}{2}]$, we find that
%
\begin{equation}\label{FP13-2}
S_{1,3}\big|_{[OM]}(u) = 2\sin(8\pi u)-2\sin(6\pi u)-2\sin(2\pi u)\,.
\end{equation}
It follows that
\begin{equation}\label{FP13-4}
S_{1,3}\big|_{[OM]}(u) = - 8 \sin(\pi u)\, \sin(3\pi u)\, \sin(4\pi u)\,.
\end{equation}


The last formula shows that there are two fixed points on the open median $[OM]$, the centroid of the triangle $F_C = (\frac{1}{3},\frac{1}{3})_{\cF}$, and the point $F_O = (\frac{1}{4},\frac{1}{4})_{\cF}$.

Taking into account the action of $G_{\cT}$ on the space $\cE_{5}$, see \eqref{E13-sym2}-\eqref{E13-sym4a}, we infer that the points $F_A = (\frac{5}{12},\frac{1}{3})_{\cF}$ and $F_B = (\frac{1}{3},\frac{5}{12})_{\cF}$ are also common zeros for the family $\Psi^{\theta}_{1,3}\,$. They are deduced from $F_O$ by applying the rotations $\rho_{\pm}\,$, and situated on the two other medians.

Using Taylor expansions, it is easy to check that the fixed points $F_{*}$ are not critical zeros of the functions $\Psi_{1,3}^{\theta}\,$. In the neighborhood of the fixed points $F_{*}$, the nodal set consists of a single regular arc. \medskip

\textbf{Remarks}. \\
(i)~ Note that we do not claim to have determined all the fixed points of the family $N(\Psi_{1,3}^{\theta})$ \ie, the set $N(C_{1,3}) \cap N(S_{1,3})$. We have so far only determined the fixed points located on the medians.\\
 (ii)~ Notice that the fixed point $F_O$ is the mid-point of the segment $[M_AM_B]$, \etc.

\begin{figure}
  \begin{minipage}[b]{0.45\linewidth}
    \centering
    \includegraphics[width=1\linewidth]{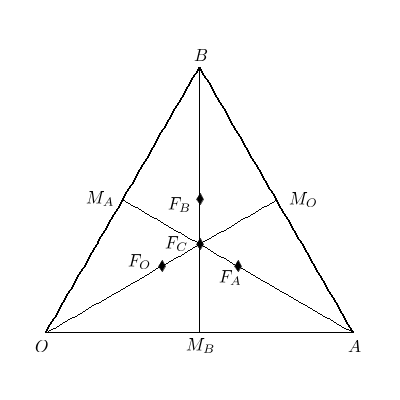}
    \par\vspace{-25pt}
  \end{minipage}%
  \begin{minipage}[b]{0.45\linewidth}
    \centering
\begin{tabular}[b]{|c|c|c|}
\hline
Point & $\E^2$ & $\cF$\\[5pt] \hline
Vertex $O$ & $(0,0)$ & $(0,0)_{\cF}$\\[5pt] \hline
Vertex $A$ & $(1,0)$ & $(\frac{2}{3},\frac{1}{3})_{\cF}$\\[5pt] \hline
Vertex $B$ & $(\frac{1}{2},\frac{\sqrt{3}}{2})$ & $(\frac{1}{3},\frac{2}{3})_{\cF}$\\[5pt] \hline
Fixed point $F_C$ & $(\frac{1}{2},\frac{\sqrt{3}}{6})$ & $(\frac{1}{3},\frac{1}{3})_{\cF}$\\[5pt] \hline
Fixed point $F_O$ & $(\frac{3}{8},\frac{\sqrt{3}}{8})$ & $(\frac{1}{4},\frac{1}{4})_{\cF}$\\[5pt] \hline
Fixed point $F_A$ & $(\frac{5}{8},\frac{\sqrt{3}}{8})$ & $(\frac{5}{12},\frac{1}{3})_{\cF}$\\[5pt] \hline
Fixed point $F_B$ & $(\frac{1}{2},\frac{\sqrt{3}}{4})$ & $(\frac{1}{3},\frac{5}{12})_{\cF}$\\[5pt] \hline
\end{tabular}
    \par\vspace{0pt}
    \end{minipage}
\caption{Triangle $\cT$ and fixed points for $\Psi_{1,3}^{\theta}\,$} \label{FP13-F1}
\end{figure}

\subsection{Partial barriers for the nodal sets}\label{S13-MB}

We have seen that the family $\Psi_{1,3}^{\theta}$ has precisely four fixed points  on the medians, the centroid $F_C$ of the triangle, and three other points $F_O$, $F_A$ and $F_B$, located respectively on the open medians $[OM]$, $[AM]$, and $[BM]$, see Figure~\ref{FP13-F1}. As a matter of fact, the medians can serve as partial barriers.

\begin{lemma}\label{MB13-L1}
For any $\theta$, the nodal set $N(\Psi_{1,3}^{\theta})$ intersects each median at exactly two points unless the function $\Psi_{1,3}^{\theta}$ is one of the functions $C_V$ for $V \in \{O,A,B\}$, in which case the corresponding median is contained in the nodal set.  In particular, if $\theta \in ]0,\frac{\pi}{6}]\,$, the nodal set $N(\Psi_{1,3}^{\theta})$ only meets the medians at the fixed points $\{F_O, F_A, F_B, F_C\}$.
\end{lemma}

\textbf{Proof}. Use the following facts:\\ (i) the families $\{C_O,S_O\}$, $\{C_A,S_A\}$ and $\{C_B,S_B\}$ span $\cE_5\,$;\\
(ii) the function $C_V$ vanishes on the median issued from the vertex $V$. Write $\Psi_{1,3}^{\theta}=\alpha C_V + \beta S_V\,$. If $x \in [VM]\cap N(\Psi_{1,3}^{\theta})$, then $\beta S_V(x)=0\,$.
If $\beta \not = 0\,$, then $x \in \{F_C,F_V\}$. If $\beta=0$, then $[VM] \subset N(\Psi_{1,3}^{\theta})\,$.\hfill $\square$\medskip

\textbf{Remark}. A consequence of Lemma~\ref{MB13-L1} and Subsection~\ref{S13-FP} is that for $\theta \in ]0,\frac{\pi}{6}]$, no critical zero of the function $\Psi_{1,3}^{\theta}$ can occur on the medians.\medskip

The medians divide the equilateral triangle $\cT$ into six isometric $H$-triangles (\ie, triangles with angles $\{30, 60, 90\}$ degrees), $T(F_C,O,M_B)$, \etc . Each one of theses triangles is homothetic to the $H$-triangle $T(O,A,M_O)$ with scaling factor $\frac{1}{\sqrt{3}}$. On the other-hand, the first Dirichlet eigenvalue for the $H$-triangle $T(O,A,M_O)$ is $\lambda_2(\cT) = 7\,$, with multiplicity $1$, and associated eigenfunction $C_{1,3}\,$.  The following lemma follows easily.

\begin{lemma}\label{MB13-23-L}
Let $\cT_M$ denote the triangle $T(F_C,O,M_B)$. The first Dirichlet eigenvalue of the six $H$-triangles determined by the medians of $\cT$ equals $$\lambda(\cT_M) = 3 \lambda_2(\cT) = 21\,.$$
\end{lemma}%

For $0 \le a \le 1$, define the line $D_{a}$ by the equation $s+t=a$ in the parametrization $\cF$.
\begin{equation}\label{MB13-23Da}
D_{a} = \cF \left( \{s+t=a\}\right)\,.
\end{equation}

\begin{lemma}\label{MB13-L2}
For $a \in \{\frac{1}{2}, \frac{2}{3}, \frac{3}{4}\}$, the intersections of the lines $D_{a}$ with the nodal sets $N(C_{1,3})$ and $N(S_{1,3})$ are as follows.
\begin{equation}\label{MB13-4}
\begin{array}{lll}
D_{\frac{1}{2}} \cap N(C_{1,3}) = \{F_{O}\}, &\text{~and~} & D_{\frac{1}{2}} \cap N(S_{1,3}) = \{F_{O}\},\\
D_{\frac{2}{3}} \cap N(C_{1,3}) = \{F_C, G_A,G_B\}, &\text{~and~} & D_{\frac{2}{3}} \cap N(S_{1,3}) = \{F_{C}\},\\
D_{\frac{3}{4}} \cap N(C_{1,3}) = \{F_{A},F_{B},G_O\}, &\text{~and~} & D_{\frac{3}{4}} \cap N(S_{1,3}) = \{F_{A},F_{B}\},\\
\end{array}
\end{equation}
where $G_A$ and $G_B$ are symmetric with respect to $[OM]$, and $G_O =(\frac{3}{8},\frac{3}{8})_{\cF}$. The lines are tangent to $N(S_{2,3})$ at the points $F_{O}$ and $F_{C}$.
\end{lemma}%

\textbf{Proof}. The segment $D_{a}\cap \cT$ is parametrized by $u \mapsto (u,a-u)$ for $u \in [\frac{a}{3},\frac{2a}{3}]$. For each value $a \in \{\frac{1}{2}, \frac{2}{3}, \frac{3}{4}\}$, define
\begin{equation}\label{MB13-6}
\begin{array}{l}
BC_{a}(u) =  C_{1,3} (u,a-u), ~~ u \in [\frac{a}{3},\frac{2a}{3}]\,,\\
BS_{a}(u) = S_{1,3} (u,a-u), ~~ u \in [\frac{a}{3},\frac{2a}{3}]\,.
\end{array}
\end{equation}

Taking $a = \frac{1}{2}$, we find
\begin{equation}\label{MB13-6-12a}
\begin{array}{l}
BC_{\frac{1}{2}}(u) = -4\sin(2\pi u)\sin(12\pi u)\,,\\
BS_{\frac{1}{2}}(u) = 4\cos(2\pi u)\sin(12\pi u)\,.
\end{array}
\end{equation}

Similarly,

\begin{equation}\label{MB13-6-23b}
\begin{array}{l}
BC_{\frac{2}{3}}(u) =  -2\sqrt{3} \sin(9\pi u) \cos(5\pi u + \frac{\pi}{3})\,,\\
BS_{\frac{2}{5}}(u) = -2\sqrt{3} \sin(9\pi u) \sin(5\pi u + \frac{\pi}{3})\,.
\end{array}
\end{equation}

\begin{equation}\label{MB13-6-34b}
\begin{array}{l}
BC_{\frac{3}{4}}(u) =  -2 \sin(12\pi u)\left( \sin(2\pi u) + \cos (2\pi u)\right) \,,\\
BS_{\frac{3}{4}}(u) = -2 \sin(12\pi u)\left( \sin(2\pi u) - \cos (2\pi u)\right) \,.
\end{array}
\end{equation}

Looking at the zeros of the above functions in the respective intervals, the lemma follows. \hfill $\square$

\begin{figure}
  \centering
  \includegraphics[width=14cm]{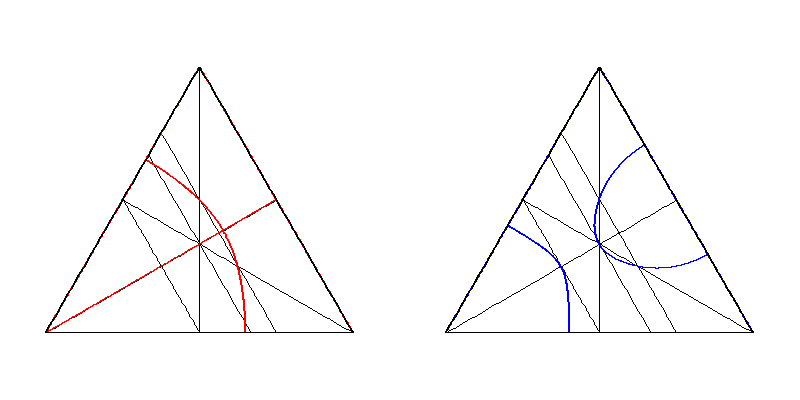}
  \caption{Partial barriers for $N(C_{1,3})$ and $N(S_{1,3})$}\label{MB13-F1}
\end{figure}

The lemmas are illustrated by Figure~\ref{MB13-F1} which displays the partial barriers (thin segments), and the nodal sets computed with Maple (thicker lines).

\subsection{Critical zeros of $C_{1,3}$ and $S_{1,3}$ on the sides of $\cT$, \\and on the median $[OM]$}\label{S13-C}

Define the functions
\begin{equation}\label{13C-FCS}
\begin{array}{ll}
FC(u) &:= - \sin(7 \pi u) + 3 \sin(5\pi u) - 4 \sin(2\pi u)\,,\\
FS(u) &:= - \cos(7 \pi u) - 3 \cos(5\pi u) + 4 \cos(2\pi u)\,,\\
\end{array}
\end{equation}
and the polynomials
\begin{equation}\label{13C-PCS}
\begin{array}{ll}
P_C(x) &:= 4x^2+4x-1\,,\\
P_S(x) &:= 4x^4-x^2+x-1\,.\\
\end{array}
\end{equation}

\begin{lemma}\label{13C-L1}
The functions $FC$ and $FS$ satisfy,
\begin{equation}\label{13C-FG}
\begin{array}{ll}
FC(u) & = -4\sin(\pi u) \left( \cos(\pi u) - 1\right)^2 \left( 2 \cos(\pi u) + 1\right)^2 \, P_C\left( \cos(\pi u) \right)\,,\\
FS(u) & = - 4 \left( \cos(\pi u) - 1\right) \left( 2 \cos(\pi u) + 1\right)^2 \, P_S\left( \cos(\pi u) \right)\,.
\end{array}
\end{equation}
\end{lemma}%

\textbf{Proof}. Use the Chebyshev polynomials.
\hfill $\square$

\begin{properties}\label{13C-P1}
The partial derivatives of the functions $C_{1,3}$ and $S_{1,3}$ satisfy the following relations.\vspace{-4mm}
\begin{enumerate}
\item Parametrize the edge $[OA]$ by $u \mapsto (u,u/2)$, with $u \in [0,2/3]$. Then,
\begin{equation}\label{13C-P1-OA}
\begin{array}{l}
\partial_s C_{1,3}(u,u/2) = 2 \pi FC(u)\,,~~ \partial_t C_{1,3}(u,u/2) = -4 \pi FC(u)\,,\\
\partial_s S_{1,3}(u,u/2) = 2 \pi FS(u)\,,~~ \partial_t S_{1,3}(u,u/2) = -4 \pi FS(u)\,.
\end{array}
\end{equation}
\item Parametrize the edge $[OB]$ by $u \mapsto (u/2,u)$, with $u \in [0,2/3]$. Then,
\begin{equation}\label{13C-P1-OB}
\begin{array}{l}
\partial_s C_{1,3}(u/2,u) = 4 \pi FC(u)\,,~~ \partial_t C_{1,3}(u/2,u) = -2 \pi FC(u)\,,\\
\partial_s S_{1,3}(u/2,u) = -4 \pi FS(u)\,,~~ \partial_t S_{1,3}(u/2,u) = 2 \pi FS(u)\,.
\end{array}
\end{equation}
\item Parametrize the edge $[BA]$ by $u \mapsto (u/2,1-u/2)$, with $u \in [2/3,4/3]$. Then,
\begin{equation}\label{13C-P1-BA}
\begin{array}{l}
\partial_s C_{1,3}(u/2,1-u/2) = -2 \pi FC(u)\,,~~ \partial_t C_{1,3}(u/2,1-u/2) = -2 \pi FC(u)\,,\\
\partial_s S_{1,3}(u/2,1-u/2) = 2 \pi FS(u)\,,~~ \partial_t S_{1,3}(u/2,1-u/2) = 2 \pi FS(u)\,.
\end{array}
\end{equation}
\end{enumerate}\vspace{-3mm}
As a consequence, the critical zeros of $C_{1,3}$ or $S_{1,3}$ on the edges $[OA]$ and $[OB]$, \resp on the edge $[BA]$, are determined by the zeros of $FC$ or $FS$ in the intervals $[0,2/3]$, \resp $[2/3,4/3]$.
\end{properties}%

\textbf{Proof}. It suffices to compute the partial derivatives of $C_{1,3}$ and $S_{1,3}\,$, and to make the substitutions corresponding to the parametrization of the edges.
\hfill $\square$

It follows from the above results that the critical zeros of $C_{1,3}\,$, \resp $S_{1,3}\,$, on $\partial \cT$ are determined by the zeros of $FC$, \resp $FS$, in $[0,4/3]$.

\begin{lemma}\label{13C-L2}
Zeros of $P_C$ and $P_S$.\vspace{-3mm}
\begin{enumerate}
\item The polynomial $P_C$ has precisely one root in the interval $[-1,1]\,$, namely $(\sqrt{2}-1)/2\,$.
\item The polynomial $P_S$ has exactly two roots in the interval $[-1,1]$, namely   $\xi_{-} \approx -0.9094691258$ and $\xi_{+} \approx 0.6638481772\,$.
\end{enumerate}
\end{lemma}%

Define the numbers
\begin{equation}\label{13C-GZ}
\begin{array}{c}
u_{1,C} := \frac{1}{\pi} \arccos(\frac{\sqrt{2}-1}{2}) \approx 0.433595245\,,\\[5pt]
u_{1,S} := \frac{1}{\pi}  \arccos(\xi_{+})  \approx 0.2689221041\,,\\[5pt]
u_{2,S} := \frac{1}{\pi}  \arccos(\xi_{-}) \approx 0.8635116189\,,\\[5pt]
u_{3,S} := 2- u_{2,S} \approx 1.136488381\,.\\
\end{array}
\end{equation}

The function $FC$ vanishes at $0, \frac{2}{3}$, $1$ and $\frac{4}{3}$, and has one simple zero $u_{1,C} \in ]0,\frac{2}{3}[$. The function $FS$ vanishes at $0, \frac{2}{3}$ and $\frac{4}{3}$. It has one simple zero $u_{1,S} \in ]0,\frac{2}{3}[\,$, and two simple zeros $u_{2,S}, u_{3,S} \in ]\frac{2}{3},\frac{4}{3}[\,$.

\begin{properties}\label{13C-P3} Critical zeros of the functions $C_{1,3}$ and $S_{1,3}$ on the open edges of $\cT$.\vspace{-3mm}
\begin{enumerate}
\item The function $C_{1,3}$ has one critical zero $Z_{1,C}=(u_{1,C},u_{1,C}/2)_{\cF}$ on the open edge $[OA]$; one critical zero $Z_{2,C} = (u_{1,C}/2,u_{1,C})_{\cF}$ on the open edge $[OB]$; one critical zero $Z_{3,C} = M_O = (1/2,1/2)_{\cF}$ on the open edge $[BA]$. These critical zeros have order $2$.
\item The function $S_{1,3}$ has one critical zero $Z_{1,S} = (u_{1,S},u_{1,S}/2)_{\cF}$ on the open edge $[OA]$; one critical zero $Z_{2,S} = (u_{1,S}/2,u_{1,S})_{\cF}$ on the open edge $[OB]$; two critical zeros $Z_{3,S} = (u_{2,S}/2, 1-u_{2,S}/2)_{\cF}$ and $Z_{4,S} = (u_{3,S}/2, 1-u_{3,S}/2)_{\cF}$ on the open edge $[BA]$, these points are symmetric with respect to the point $M_O$. These critical zeros have order $2\,$.
\end{enumerate}
\end{properties}%

\textbf{Proof}. Use Properties~\ref{13C-P1} and Lemma~\ref{13C-L2}. \hfill $\square$ \medskip

\textbf{Remark}. The vertex $O$ is a critical zero of order $6$ of $C_{1,3}\,$, and a critical zero of order $3$ of $S_{1,3}\,$. The vertices $A$ and $B$ are critical zeros of order $3$ of both $C_{1,3}$ and $S_{1,3}\,$, see Properties~\ref{S13-vert-P1}.

\begin{figure}
  \begin{minipage}[b]{0.40\linewidth}
    \centering
    \includegraphics[width=1\linewidth]{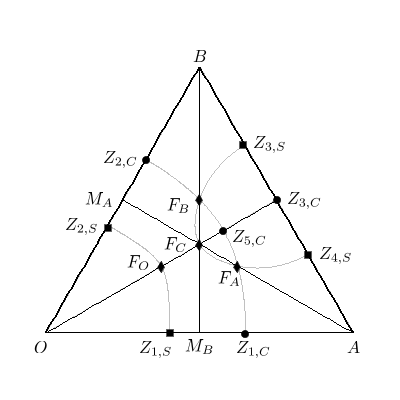}
    \par\vspace{-25pt}
  \end{minipage}%
  \begin{minipage}[b]{0.55\linewidth}
    \centering
\begin{tabular}[b]{|c|c|c|}
\hline
Point & $\E^2$ & $\cF$\\[5pt] \hline
$Z_{1,C}$ & $\approx (0.6504,0)$ & $(u_{1,C}\,,\frac{u_{1,C}}{2})_{\cF}$\\[5pt] \hline
 $Z_{2,C}$ & $\approx (0.3252,0.5633)$ & $(\frac{u_{1,C}}{2}\,,u_{1,C})_{\cF}$\\[5pt] \hline
$Z_{3,C}$ & $(\frac{3}{4},\frac{\sqrt{3}}{4})$ & $(\frac{1}{2},\frac{1}{2})_{\cF}$\\[5pt] \hline
$Z_{5,C}$ & $\approx (0.5775,0.3334)$ & $(\frac{u_{5,C}}{2}\,,\frac{u_{5,C}}{2})_{\cF}$\\[5pt] \hline
$Z_{1,S}$ & $\approx (0.4034,0)$ & $(u_{1,S}\,,\frac{u_{1,S}}{2})_{\cF}$\\[5pt] \hline
$Z_{2,S}$ & $\approx (0.2017,0.3494)$ & $(\frac{u_{1,S}}{2},u_{1,S})_{\cF}$\\[5pt] \hline
$Z_{3,S}$ & $\approx (0.6476,0.6104)$ & $(\frac{u_{2,S}}{2},1-\frac{u_{2,S}}{2})_{\cF}$\\[5pt] \hline
$Z_{4,S}$ & $\approx (0.8524,0.2507)$ & $(\frac{u_{3,S}}{2},1-\frac{u_{3,S}}{2})_{\cF}$\\[5pt] \hline
\end{tabular}
    \par\vspace{0pt}
    \end{minipage}
\caption{Fixed points and critical zeros for $C_{1,3}$ and $S_{1,3}$} \label{13C-F1}
\end{figure}

\begin{properties}\label{13C-P5}[ Critical zeros of the functions $C_{1,3}$ and $S_{1,3}$ on the median $[OM_O]$, illustrated by Figure~\ref{13C-F1}] \vspace{-3mm}
\begin{enumerate}
\item The function $C_{1,3}$ has one critical zero at $O$; one critical zero $Z_{5,C}$ of order $2\,$, where
    $$
    Z_{5,C} = (u_{5,C}/2,u_{5,C}/2)_{\cF}, \text{~with~}u_{5,C} := 1-\arccos(3/4)/\pi \approx 0.7699465439 \,;
    $$
one critical zero $M_O = Z_{3,C}$ of order $2\,$.
\item The function $S_{1,3}$ has no critical zero on the median $[OM]$, except the point $O$.
\end{enumerate}
\end{properties}%

\textbf{Proof}. Since $C_{1,3}$ vanishes on the median, its critical zeros on the median are the common zeros of its partial derivatives. They are precisely the zeros of the function
$$
2\sin(4\pi u)-5\sin(3\pi u)+7\sin(\pi u)\,,
$$
if we parametrize the median by $u \mapsto (u/2,u/2)$ for $u \in [0,1]$. The above function can be factorized as
$$
4\sin(\pi u)\big( 4\cos(\pi u)+3\big) \big( \cos(\pi u)-1\big)^2\,,
$$
and the first assertion follows. 

For the second assertion, we have to look for the common zeros of the function $S_{1,3}$ and its derivatives on the median. This amounts to finding the common zeros of the functions
$$
2\sin(4\pi u)-2\sin(3\pi u)-2\sin(\pi u)\,,
$$
and
$$
4\cos(4\pi u)-3\cos(3\pi u)-\cos(\pi u)\,.
$$

The first function factorizes as
$$
-8 \sin(\frac{\pi u}{2}) \sin(\frac{3\pi u}{2})\sin(2\pi u)\,.
$$
It is easy to check that the only common zero is $u=0\,$. \hfill $\square$

We state the following corollary of Proposition~\ref{13C-P5} for later reference.

Recall (see the notation \eqref{E13-sym4}) that $\Psi_{1,3}^{\frac{\pi}{6}} = - S_B$, and that $S_B = \rho_{+}^{*}S_O\,$.

\begin{corollary}\label{13C-C1}
Critical zeros of the function $S_B\,$. \vspace{-3mm}
\begin{enumerate}
  \item The function $S_B$ has two critical zeros of order $2$ on the side $[OA]$,
    $$
    Z_{3,S_B} = \rho_{-}(Z_{3,S}) \text{~and~} Z_{4,S_B} =
    \rho_{-}(Z_{4,S})\,.
    $$
  \item The function $S_B$ has one critical zero of order $2$ on the side $[OB]$,
    $$
    Z_{1,S_B} =\rho_{-}(Z_{1,S})\,.
    $$
  \item The function $S_B$ has one critical zero of order $2$ on the side $[BA]$,
    $$
    Z_{2,S_B} = \rho_{-}(Z_{2,S})\,.
    $$
\end{enumerate}
\end{corollary}%

\subsection{The nodal sets of $C_{1,3}$ and $S_{1,3}$}\label{S13-NCS}

\begin{properties}\label{NCS13-P2} Nodal sets of $C_{1,3}$ and $S_{1,3}\,$.\vspace{-4mm}
\begin{enumerate}
  \item The function $C_{1,3}$ has only one critical zero $Z_{5,C}$ in the interior of the triangle. Its nodal set consists of the diagonal $[OM_O]$, and a regular arc from $Z_{1,C}$ to $Z_{2,C}$ which intersects $[OM_O]$ orthogonally at $Z_{5,C}$, and passes through $F_A$ and $F_B$.
  \item The function $S_{1,3}$  has no critical zero in the interior of the triangle. Its nodal set consists of two disjoint regular arcs, one from $Z_{1,S}$ to $Z_{2,S}$, passing through $F_O$; one from $Z_{3,S}$ to $Z_{4,S}$, passing through $F_B$ and $F_A$.
\end{enumerate}
\end{properties}%

\begin{figure}[b!]
  \centering
  \includegraphics[width=12cm]{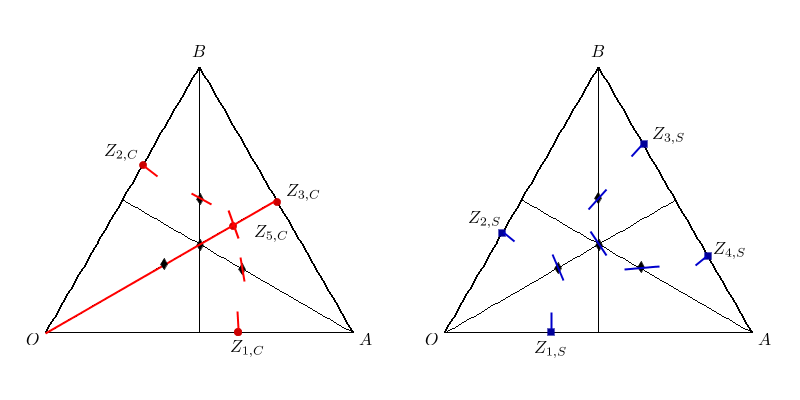}
  \caption{Local nodal patterns for $C_{1,3}$ (left) and $S_{1,3}$ (right)}\label{NCS13-F1}
\end{figure}

\textbf{Proof}. We have determined the common zeros of $C_{1,3}$ and $S_{1,3}$ located on the medians (Subsection~\ref{S13-FP}), as well as the critical zeros on the open edges of the triangle $\cT$ and on the medians $[OM]$ (Subsection~\ref{S13-C}). We already know the local behaviour at the vertices (Subsection~\ref{S13-vert}). Using Subsection~\ref{S13-MB}, we also know that the nodal set $N(C_{1,3})$ only meets the medians at the fixed points and at $Z_{5,C}$, and that the nodal set $N(S_{1,3})$ only meets the medians at the fixed points. Looking at the Taylor expansions, we can determine the local nodal patterns of $C_{1,3}$ and $S_{1,3}$ near the fixed points and near the critical zeros, see Figure~\ref{NCS13-F1}. This figure also displays the medians, and takes into account the fact that $[OM_0] \subset N(C_{1,3})$.

The medians divide $\cT$ into six isometric $H$-triangles, Figure~\ref{13C-F1}. The nodal sets of $C_{1,3}$ and $S_{1,3}$ consist of finitely many nodal arcs which are smooth except at the critical zeros.  According to our preliminary results, such arcs can only exit the interior of an $H$-triangle at a fixed point or at a critical zero, with the local nodal patterns shown in Figure~\ref{NCS13-F1}. \medskip

\emph{Claim}. The functions $C_{1,3}$ and $S_{1,3}$ cannot have any critical zero in the interiors of the $H$-triangles.\smallskip

Indeed, assume there is one critical zero $Z$ in the interior of some $H$-triangle $\cH$. At this point, the nodal set would consist of at least four semi-arcs. Following any such semi-arc, we either obtain a simply closed nodal arc, or exit the triangle. Since there are at most three exit points (with only one exit direction at each point), there would be at least one simply closed nodal component in the interior of the triangle $\cH$. This component would bound at least one nodal domain $\omega$. The first Dirichlet eigenvalue $\lambda(\omega)$ would satisfy $\lambda(\omega) = \lambda_5(\cT) = 13$. On the other-hand, since $\omega$ is contained in the interior of $\cH$, we would have $\lambda(\omega) > \lambda(\cH) = 21$ according to Lemma~\ref{MB13-L2}. This proves the claim by contradiction. \medskip

The last argument in the proof of the claim also shows that the interiors of the $H$-triangles cannot contain any closed nodal component. This shows that the nodal sets of $C_{1,3}$ and $S_{1,3}$ are indeed as shown in Figure~\ref{ETR-F4}. \hfill $\square$ \medskip

\textbf{Remark}.\\
Here is another argument to determine the nodal set $N(C_{1,3})$. The claim is that $C_{1,3}$ cannot have a second interior critical zero. Indeed, such a critical zero would be of order at least $2$, and cannot belong to the median $[OM]$, so that it must belong to one of the two $H$-triangles determined by this median. For symmetry reasons, we would have at least one critical zero in each of these $H$-triangles. The eigenvalue $\lambda_5(\cT)$ which corresponds to the pairs $[1,3]$ and $[3,1]$ is also the second Dirichlet eigenvalue of these triangles, so that it has two nodal domains, $\Omega_1$ and $\Omega_2$. Applying the Euler formula \cite[(2.11)]{BeHe0} to the nodal partition $\mathcal D =(\Omega_1,\Omega_2)$, we would have
$$
1 + \frac{1}{2}(1+1+2) \le \chi(\frac{1}{2}\cT) + \frac{1}{2}\sigma(\Omega_1,\Omega_2) = \chi(\Omega_1) + \chi(\Omega_2) \le 2\,,
$$
a contradiction.

\subsection{Critical zeros of $\Psi_{1,3}^{\theta}$ on the sides of $\cT$}\label{S13-CP}

As a consequence of Properties~\ref{13C-P1}, the critical zeros of the functions $\Psi_{1,3}^{\theta}\,$, for $\theta \in ]0,\frac{\pi}{6}]$, on the sides of the triangle $\cT$ are determined by one of the equations
$$
\cos\theta \, FC(u) \pm \sin\theta \, FS(u) = 0\,.
$$

Since $\theta \in ]0,\frac{\pi}{6}]$, the vertices of $\cT$ are critical zeros of order $3$ of $\Psi_{1,3}^{\theta}$ (Properties~\ref{S13-vert-P1}). Since we are interested in the critical zeros on the open edges, we can substitute $FC$, \resp $FS$, by the functions $GC$, \resp $GS$, defined as follows.

\begin{equation}\label{13C-GCS}
\begin{array}{ll}
GC(u) &:= \sin(\pi u) \left( cos(\pi u) - 1\right) \left( 4 \cos^2(\pi u) + 4 \cos(\pi u) -1\right)\,,\\
GS(u) &:= 4 \cos^4(\pi u) - \cos^2(\pi u) + \cos(\pi u) -1\,.
\end{array}
\end{equation}

\begin{properties}\label{13CP-P2}
The critical zeros of the function $\Psi_{1,3}^{\theta}$ on the open edges of the triangle $\cT$ are determined by the following equations. \vspace{-3mm}
\begin{enumerate}
\item On the edge $[OA]$ parametrized by $u \mapsto (u,u/2)$,
  \begin{equation}\label{13C-P2-OA}
  \cos\theta \, GC(u) + \sin\theta \, GS(u) = 0\,, ~\text{~for~} u \in [0,2/3]\,.
  \end{equation}
\item On the edge $[OB]$ parametrized by $u \mapsto (u/2,u)$,
  \begin{equation}\label{13C-P2-OB}
  \cos\theta \, GC(u) - \sin\theta \, GS(u) = 0\,, ~\text{~for~} u \in [0,2/3]\,.
  \end{equation}
\item On the edge $[BA]$ parametrized by $u \mapsto (u/2,1-u/2)$,
  \begin{equation}\label{13C-P2-BA}
  \cos\theta \, GC(u) - \sin\theta \, GS(u) = 0\,, ~\text{~for~} u \in [2/3,4/3]\,.
  \end{equation}
\end{enumerate}
\end{properties}%

For convenience, we introduce the functions,
\begin{equation}\label{13C-H0}
H_{\pm}^{\theta}(u) = \cos\theta \, GC(u) \pm \sin\theta \, GS(u)\,.
\end{equation}

\begin{properties}\label{13CP-P3}
Recall the notation \eqref{13C-GZ}. \vspace{-3mm}
\begin{enumerate}
  \item In the interval $]0,2/3[$ (corresponding to critical zeros on the open side $[OA]$), the function $H_{+}^{\theta}$ has two simple zeros,
$$
\beta_1(\theta) \in\, ]0\,,\,u_{1,S}[ \text{~and~} \beta_2(\theta) \in\, ]u_{1,C}\,,\, \frac{2}{3}[\,.
$$
They are smooth increasing functions in $\theta$.
  \item In the interval $]0\,,\,2/3[$ (corresponding to critical zeros on the open side $[OB]$), the function $H_{-}^{\theta}$ has one simple zero
$$
\alpha_1(\theta) \in\, ]\frac{1}{3}\,, \,u_{1,C}[\,.
$$
This is a smooth decreasing function of $\theta$.
  \item In the interval $]2/3,4/3[$ (corresponding to critical zeros on the open side $[BA]$), the function $H_{-}^{\theta}$ has one simple zero,
$$
\omega_1(\theta) \in\, ]u_{2,S}\,,\,1[\,.
$$
This is a smooth decreasing function of $\theta$.
\end{enumerate}
\end{properties}%

\textbf{Proof}. Notice that the zeros of $H_{\pm}^{\theta}$ are continuous functions of $\theta$ because the equations $H_{\pm}^{\theta}(u)=0$ can be written as polynomial equations in $\tan(\frac{\pi u}{2})$ with coefficients depending continuously on $\theta$. We study the functions $H_{\pm}^{\theta}$ in the interval $]-1/6,3/2[$ which contains the interval $[0,4/3]$.

First of all, taking into account the fact that $\theta \in ]0,\pi/6]$, we look at the values of the functions $H_{\pm}^{\theta}$ at the points
$$
0 < u_{1,S} < \frac{1}{3} < u_{1,C} < \frac{2}{3} < u_{2,S} < 1 < u_{3,S} < \frac{4}{3}\,,
$$
and infer the existence of at least one zero in each of the intervals mentioned in the statements. Note that the zero $\beta_1(\theta)$ comes from the fact that the vertex $O$ has order $6$ for $C_{1,3}$ and order $3$ for $\Psi_{1,3}^{\theta}$ as soon as $0 < \theta \le \pi/6$. The details appear in Table~\ref{13CP-T1} (in which we have only indicated the useful values). The values listed above correspond to the vertices of the triangles ($0,  \frac{2}{3}$ and $\frac{4}{3}$), the mid-points on the edges ($\frac{1}{3}$ and $1$), and the critical zeros on the open edges.

\begin{table}[hbt]
\centering
\begin{tabular}{|c|c|c|c|c|}
  \hline
$u \in [0,4/3]$ & $GC(u)$ & $GS(u)$ & $H_{+}^{\theta}(u)$ & $H_{1}^{\theta}(u)$ \\ \hline
$0$ & $0$ & $3$ & $3\sin\theta$ & $-3\sin\theta$ \\ \hline
$\color{red}{\beta_1(\theta)}$ & -- & -- & $\color{red}{0}$ & -- \\ \hline
$u_{1,S}$ & $GC(u_{1,S})$ & $0$ & $GC(u_{1,S})\cos\theta$ & $GC(u_{1,S})\cos\theta$ \\
$\approx 0.2689$  & $\approx -0.8660$ &  & $\approx -0.8660\cos\theta$ & $\approx -0.8660\cos\theta$ \\ \hline
 $\frac{1}{3}$ & $GC(\frac{1}{3}) = -\frac{\sqrt{3}}{2}$ & $GS(\frac{1}{3}) = -\frac{1}{2}$ & $-\cos(\theta-\frac{\pi}{6})$ & $-\cos(\theta+\frac{\pi}{6})$ \\
  & $\approx -0.866025$ &  &  &  \\ \hline
$\color{red}{\alpha_1(\theta)}$ & -- & -- & -- & $\color{red}{0}$ \\ \hline
$u_{1,C}$ & $0$ & $GS(u_{1,C})$ & $GS(u_{1,C})\sin\theta$ & $-GS(u_{1,C})\sin\theta$ \\
$\approx 0.4336$ &  & $\approx -0.8284$ & $\approx -0.8284 \sin\theta$ & $\approx 0.8284 \sin\theta$ \\ \hline
$\color{red}{\beta_2(\theta)}$ & -- & -- & $\color{red}{0}$ & -- \\ \hline
$\frac{2}{3}$ & $GC(\frac{2}{3}) = \frac{3\sqrt{3}}{2}$ & $GS(\frac{2}{3}) = -\frac{3}{2}$ & $3\cos(\theta+\frac{\pi}{6})$ & $3\cos(\theta-\frac{\pi}{6})$ \\
  & $\approx 2.5981$ &  &  &  \\ \hline
$u_{2,S}$ & $GC(u_{2,S})$ & $0$ & -- & $GC(u_{2,S})\cos\theta$ \\
$\approx 0.8635$ & $\approx 1.0554$ &  & -- & $\approx 1.0554\cos\theta$ \\ \hline
$\color{red}{\omega_1(\theta)}$ & -- & -- & -- & $\color{red}{0}$ \\ \hline
$1$ & $0$ & $1$ & -- & $-\sin\theta$ \\ \hline
$u_{3,S}$ & $GC(u_{3,S})$ & $0$ & -- & $GC(u_{3,S})\cos\theta$ \\
$1.1365$ & $\approx -1.0554$ &  & -- & $\approx -1.0554\cos\theta$ \\ \hline
$\frac{4}{3}$ & $GC(\frac{4}{3}) = -\frac{3\sqrt{3}}{3}$ & $GS(\frac{4}{3})=-\frac{3}{2}$ & -- & $-3\cos(\theta+\frac{\pi}{6})$ \\
 & $\approx -2.5981$ &  & -- &  \\ \hline
\end{tabular}
\caption{Values of $H_{\pm}^{\theta}\,$}\label{13CP-T1}
\end{table}

We then investigate whether the zeros of $H_{\pm}^{\theta}$ can have order at least $2$. More precisely, we investigate whether there exists a pair $(\theta,u)$ such that
\begin{equation}\label{13C-H1}
\begin{array}{l}
\cos\theta \, GC(u) \pm \sin\theta \, GS(u) = 0\,,\\
\cos\theta \, GC'(u) \pm \sin\theta \, GS'(u) = 0\,.
\end{array}
\end{equation}

For this purpose, we define the function
\begin{equation}\label{13C-HW1}
WCS (u) := GC(u) GS'(u) - GS(u) GC'(u).
\end{equation}

\begin{lemma}\label{13C-LW}
The function $WCS$ satisfies the relation
\begin{equation}\label{13C-HW2}
WCS(u) = \pi \left( 1- \cos(\pi u) \right) \left( 2\cos(\pi u) +1 \right)^2 \left( 13 \cos^3(\pi u) - 9 \cos(\pi u) + 4\right)\,.
\end{equation}
Furthermore, the function $WCS$ is non-negative, and vanishes in the interval $]-1/6,3/2[$ if and only if $u \in \{0, \frac{2}{3}, \frac{4}{3}\}$.
\end{lemma}%

\emph{Proof of the Lemma}. Compute the derivatives, make use of the Chebyshev polynomials, and notice that the polynomial of degree $3$ in $\cos(\pi u)$ is always bigger than or equal to~$1$. \hfill $\square$

Lemma~\ref{13C-LW} implies that for $\theta \in ]0, \pi/6]$, the zeros of the functions $H_{\pm}^{\theta}$ in the interval $]0,\frac{2}{3}[$ and $]\frac{2}{3},\frac{4}{3}[$ are simple, so that they are smooth in $\theta$. If $u(\theta)$ is such a zero, its derivative with respect to $\theta$ satisfies the relation
$$
1+\tan^2(\theta) = \pm \frac{WCS(u(\theta))}{GS^2(u(\theta))} u'(\theta)\,.
$$

We can now start from the function $\Psi_{1,3}^{\frac{\pi}{6}} = -S_B$, and follow the zeros by continuity, using Corollary~\ref{13C-C1}. This proves Properties~\ref{13CP-P3}. \hfill $\square$\medskip

Recall the notation of Corollary~\ref{13C-C1}.

\begin{corollary}\label{13CP-C1}
Critical zeros of the function $\Psi_{1,3}^{\theta}$ for $\theta \in ]0,\frac{\pi}{6}]$.\vspace{-3mm}
\begin{enumerate}
  \item There are two critical zeros of order $2$ on the side $[OA]$, one $Z_{4,\theta}$ in the segment $]O,Z_{4,S_B}]$, and one $Z_{3,\theta}$ in the segment $]Z_{1,C},Z_{3,S_B}]$.
  \item There is one critical zero $Z_{1,\theta}$ of order $2$ on the segment $[Z_{1,S_B},Z_{2,C}[$.
  \item There is one critical zero $Z_{2,\theta}$ of order $2$ on the segment $[Z_{2,S_B},Z_{3,C}[$ (recall that $Z_{3,C}=M_O$).
\end{enumerate}
\end{corollary}%

\subsection{Nodal set of $\Psi_{1,3}^{\theta}$}\label{S13-NSP}

\begin{proposition}\label{NSP13-P1}
For $\theta \in ]0,\frac{\pi}{6}]$ the nodal set of $\Psi_{1,3}^{\theta}$ consists of two disjoint simple arcs, one from $Z_{1,\theta}$ to $Z_{2,\theta}$, through $F_{B}$; another from $Z_{4,\theta}$ to $Z_{3,\theta}\,$, through the points $\{F_{0}, F_{C}$ and $F_{A}\}$.  In particular, the function $\Psi_{1,3}^{\theta}$ has three nodal domains. As a consequence, the eigenvalue $\lambda_5(\cT)=13$ is not Courant-sharp.
\end{proposition}%

\begin{figure}[htb!]
  \centering
  \includegraphics[width=12cm]{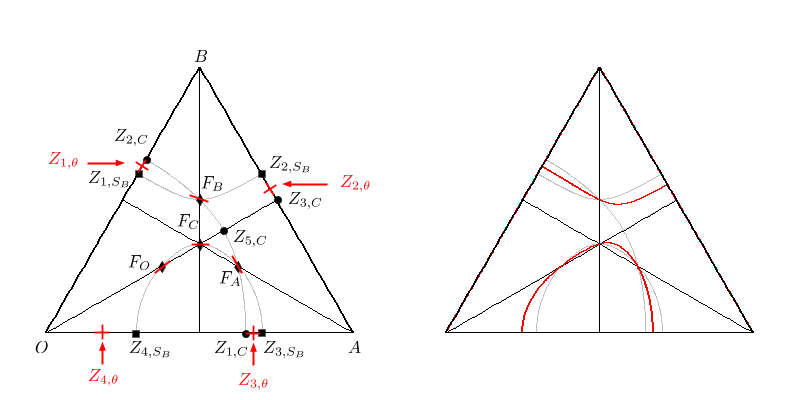}
  \caption{Nodal patterns and nodal set for $\Psi_{1,3}^{\theta}\,$}\label{NSP13-F1}
\end{figure}

\begin{figure}[htb!]
  \centering
  \includegraphics[width=12cm]{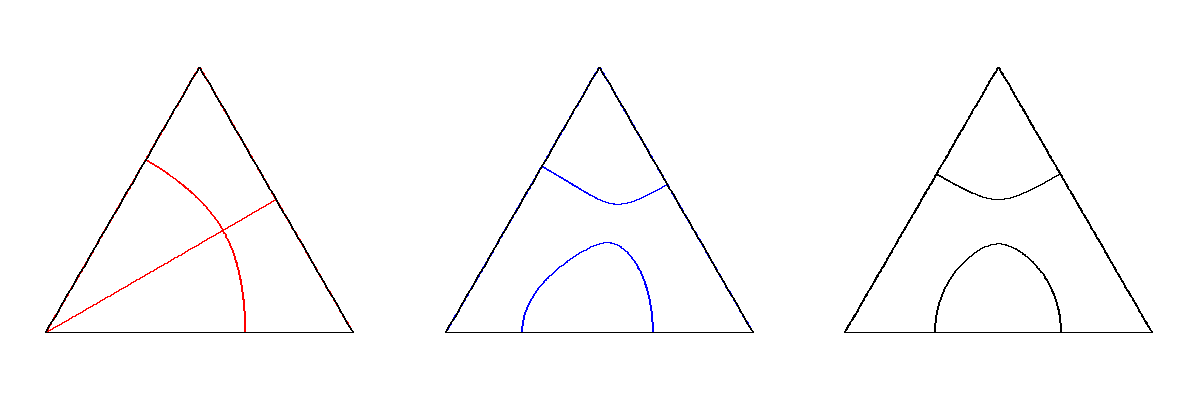}
  \caption{Nodal set of $\Psi_{1,3}^{\theta}\,$, for $\theta = 0\,, \frac{\pi}{12}$ and $\frac{\pi}{6}\,$}\label{NSP13-F2}
\end{figure}

\textbf{Proof}. It is similar to the proof of Properties~\ref{NCS13-P2}, and illustrated by Figure~\ref{NSP13-F1}. The picture on the left-hand side displays the local nodal patterns, from which we can deduce that there are no interior critical zeros in any of the six $H$-triangles determined by the medians. The picture on the right-hand side displays the nodal set $N(\Psi_{1,3}^{\frac{\pi}{12}})$ computed with Maple. The pictures in Figure~\ref{NSP13-F2} show the nodal sets of the eigenfunctions $C_{1,3}\,$, $\Psi_{1,3}^{\frac{\pi}{12}}$ and $\Psi_{1,3}^{\frac{\pi}{6}}$ (from left to right). \hfill $\square$

 \textbf{Remark}. Let $\Psi$ be an eigenfunction. Once one knows the critical zeros of $\Psi$, together with their orders, and the number of connected components of $N(\Psi) \cup \partial \cT$, one can apply the Euler-type formula of \cite[Proposition~2.8]{HOMN} to obtain the number of nodal domains of $\Psi$. Using the proof of Proposition~\ref{NSP13-P1} and this formula, one can recover the number of nodal domains of $\Psi_{1,3}^{\theta}$\,: three for any $\theta \in ]0,\frac{\pi}{6}]$.

\section{The eigenvalue $\lambda_7(\cT)$ and its eigenspace}\label{SS-E23}

We call $\cE_{7}$ the $2$-dimensional eigenspace associated with the eigenvalue $\lambda_7(\cT)$, \ie, with the pairs $[2,3]$ and $[3,2]$. Recall the eigenfunctions,
\begin{equation}\label{S23-0}
\begin{array}{ll}
C_{2,3}(s,t) = & \cos 2\pi(2s+3t) - \cos 2\pi(-2s+5t)
- \cos 2\pi(5s-3t)\\
& - \cos 2\pi(-3s-2t) + \cos 2\pi(3s-5t) + \cos 2\pi(-5s+2t)\,,\\[5pt]
S_{2,3}(s,t) = & \sin2\pi(2s+3t) - \sin2\pi(-2s+5t) - \sin 2\pi(5s-3t)\\
& - \sin 2\pi(-3s-2t) + \sin 2\pi(3s-5t) + \sin 2\pi(-5s+2t)\,.
\end{array}
\end{equation}

 In this section, we study  more carefully the nodal sets of the family $\Psi_{2,3}^{\theta} = \cos\theta \, C_{2,3} + \sin\theta \, S_{2,3}$, following the strategy sketched in Section~\ref{sProl}. The first five subsections contain preliminary results. In Subsection~\ref{S23-NCS}, we determine rigorously the nodal sets of the functions $C_{2,3}$ and $S_{2,3}$. In Subsection~\ref{S23-NSP}, we determine the nodal sets of the function $\Psi^{\theta}_{2,3}$ when $\theta \in ]0,\frac{\pi}{6}]$. As a consequence, we obtain that the number of nodal domains is at most $4$, so that the eigenvalue $\lambda_7(\cT)$ is not Courant-sharp.

\subsection{Symmetries}\label{S23-sym}

Taking Subsection~\ref{S-sym} into account, we find that
\begin{equation}\label{E23-sym2}
\begin{array}{lll}
\sigma_1^{*}\Psi_{2,3}^{\theta} &=& \Psi_{2,3}^{\pi-\theta}\,,\\[5pt]
\sigma_2^{*}\Psi_{2,3}^{\theta} &=& \Psi_{2,3}^{\frac{5\pi}{3}-\theta}\,,\\[5pt]
\sigma_3^{*}\Psi_{2,3}^{\theta} &=&
\Psi_{2,3}^{\frac{\pi}{3}-\theta} \,,\\[5pt]
(\sigma_2 \circ \sigma_1)^{*}\Psi_{2,3}^{\theta} &=& \Psi_{2,3}^{\frac{4\pi}{3}+\theta}\,,\\[5pt]
(\sigma_1 \circ \sigma_2)^{*}\Psi_{2,3}^{\theta} &=& \Psi_{2,3}^{\frac{2\pi}{3}+\theta}\,.
\end{array}
\end{equation}

It follows that, up to multiplication by scalars, and for $V \in \{O,A,B\}$, the eigenspace $\cE_7$ contains a unique eigenfunction $S_V$, \resp a unique eigenfunction $C_V$, which is invariant, \resp anti-invariant, under the symmetry $\sigma_{i(V)}$ with respect to the median $[VM]$ issued from the vertex $V$, where $i(O)=1, i(A)=3$, and $i(B)=2$. More precisely, we choose
\begin{equation}\label{E23-sym4}
\begin{array}{l}
S_O = \Psi_{2,3}^{\frac{\pi}{2}} \text{~and~} C_O = \Psi_{2,3}^{0}\,,\\[5pt]
S_A = \Psi_{2,3}^{\frac{7\pi}{6}} \text{~and~} C_A = \Psi_{2,3}^{\frac{2\pi}{3}}\,,\\[5pt]
S_B = \Psi_{2,3}^{\frac{11\pi}{6}} \text{~and~} C_B = \Psi_{2,3}^{\frac{4\pi}{3}}\,.
\end{array}
\end{equation}

Note that $C_O$, \resp $S_O$, are the functions $C_{2,3}$, \resp $S_{2,3}$, and that,
\begin{equation}\label{E23-sym4a}
\begin{array}{l}
S_A = \rho_{-}^{*} S_O \text{~and~} C_A = \rho_{-}^{*} C_O\,,\\[5pt]
S_B = \rho_{+}^{*} S_O \text{~and~} C_B = \rho_{+}^{*} C_O\,.
\end{array}%
\end{equation}

The eigenfunctions $S_V$, \resp $C_V$, are permuted under the action of the rotations  $\rho_{+}=\sigma_2 \circ \sigma_1$ and $\rho_{-}=\sigma_1 \circ \sigma_2\,$. The eigenspace $\cE_7$ does not contain any non trivial rotation invariant eigenfunction.\medskip

Since $\Psi_{2,3}^{\theta+\pi} = - \Psi_{2,3}^{\theta}\,$, it follows from \eqref{E23-sym2} that, up to the symmetries $\sigma_i$, the nodal sets of the family $\Psi_{2,3}^{\theta}$, $\theta \in [0,2\pi]$ are determined by the nodal sets of the sub-family $\theta \in [0,\frac{\pi}{6}]$.\medskip

\textbf{From now on, we assume that $\theta \in [0,\frac{\pi}{6}]\,$.}

\subsection{Behaviour at the vertices}\label{S23-vert}

The vertices of the equilateral triangle $\cT$ belong to the nodal set $N(\Psi_{2,3}^{\theta})$ for all $\theta$. For geometric reasons, the order of vanishing at a vertex is at least $3$. More precisely,

\begin{properties}\label{S23-vert-P1}
Behaviour of $\Psi_{2,3}^{\theta}$ at the vertices.\vspace{-4mm}
\begin{enumerate}
  \item  The function $\Psi_{2,3}^{\theta}$ vanishes at order $6$ at $O$ if and only if $\theta \equiv 0 \pmod{\pi}$; otherwise it vanishes at order $3$.
  \item  The function $\Psi_{2,3}^{\theta}$ vanishes at order $6$ at $A$ if and only if $\theta \equiv \frac{2\pi}{3} \pmod{\pi}$; otherwise it vanishes at order $3$.
  \item  The function $\Psi_{2,3}^{\theta}$ vanishes at order $6$ at $B$ if and only if $\theta \equiv \frac{4\pi}{3} \pmod{\pi}$; otherwise it vanishes at order $3$.
\end{enumerate}\vspace{-4mm}
In other words, up to scaling, the only eigenfunction $\Psi_{2,3}^{\theta}$ which vanishes at higher order at a vertex $V \in \{O,A,B\}$ is $C_V$, the anti-invariant eigenfunction with respect to the median issued from the vertex $V$. In particular, when $\theta \in ]0,\frac{\pi}{6}]$, the three vertices are critical zeros of order three for the eigenfunction $\Psi_{2,3}^{\theta}$, and no interior nodal curve of such an eigenfunction can arrive at a vertex.
\end{properties}

\textbf{Proof}. Compute the Taylor expansions of $\Psi_{2,3}^{\theta}$ at the point  $(0,0)$, and use \eqref{E23-sym4a}.\hfill $\square$

\subsection{Fixed points on the medians}\label{S23-FP}

Since the median $[OM]$ is contained in the nodal set $N(C_{2,3})$, the intersection points of $[OM]$ with  $N(S_{2,3})$ are fixed points of the family $N(\Phi_{2,3}^{\theta})$, \ie, common zeros of the functions $\Psi_{2,3}^{\theta}\,$. If we parametrize $[OM]$ by $u \mapsto (u,u)$ with $u \in [0,\frac{1}{2}]$, we find that
\begin{equation}\label{FP23-2}
S_{2,3}\big|_{[OM]}(u) = 2\sin(10\pi u)-2\sin(6\pi u)-2\sin(4\pi u)\,.
\end{equation}
This function can be factored as,
\begin{equation}\label{FP23-4}
\begin{array}{ll}
S_{2,3}\big|_{[OM]}(u) & = -8\sin(2\pi u)\, \sin(3\pi u)\, \sin(5\pi u) \,.
\end{array}
\end{equation}

This formula shows that there are three fixed points on the open median $[OM]$, the centroid of the triangle $F_C = (\frac{1}{3},\frac{1}{3})_{\cF}$, the point $F_{1,O} = (\frac{1}{5},\frac{1}{5})_{\cF}$ and the point  $F_{2,O} = (\frac{2}{5},\frac{2}{5})_{\cF}\,$. The points $O$ and $M_O$ are ``obvious'' fixed points (corresponding to the values $0$ and $1/2$).

Taking into account the action of $G_{\cT}$ on the space $\cE_{7}$, see \eqref{E23-sym2}-\eqref{E23-sym4a}, we infer that the points $F_{1,A} = (\frac{7}{15},\frac{1}{3})_{\cF}\,$, $F_{2,A} = (\frac{4}{15},\frac{1}{3})_{\cF}\,$, $F_{1,B} = (\frac{1}{3},\frac{7}{15})_{\cF}\,$, and $F_{2,B} = (\frac{1}{3},\frac{4}{15})_{\cF}\,$,  are also common zeros for the family $\Psi^{\theta}_{2,3}\,$. They are deduced from $F_{i,O}$ by applying the rotations $\rho_{\pm}$, and situated on the two other medians.

Using Taylor expansions, it is easy to check that the fixed points $F_{*}$ are not critical zeros of the functions $\Psi_{2,3}^{\theta}$. In the neighborhood of the fixed points $F_{*}$, the nodal set consists of a single regular arc. \medskip 

\textbf{Remarks}. Note that we do not claim to have determined all the fixed points of the family $N(\Psi_{2,3}^{\theta})$ \ie, the set $N(C_{2,3}) \cap N(S_{2,3})$. We have only determined the fixed points located on the medians.
 Note also that the points $F_{1,B}, F_{2,O}$ and $F_{1,A}$ belong to the same line parallel to $[BA]$, and that the line through $F_{2,A}$ and $F_{2,B}$ is parallel to $[BA]\,$.

\begin{figure}
  \begin{minipage}[b]{0.45\linewidth}
    \centering
    \includegraphics[width=1\linewidth]{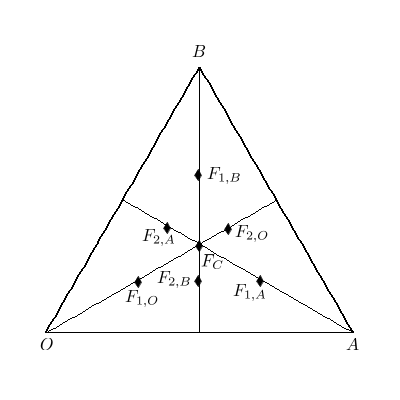}
    \par\vspace{0pt}
  \end{minipage}%
  \begin{minipage}[b]{0.45\linewidth}
    \centering
\begin{tabular}[b]{|c|c|c|}
\hline
Point & $\E^2$ & $\cF$\\[5pt] \hline
Vertex $O$ & $(0,0)$ & $(0,0)_{\cF}$\\[5pt] \hline
Vertex $A$ & $(1,0)$ & $(\frac{2}{3},\frac{1}{3})_{\cF}$\\[5pt] \hline
Vertex $B$ & $(\frac{1}{2},\frac{\sqrt{3}}{2})$ & $(\frac{1}{3},\frac{2}{3})_{\cF}$\\[5pt] \hline
Fixed point $F_C$ & $(\frac{1}{2},\frac{\sqrt{3}}{6})$ & $(\frac{1}{3},\frac{1}{3})_{\cF}$\\[5pt] \hline
Fixed point $F_{1,O}$ & $(\frac{3}{10},\frac{\sqrt{3}}{10})$ & $(\frac{1}{5},\frac{1}{5})_{\cF}$\\[5pt] \hline
Fixed point $F_{2,O}$ & $(\frac{3}{5},\frac{\sqrt{3}}{5})$ & $(\frac{2}{5},\frac{2}{5})_{\cF}$\\[5pt] \hline
Fixed point $F_{1,A}$ & $(\frac{7}{10},\frac{\sqrt{3}}{10})$ & $(\frac{7}{15},\frac{1}{3})_{\cF}$\\[5pt] \hline
Fixed point $F_{2,A}$ & $(\frac{2}{5},\frac{\sqrt{3}}{5})$ & $(\frac{4}{15},\frac{1}{3})_{\cF}$\\[5pt] \hline
Fixed point $F_{1,B}$ & $(\frac{1}{2},\frac{3\sqrt{3}}{10})$ & $(\frac{1}{3},\frac{7}{15})_{\cF}$\\[5pt] \hline
Fixed point $F_{2,B}$ & $(\frac{1}{2},\frac{\sqrt{3}}{10})$ & $(\frac{1}{3},\frac{4}{15})_{\cF}$\\[5pt] \hline
\end{tabular}
    \par\vspace{0pt}
    \end{minipage}
\caption{Fixed points for $N(\Psi_{2,3}^{\theta})$} \label{FP23-F1}
\end{figure}

\subsection{Partial barriers for the nodal sets}\label{S23-MB}

We have seen that the family $\Psi_{2,3}^{\theta}$ has  (at least) seven fixed points, the centroid $F_C$ of the triangle, and six other points $F_{1,O}\,$, $F_{2,O}\,$, $F_{1,A}\,$, $F_{2,A}\,$, $F_{1,B}$ and $F_{2,B}\,$, located respectively on the open medians $[OM]$, $[AM]$, and $[BM]$, see Figure~\ref{FP23-F1}.  More precisely, we have

\begin{lemma}\label{MB23-L1}
For any $\theta$, the nodal set $N(\Psi_{2,3}^{\theta})$ intersects each median at exactly three points unless the function $\Psi_{2,3}^{\theta}$ is one of the functions $C_V$ for $V \in \{O,A,B\}$, in which case the corresponding median is contained in the nodal set.  In particular, if $\theta \in ]0,\frac{\pi}{6}]$, the nodal set $N(\Psi_{2,3}^{\theta})$ only meets the medians at the fixed points $\{F_{i,O}, F_{i,A}, F_{i,B}, F_C\}$, for $i=1, 2$.
\end{lemma}

\textbf{Proof}. Use the following facts:\\ (i) the families $\{C_O,S_O\}$, $\{C_A,S_A\}$ and $\{C_B,S_B\}$ span $\cE_7$; \\ (ii) the function $C_V$ vanishes on the median issued from the vertex $V$.\\ Then write $\Psi_{2,3}^{\theta} = \alpha C_V + \beta S_V$. If $x \in [VM] \cap N(\Psi_{2,3}^{\theta})$, then $\beta S_V(x)=0$. If $\beta \not = 0$, then $x \in \{F_C, F_{1,V}, F_{2,V}\}$. If $\beta=0$, then $[VM] \subset N(\Psi_{2,3}^{\theta})$. \hfill $\square$\medskip

\textbf{Remark}. A consequence of Lemma~\ref{MB23-L1} and Subsection~\ref{S23-FP} is that for if $\theta \in ]0,\frac{\pi}{6}]$, no critical zero of the function $\Psi_{2,3}^{\theta}$ can occur on the medians.\medskip

For $0 \le a \le 1$, let $D_{a}$ denote the line whose equation in the parametrization $\cF$ is $s+t=a$.

We now investigate the intersections of the nodal sets with the medians, and with the lines $D_{a}$ for $a \in \{ \frac{2}{5}, \frac{3}{5}, \frac{2}{3}, \frac{4}{5}\}$ \ie, the lines through the fixed points, parallel to the edge $[BA]$.

\begin{lemma}\label{MB23-L2}
The intersections of the lines $D_{a}$ for $a \in \{\frac{2}{5}, \frac{3}{5}, \frac{2}{3}, \frac{4}{5}\}$ with the nodal sets $N(C_{2,3})$ and $N(S_{2,3})$ are as follows,
\begin{equation}\label{MB23-4}
\begin{array}{lll}
D_{\frac{2}{5}} \cap N(C_{2,3}) = \{F_{1,O}\}, &\text{~and~} & D_{\frac{2}{5}} \cap N(S_{2,3}) = \{F_{1,O}\},\\
D_{\frac{3}{5}} \cap N(C_{2,3}) = \{G_{O},F_{2,A},F_{2,B}\}, &\text{~and~} & D_{\frac{3}{5}} \cap N(S_{2,3}) = \{F_{2,A},F_{2,B}\},\\
D_{\frac{2}{3}} \cap N(C_{2,3}) = \{F_C, G_A,G_B\}, &\text{~and~} & D_{\frac{2}{3}} \cap N(S_{2,3}) = \{F_{C}\},\\
D_{\frac{4}{5}} \cap N(C_{2,3}) = \{F_{2,O},F_{1,A},F_{1,B}\}, &\text{~and~} & D_{\frac{4}{5}} \cap N(S_{2,3}) = \{F_{2,O},F_{1,A},F_{1,B}\},
\end{array}
\end{equation}
where $G_A$ and $G_B$ are symmetric with respect to $[OM]$, and $G_O =(0.3,0.3)_{\cF}$. These lines are tangent to $N(S_{2,3})$ at the points $F_{1,O}\,, F_{C}$ and $F_{2,O}\,$.
\end{lemma}%

\textbf{Proof}. The segment $D_{a}\cap \cT$ is parametrized by $u \mapsto (u,a-u)$ for $u \in [\frac{a}{3},\frac{2a}{3}]$. For each value $a \in \{\frac{2}{5}, \frac{3}{5}, \frac{2}{3}, \frac{4}{5}\}$, define
\begin{equation}\label{MB23-6}
\begin{array}{l}
BC_{a}(u) = C_{2,3}(u,a-u), ~~ u \in [\frac{a}{3},\frac{2a}{3}]\,,\\[5pt]
BS_{a}(u) = S_{2,3}(u,a-u), ~~ u \in [\frac{a}{3},\frac{2a}{3}]\,.
\end{array}
\end{equation}

Taking $a = \frac{2}{5}$, we find
\begin{equation}\label{MB23-6-25a}
\begin{array}{l}
BC_{\frac{2}{5}}(u) = -\cos(14 \pi u) + \cos(14\pi u + 2\pi/5) +
\cos(16\pi u) - \cos(16\pi u +3\pi/5)\,,\\
BS_{\frac{2}{5}}(u) = \sin(14 \pi u) - \sin(14\pi u + 2\pi/5) +
\sin(16\pi u) + \sin(16\pi u - 2\pi/5)\,.
\end{array}
\end{equation}

These equations simplify to

\begin{equation}\label{MB23-6-25b}
\begin{array}{l}
BC_{\frac{2}{5}}(u) =  -4 \sin(\pi/5) \sin(15\pi u) \sin(\pi u + 3\pi/10)\,,\\
BS_{\frac{2}{5}}(u) = 4 \sin(\pi/5) \sin(15\pi u) \cos(\pi u + 3\pi/10)\,.
\end{array}
\end{equation}

Similarly, we obtain

\begin{equation}\label{MB23-6-35b}
\begin{array}{l}
BC_{\frac{3}{5}}(u) =  4 \sin(\pi/5) \sin(15\pi u) \cos(\pi u + \pi/5)\,,\\
BS_{\frac{3}{5}}(u) = 4 \sin(\pi/5) \sin(15\pi u) \sin(\pi u + \pi/5)\,.
\end{array}
\end{equation}

\begin{equation}\label{MB23-6-23b}
\begin{array}{l}
BC_{\frac{2}{3}}(u) =  2\sqrt{3}\sin(9\pi u) \cos(7\pi u -\pi/3)\,,\\
BS_{\frac{2}{3}}(u) = 2\sqrt{3}\sin(9\pi u) \sin(7\pi u -\pi/3)\,.
\end{array}
\end{equation}

\begin{equation}\label{MB23-6-45b}
\begin{array}{l}
BC_{\frac{4}{5}}(u) =  -4 \cos(\pi/10) \sin(15\pi u) \sin(\pi u +\pi/10)\,,\\
BS_{\frac{4}{5}}(u) = 4 \cos(\pi/10) \sin(15\pi u) \cos(\pi u +\pi/10)\,.
\end{array}
\end{equation}

Looking at the zeros of the above functions in the respective intervals, the lemma follows. \hfill $\square$

\begin{figure}
  \centering
  \includegraphics[width=12cm]{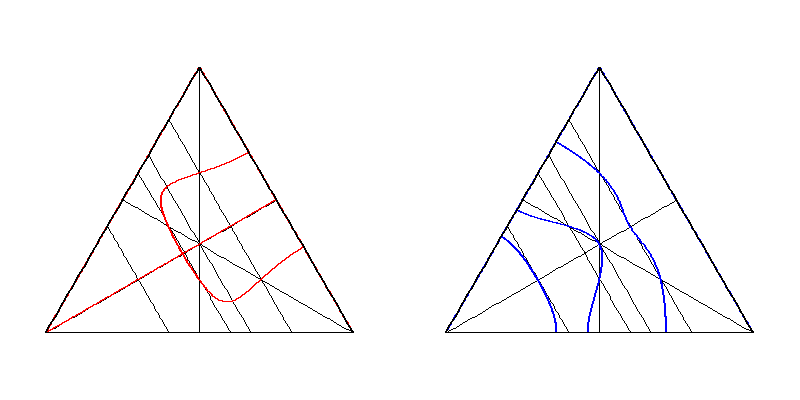}
  \caption{Partial barriers for $N(C_{2,3})$ and $N(S_{2,3})$}\label{MB23-F1}
\end{figure}

 The above lemmas are illustrated by Figure~\ref{MB23-F1}. It displays the partial barriers (thin segments), and the nodal sets as computed by Maple (thicker lines). 

\subsection{Critical zeros of $C_{2,3}$ and $S_{2,3}$ on the sides of $\cT$, \\and on the median $[OM]$}\label{S23-C}

Define the functions
\begin{equation}\label{23C-FCS}
\begin{array}{ll}
FC(u) &:= - 2\sin(8 \pi u) + 3 \sin(7\pi u) - 5 \sin(\pi u)\,,\\
FS(u) &:= - 2\cos(8 \pi u) - 3 \cos(7\pi u) + 5 \cos(\pi u)\,.\\
\end{array}
\end{equation}

\begin{properties}\label{23C-P1}
The partial derivatives of the functions $C_{2,3}$ and $S_{2,3}$ satisfy the fol\-lo\-wing relations.\vspace{-4mm}
\begin{enumerate}
\item Parametrize the edge $[OA]$ by $u \mapsto (u,u/2)$, with $u \in [0,2/3]$. Then,
\begin{equation}\label{23C-P1-OA}
\begin{array}{l}
\partial_s C_{2,3}(u,u/2) = 2 \pi FC(u),~~ \partial_t C_{2,3}(u,u/2) = -4 \pi FC(u)\,,\\
\partial_s S_{2,3}(u,u/2) = 2 \pi FS(u),~~ \partial_t S_{2,3}(u,u/2) = -4 \pi FS(u)\,.
\end{array}
\end{equation}
\item Parametrize the edge $[OB]$ by $u \mapsto (u/2,u)$, with $u \in [0,2/3]$. Then,
\begin{equation}\label{23C-P1-OB}
\begin{array}{l}
\partial_s C_{2,3}(u/2,u) = 4 \pi FC(u),~~ \partial_t C_{2,3}(u/2,u) = -2 \pi FC(u)\,,\\
\partial_s S_{2,3}(u/2,u) = -4 \pi FS(u),~~ \partial_t S_{2,3}(u/2,u) = 2 \pi FS(u)\,.
\end{array}
\end{equation}
\item Parametrize the edge $[BA]$ by $u \mapsto (u/2,1-u/2)$, with $u \in [2/3,4/3]$. Then,
\begin{equation}\label{23C-P1-BA}
\begin{array}{l}
\partial_s C_{2,3}(u/2,1-u/2) = -2 \pi FC(u),~~ \partial_t C_{2,3}(u/2,1-u/2) = -2 \pi FC(u)\,,\\
\partial_s S_{2,3}(u/2,1-u/2) = 2 \pi FS(u),~~ \partial_t S_{2,3}(u/2,1-u/2) = 2 \pi FS(u)\,.
\end{array}
\end{equation}
\end{enumerate}\vspace{-3mm}
As a consequence, the critical zeros of $C_{2,3}$ and $S_{2,3}$, on the edges $[OA]$ and $[OB]$, \resp on the edge $[BA]$, are determined by the zeros of $FC$ and $FS$ in the intervals $[0,2/3]$, \resp $[2/3,4/3]$.
\end{properties}%

\textbf{Proof}. It suffices to compute the partial derivatives of $C_{2,3}$ and $S_{2,3}$, and to make the substitutions corresponding to the parametrization of the edges. \hfill $\square$\medskip

Define the polynomials
\begin{equation}\label{23C-PCS}
\begin{array}{ll}
P_C(x) &:= 8x^3 + 2x^2 - 4x + 1\,,\\
P_S(x) &:= 8x^5 + 6x^4 - 10x^3 - 4x^2 + 4x - \frac{1}{4}\,.
\end{array}
\end{equation}

\begin{lemma}\label{23C-L1}
The functions $FC$ and $FS$ satisfy,
\begin{equation}\label{23C-FG}
\begin{array}{ll}
FC(u) & = - 8 \sin(\pi u) \left( \cos(\pi u) - 1\right)^2 \left( 2 \cos(\pi u) + 1\right)^2 \, P_C\left( \cos(\pi u)\right)\,,\\
FS(u) & = - 8 \left( \cos(\pi u) - 1\right) \left( 2 \cos(\pi u) + 1\right)^2 \, P_S\left( \cos(\pi u)\right)\,.
\end{array}
\end{equation}
\end{lemma}%

\textbf{Proof}. Use the Chebyshev polynomials. \hfill $\square$

\begin{lemma}\label{23C-L2}
Roots of the polynomials $P_C$ and $P_S$.\vspace{-3mm}
\begin{enumerate}
\item The polynomial $P_C$ has  exactly one root $-\xi_1 \in [-1,1]$, where
$\xi_1 \approx 0.9311441818\,$.
\item The polynomial $P_S$ has  exactly three roots in the interval $[-1,1]$,
    $$\eta_1 \approx 0.7261887036\,, ~\eta_2 \approx 0.5658979255\,, \text{~and~} \eta_3 \approx 0.06784981490\,.$$
\end{enumerate}
\end{lemma}%

 \textbf{Proof.} We leave the proof to the reader. \hfill $\square$ \medskip

Define the numbers
\begin{equation}\label{23C-GZ}
\begin{array}{l}
u_{1,C} := 1-\frac{\arccos{\xi_1}}{\pi} \approx 0.8811882234 \,,\\[5pt]
u_{2,C} := 1\,,\\[5pt]
u_{3,C} :=  1+\frac{\arccos{\xi_1}}{\pi} \approx 1.1188117766 \,,\\[5pt]
u_{1,S} := \frac{\arccos{\eta_1}}{\pi} \approx 0.2412898667 \,,\\[5pt]
u_{2,S} := \frac{\arccos{\eta_2}}{\pi} \approx 0.3085296215 \,,\\[5pt]
u_{3,S} := \frac{\arccos{\eta_3}}{\pi} \approx 0.4783861278 \,.\\[5pt]
\end{array}
\end{equation}

In the interval $[0,4/3]$, the function $FC$ vanishes at $0, 2/3$ and $4/3$ (these values of $u$ correspond to the vertices), and at the $u_{i,C}$ which correspond to critical zeros $Z_{i,C}$ of $C_{2,3}$ on the edge $[BA]$. In the interval $[0,4/3]$, the function $FS$ vanishes at $0, 2/3$ and $4/3$ (these values of $u$ correspond to the vertices), and at the $u_{i,S}$ which correspond to critical zeros $Z_{i,S}$ of $S_{2,3}$ on the edge $[OA]$, and $Z_{i+3,S}$ on the edge $[OB]$.

\begin{figure}
  \begin{minipage}[b]{0.44\linewidth}
    \centering
    \includegraphics[width=1\linewidth]{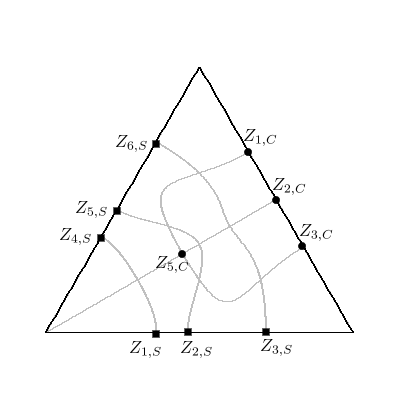}
    \par\vspace{0pt}
  \end{minipage}%
  \begin{minipage}[b]{0.55\linewidth}
    \centering
\begin{tabular}[b]{|c|c|c|}
\hline
Point & $\E^2$ & $\cF$\\[5pt] \hline
 $Z_{1,C}$ & $\approx (0.6609,0.5873)$ & $(\frac{u_{1,C}}{2},1-\frac{u_{1,C}}{2})_{\cF}$\\[5pt] \hline
$Z_{2,C}$ & $\approx (0.75,0.4330)$ & $(\frac{1}{2},\frac{1}{2})_{\cF}$\\[5pt] \hline
$Z_{3,C}$ & $\approx (0.8391,0.2787)$ & $(\frac{u_{3,C}}{2},1-\frac{u_{3,C}}{2})_{\cF}$\\[5pt] \hline
$Z_{5,C}$ & $\approx (0.4459,0.2576)$ & $(\frac{u_{5,C}}{2},\frac{u_{5,C}}{2})_{\cF}$\\[5pt] \hline
$Z_{1,S}$ & $\approx (0.3619,0)$ & $(u_{1,S},\frac{u_{1,S}}{2})_{\cF}$\\[5pt] \hline
$Z_{2,S}$ & $\approx (0.4628,0)$ & $(u_{2,S},\frac{u_{2,S}}{2})_{\cF}$\\[5pt] \hline
$Z_{3,S}$ & $\approx (0.7176,0)$ & $(u_{3,S},\frac{u_{3,S}}{2})_{\cF}$\\[5pt] \hline
$Z_{4,S}$ & $\approx (0.1810,0.3135)$ & $(\frac{u_{1,S}}{2},u_{1,S})_{\cF}$\\[5pt] \hline
$Z_{5,S}$ & $\approx (0.2314,0.4008)$ & $(\frac{u_{2,S}}{2},u_{2,S})_{\cF}$\\[5pt] \hline
$Z_{6,S}$ & $\approx (0.3588,0.6215)$ & $(\frac{u_{3,S}}{2},u_{3,S})_{\cF}$\\[5pt] \hline
\end{tabular}
    \par\vspace{0pt}
    \end{minipage}
\caption{Critical zeros of $C_{2,3}$ and $S_{2,3}$} \label{23C-F1}
\end{figure}

\begin{properties}\label{23C-P3}[Critical zeros of the functions $C_{2,3}$ and $S_{2,3}$ on the open edges of $\cT$, illustrated by Figure~\ref{23C-F1}] \vspace{-3mm}
\begin{enumerate}
\item The function $C_{2,3}$ has three critical zeros of order $2$ on the open edge $[BA]$: $Z_{1,C}$ between $B$ and $M_O$, $Z_{2,C} = M_O$, and $Z_{3,C}$ between $M_O$ and $A$. It has no critical zero on the open edges $[OA]$ and $[OB]$.
\item The function $S_{2,3}$ has three critical zeros of order $2$ on the open edge $[OA]$: $Z_{i,S}$, for $i \in \{1, 2, 3\}$. The points $Z_{1,S}$ and $Z_{2,S}$ lie between $O$ and $M_B$; the point $Z_{3,S}$ between $M_B$ and $A$. The function $S_{2,3}$ has three critical zeros of order $2$ on the open edge $[OB]$, $Z_{i+3,S}$ for $i\in \{1, 2, 3\}$, where $Z_{i,S}$ and $Z_{i+3,S}$ are symmetric with respect to the median $[OM]$. The function $S_{2,3}$ has no critical zero on the open edge $[BA]$.
\end{enumerate}
\end{properties}%

\textbf{Proof}. Use Properties~\ref{23C-P1} and Lemma~\ref{23C-L2}. \hfill $\square$ \medskip

\textbf{Remark}. The vertex $O$ is a critical zero of order $6$ of $C_{2,3}$, and a critical zero of order~$3$ of $S_{2,3}$. The vertices $A$ and $B$ are critical zeros of order $3$ of both $C_{2,3}$ and $S_{2,3}$, see Properties~\ref{S23-vert-P1}.

\begin{properties}\label{23C-P5} Critical zeros of the functions $C_{2,3}$ and $S_{2,3}$ on the median $[OM]$.\vspace{-3mm}
\begin{enumerate}
\item The function $C_{2,3}$ has one critical zero at $O$ of order $6$; two critical zeros of order~$2$, $M_O = Z_{2,C}$ and $Z_{5,C}$, corresponding to%
    $$u_{5,C} := 1-\arccos(1-1/\sqrt{2})/\pi \approx 0.5946180472\,.$$%
\item The function $S_{2,3}$ has no critical zero on the median $[OM]$, except the point $O$.
\end{enumerate}
\end{properties}%

\textbf{Proof}. Since $C_{2,3}$ vanishes on the median, its critical zeros on the median are the common zeros of its partial derivatives. They are precisely the zeros of the function
$$
\sin(5\pi u)-7\sin(3\pi u)+8\sin(2\pi u)\,,
$$
if we parametrize the median by $u \mapsto (u/2,u/2)$ for $u \in [0,1]$. The above function can be factorized as
$$
8\sin(\pi u) \big( \cos(\pi u)-1\big)^2 \big( 2\cos^2(\pi u)+4\cos(\pi u)+1\big)\,,
$$
and the first assertion follows.

For the second assertion, we have to look for the common zeros of the function $S_{2,3}$ and its derivatives on the median. This amounts to finding the common zeros of the functions
$$
\sin(5\pi u)-\sin(3\pi u)-\sin(2\pi u)\,,
$$
and
$$
5\cos(5\pi u)-3\cos(3\pi u)-2\cos(2\pi u)\,,
$$
in the interval $[0,1]$. The first function can be factorized as
$$
-4\sin(\frac{5\pi}{2}u)\sin(\frac{3\pi}{2}u)\sin(\pi u)\,,
$$

and it is easy to check that the only common zero is $u=0\,$. \hfill $\square$

Figure~\ref{23C-F1} displays the critical zeros of $C_{2,3}$ and $S_{2,3}$ (the nodal sets computed with Maple appear in grey).

In view of later reference, we mention the following corollary of Properties~\ref{23C-P5}.

Recall (see the notation \eqref{E23-sym4}) that $\Psi_{2,3}^{\frac{7\pi}{6}} = S_A\,$, and that $S_A = \rho_{-}^{*}S_O\,$.

\begin{corollary}\label{23C-C1}
Critical zeros of the function $S_A$. \vspace{-3mm}
\begin{enumerate}
  \item The function $S_A$ has three critical zeros of order $2$ on the open side $[OA]$,  $Z_{i,S_A} = \rho_{+}(Z_{i,S})$ for $i \in \{4,5,6\}$.
  \item The function $S_A$ has three critical zeros of order $2$ on the open side $[BA]\,$,  $Z_{i,S_A} = \rho_{+}(Z_{i,S})$ for $i \in \{1,2,3\}\,$.
  \item The function $S_A$ has no critical zero on the open side $[OB]$, and no critical zero on the open median $[AM_A]$.
\end{enumerate}
\end{corollary}%

\textbf{Remark}. With the usual parametrization of the edges, the critical zeros $Z_{i,S_A}$ are associated with the values $u_{i,S_A}$ defined by
\begin{equation}\label{23C-sym10}
\begin{array}{l}
u_{i,S_A} = \frac{4}{3} - u_{i,S}\,, \text{~for~} i \in \{1, 2, 3\}\,,\\[5pt]
u_{3+i,S_A} = \frac{2}{3} - u_{i,S}\,, \text{~for~} i \in \{1, 2, 3\}\,.
\end{array}
\end{equation}

\subsection{The nodal sets of $C_{2,3}$ and $S_{2,3}$}\label{S23-NCS}

\begin{properties}\label{NCS23-P2} Nodal sets of $C_{2,3}$ and $S_{2,3}$.\vspace{-4mm}
\begin{enumerate}
\item The function $C_{2,3}$ has only one critical zero $Z_{5,C}$ in the interior of the triangle $\cT$. Its nodal set consists of the diagonal $[OM_O]$, and an injective regular arc from $Z_{1,C}$ to $Z_{3,C}$ which intersects $[OM_O]$ orthogonally at $Z_{5,C}$, and passes through the fixed points $F_{1,A}$, $F_{2,A}$, $F_{1,B}$ and $F_{2,B}$.
\item The function $S_{2,3}$  has no critical zero in the interior of the triangle. Its nodal set consists of three disjoint injective regular arcs: one from $Z_{1,S}$ to $Z_{4,S}$, passing through $F_{1,O}$; one from $Z_{2,S}$ to $Z_{5,S}$, passing through $F_{2,A}$, $F_{2,B}$, and $F_C$; one from $Z_{3,C}$ to $Z_{6,C}$, passing through $F_{1,A}$, $F_{2,O}$, and $F_{1,B}$.
\end{enumerate}
\end{properties}%

\begin{figure}
\centering
\includegraphics[width=12cm]{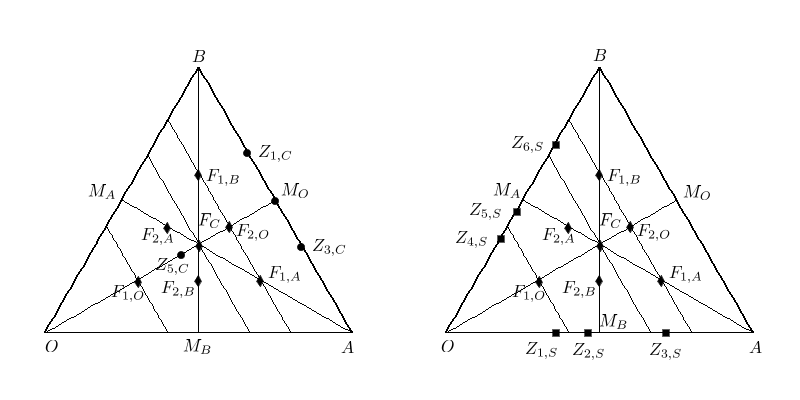}
\caption{Barriers, fixed points and critical zeros for $N(C_{2,3})$ and $N(S_{2,3})$}\label{NCS23-F1}
\end{figure}

\textbf{Proof}. We already know that the nodal set $N(C_{2,3})$ contains the median $[OM_O]$. Except for this segment, no other nodal arc for either $C_{2,3}$ or $S_{2,3}$ can arrive at a vertex. Figure~\ref{NCS23-F1} displays the fixed points and the critical zeros for $N(C_{2,3})$ (left) and $N(S_{2,3})$ (right). We know the local behaviour of the nodal sets near the fixed points $F_{*}$, and near the critical zeros $Z_{*}$.  Three  partial barriers (see Subsection~\ref{S23-MB}) are also shown in Figure~\ref{NCS23-F1}. The nodal sets $N(C_{2,3})$ and $N(S_{2,3})$  only meet the barriers at the points given in Lemma~\ref{MB23-L2}.

The medians divide $\cT$ into six isometric $H$-triangles. The nodal sets of $C_{2,3}$ and $S_{2,3}$ consist of finitely many nodal arcs which are smooth except at the critical zeros. They can only exit the interior of an $H$-triangle at a fixed point or at a critical zero. \medskip

\emph{Nodal set $N(C_{2,3})$.}\\
 It suffices to look at the subset $N':=\left( N(C_{2,3}) \setminus [OM_O] \right) \cup \{Z_{5,C}\}$.

First of all, notice that for each $H$-triangle determined by the medians of $\cT$, there are exactly two exit points for an arc belonging to $N'$, with only one exit direction at each point. Indeed, except for the edges contained in the median $[OM]$, each edge contains at most one exit point (a fixed point or a critical zero). Furthermore, the vertices of the $H$-triangles are not exit points.

We claim that the function $C_{2,3}$ cannot have any critical zero in the interiors of the $H$-triangles. Indeed, assume that there is one critical zero $Z \in N'$ in the interior of some $H$-triangle $\cH$. At this point, the nodal set would consist of at least four semi-arcs. Following any such semi-arc, we either obtain a simply closed nodal arc, or exit the triangle. Since there are at most two exit points, there would be at least one simply closed nodal component in the interior of the triangle $\cH$. This component would bound at least one nodal domain $\omega$. The first Dirichlet eigenvalue $\lambda(\omega)$ would satisfy $\lambda(\omega) = \lambda_7(\cT) = 19$. On the other-hand, since $\omega$ is contained in the interior of $\cH$, we would have $\lambda(\omega) > \lambda(\cH) = 21$ according to Lemma~\ref{MB13-23-L}, a contradiction. \medskip

The last argument in the proof of the claim also shows that the interiors of the $H$-triangles cannot contain any closed nodal component of $C_{2,3}$. This shows that the nodal set of $C_{2,3}$ is indeed as shown in Figure~\ref{ETR-F4}. \medskip

\emph{Nodal set $N(S_{2,3})$.}  First of all, notice that for the $H$-triangle determined by the medians of $\cT$, there are either two or four exit points, with only one exit direction at each point. More precisely, there are three cases.

Case (i). The triangles $T(F_C,B,M_O)$ and $T(F_C,A,M_O)$ have two exit points which are fixed points. The arguments used for $C_{2,3}$ apply for these triangles. \smallskip

Case (ii). The triangles $T(F_C,O,M_A)$ and $T(F_C,O,M_B)$ have four exit points, two fixed points and two critical zeros on the open edges. We claim that the function $S_{2,3}$ cannot have any critical zero in the interiors of these $H$-triangles. Indeed, assume that there is one critical zero $Z$ in the interior of such an $H$-triangle $\cH$. At this point, the nodal set would consist of at least four semi-arcs. Following any such semi-arc, we either obtain a simply closed nodal arc, or exit the triangle. The preceding arguments show that we cannot have any simply closed nodal component inside $\cH$. Since there are four exit points, each of them with a single exit direction, we would have two arcs joining the exit points and meeting at the critical zero $Z$. This would yield a nodal domain $\omega$ bounded by two semi-arcs and a segment in one of the edges of $\cT$, or at least one simply closed nodal component in the interior of the triangle $\cH$. The first Dirichlet eigenvalue $\lambda(\omega)$ would satisfy $\lambda(\omega) = \lambda_7(\cT) = 19$. On the other-hand, since $\omega$ is contained in the interior of $\cH$, we would have $\lambda(\omega) > \lambda(\cH) = 21$, according to Lemma~\ref{MB13-23-L},  a contradiction.  \smallskip

Case (iii). The triangles $T(F_C,B,M_A)$ and $T(F_C,A,M_B)$ have four exit points, the vertex $F_C$, two fixed points on the open edges, and one critical zero on an open edge. We begin with the same argument as in the previous case. To conclude, we use the barriers given by Lemma~\ref{MB23-L2}.

The above proofs also show that the interiors of the $H$-triangles cannot contain any closed nodal component of $S_{2,3}$. This shows that the nodal set of $S_{2,3}$ is indeed as shown in Figure~\ref{ETR-F4}.
\hfill $\square$

\subsection{Critical zeros of $\Psi_{2,3}^{\theta}$ on the sides of $\cT$, for $\theta \in ]0,\frac{\pi}{6}]$}\label{S23-CP}

As a consequence of Properties~\ref{23C-P1}, the critical zeros of the functions $\Psi_{2,3}^{\theta}$ on the sides of the triangle $\cT$ are determined by one of the equations
$$
\cos\theta \, FC(u) \pm \sin\theta \, FS(u) = 0\,.
$$

Recall that for $\theta \in ]0,\frac{\pi}{6}]$, the vertices of $\cT$ are critical zeros of order $3$ of $\Psi_{2,3}^{\theta}$ (Properties~\ref{S23-vert-P1}).

\begin{properties}\label{23CP-P2}
The critical zeros of the function $\Psi_{2,3}^{\theta}$ on the open edges of the triangle $\cT$ are determined by the following equations. \vspace{-3mm}
\begin{enumerate}
\item On the edge $[OA]$ parametrized by $u \mapsto (u,u/2)$,
  \begin{equation}\label{23C-P2-OA}
  \cos\theta \, FC(u) + \sin\theta \, FS(u) = 0\,, ~\text{~for~} u \in [0,2/3]\,.
  \end{equation}
\item On the edge $[OB]$ parametrized by $u \mapsto (u/2,u)$,
  \begin{equation}\label{23C-P2-OB}
  \cos\theta \, FC(u) - \sin\theta \, FS(u) = 0\,, ~\text{~for~} u \in [0,2/3]\,.
  \end{equation}
\item On the edge $[BA]$ parametrized by $u \mapsto (u/2,1-u/2)$,
  \begin{equation}\label{23C-P2-BA}
  \cos\theta \, FC(u) - \sin\theta \, FS(u) = 0\,, ~\text{~for~} u \in [2/3,4/3]\,.
  \end{equation}
\end{enumerate}
\end{properties}%

For convenience, we introduce the functions,
\begin{equation}\label{23C-K0}
K_{\pm}^{\theta}(u) = \cos\theta \, FC(u) \pm \sin\theta \, FS(u)\,.
\end{equation}

\begin{figure}[ht]
\centering
\includegraphics[width=12cm]{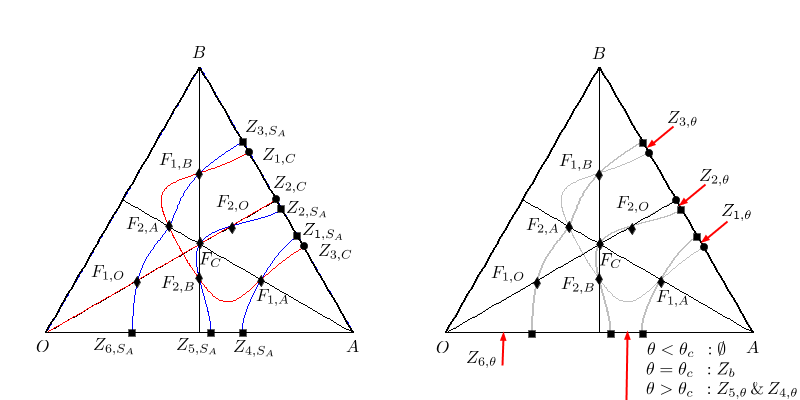}
\caption{Boundary critical zeros for $C$, $S_A\,$, and $\Psi_{2,3}^{\theta}$} \label{NSP23-F1}
\end{figure}

\begin{properties}\label{23CP-P3}
[ Illustrated by Figure~\ref{NSP23-F1} and using the notation \eqref{23C-GZ}]\vspace{-3mm}
\begin{enumerate}
\item Zeros of $K_{+}^{\theta}$ in the interval $]0,2/3[$ (corresponding to critical zeros on the open side $[OA]$). There exists $u_{b} \in\, ]\frac{1}{3},u_{3,S}[$, and $\theta_c \in\, ]0,\frac{\pi}{6}[$ such that:
\begin{enumerate}
  \item if $0 < \theta < \theta_c\,$, the function $K_{+}^{\theta}$ has only one simple zero $\alpha_1(\theta) \in \, ]0,u_{6,S_A}]$\,;
  \item if $\theta = \theta_c\,$, the function $K_{+}^{\theta}$ has a simple zero  $\alpha_1(\theta_c) \in\,  ]0,u_{6,S_A}]$, and a double zero at $u_{b}$\,.
  \item if $\theta_c < \theta \le \frac{\pi}{6}\,$, the function $K_{+}^{\theta}$ has three simple zeros, $\alpha_1(\theta) \in\, ]0,u_{6,S_A}]$, $\alpha_2(\theta) \in [u_{5,S_A},u_{b}[$, and $\alpha_3(\theta) \in\, ]u_{b},u_{4,S_A}]$.
\end{enumerate}
    The function $\alpha_1(\theta)$ and $\alpha_3(\theta)$ are increasing. The function $\alpha_2(\theta)$ is decreasing.
  \item Zeros of $K_{-}^{\theta}$ in the interval $]0,2/3[$ (corresponding to critical zeros on the open side $[OB]$). In this interval, the function $K_{-}^{\theta}$ does not vanish.
  \item Zeros of $K_{-}^{\theta}$ in the interval $]2/3,4/3[$ (corresponding to critical zeros on the open side $[BA]$).The function $ K_{-}^{\theta}$ has three simple zeros $\omega_1(\theta) \in [u_{3,S_A},u_{1,C}[\,$,  $\omega_2(\theta) \in ]1,u_{2,S_A}]\,$, and $\omega_3(\theta) \in [u_{1,S_A},u_{3,C}[\,$. The functions $\omega_1$ and $\omega_3$ are decreasing, the function $\omega_2$ is increasing.
\end{enumerate}
\end{properties}%

\textbf{Proof}. Notice that the zeros are continuous with respect to $\theta$ (indeed, the equations can be transformed into polynomials  whose coefficients are continuous in $\theta$). An information on the possible location of zeros in given by Table~\ref{23CP-T1}, depending on the value of $\theta$, see lines 7 and 9 in the table. \medskip

We first investigate whether the zeros of $K_{\pm}^{\theta}$ can have order at least $2$. More precisely, we investigate whether there exists a pair $(\theta,u)$ such that
\begin{equation}\label{23C-H1}
\begin{array}{l}
\cos\theta \, FC(u) \pm \sin\theta \, FS(u) = 0\,,\\
\cos\theta \, FC'(u) \pm \sin\theta \, FS'(u) = 0\,.
\end{array}
\end{equation}

For this purpose, we define the function
\begin{equation}\label{23C-HW1}
WFCS (u) := FC(u) FS'(u) - FS(u) FC'(u).
\end{equation}

\begin{lemma}\label{23C-LW2}
Define the polynomial
\begin{equation}\label{23C-PW}
P_W(x) := - 6x^5 + 25x^3 - 15 x^2 - 15x + 11\,.
\end{equation}
Then,
\begin{equation}\label{23C-W2}
\begin{array}{lll}
\frac{1}{2\pi}WFCS(u) & = & 28 + 3\sin(8\pi u)\sin(7\pi u)
- 3\cos(8\pi u)\cos(7\pi u)\\
&&~~ - 35\sin(8\pi u)\sin(\pi u) + 35\cos(8\pi u)\cos(\pi u)\\
&&~~ - 60\sin(7\pi u)\sin(\pi u) - 60\cos(7\pi u)\cos(\pi u) \\[5pt]
& = & 28 - 3\cos(15\pi u) + 35\cos(9\pi u) - 60\cos(6\pi u) \\[5pt]
&=& 8 P_W\left( \cos(3\pi u)\right)\,.
\end{array}
\end{equation}
The polynomial $P_W$ factors as
\begin{equation}\label{23C-W4}
P_W(x) = -(x-1)^3(x+\xi)(x+\eta)\,,
\end{equation}
where
\begin{equation}\label{23C-W6}
\xi = \frac{9-\sqrt{15}}{6} \text{~and~} \eta = \frac{9+\sqrt{15}}{6}\,.
\end{equation}
Let
\begin{equation}\label{23C-W8}
u_{0,W} := \frac{1}{3\pi} \arccos(\xi) \approx 0.2753793461\,.
\end{equation}
In the interval $]-1/6,3/2[$, the function $WFCS$ vanishes at the points $u \in \{0, \frac{2}{3}, \frac{4}{3}\}$ (which correspond to vertices), and at the points $\{u_{1,W},u_{2,W},u_{3,W},u_{4,W}\}$,
where
\begin{equation}\label{23C-W10}
\begin{array}{lll}
u_{1,W} &=  \frac{1}{3} - u_{0,W} &\approx 0.2753793461\,,\\
u_{2,W} &=  \frac{1}{3} + u_{0,W} &\approx 0.3912873205\,,\\
u_{3,W} &= 1 - u_{0,W} &\approx 0.9420460128\,,\\
u_{4,W} &= 1 + u_{0,W} &\approx 1.057953987\,.\\
\end{array}
\end{equation}
The points $u_{i,W}$, are simple zeros of the function $WFCS$.
\end{lemma}%

\emph{Proof of the lemma}. \\
Equations \eqref{23C-W2} follow by computing the derivatives of $FC$ and $FS$, by expanding the expression of $WFCS$, and by making use of the Chebyshev polynomials. The remaining part of the lemma follows easily.\hfill $\square$

Recall that $\theta \in \, ]0,\frac{\pi}{6}]$. If $WFCS(u) \not = 0$, the system \eqref{23C-H1} has no solution $(\theta,u)$. If $(\theta,u)$ is a solution of the system \eqref{23C-H1}, then $WFCS(u) = 0$, and $\cos\theta FC(u) \pm \sin\theta FS(u)=0$, which implies that $\pm \frac{FC(u)}{FS(u)} \in\, ]0,\sqrt{3}]$. Computing this ratio for the above values $u_{i,W}$, we conclude that $u_{2,W}$ is the unique value for which the condition can be satisfied, and this only occurs for the function $K_{+}^{\theta}$.

Define
\begin{equation}\label{23C-W12}
\begin{array}{l}
u_{b} := u_{2,W} \approx 0.3912873205\,,\\[5pt]
\theta_{c} := \arctan\left( \frac{FC(u_{b})}{FS(u_{b})} \right) \approx 0.3005211736\,.
\end{array}
\end{equation}

What we have just proved is that, for all $\theta \in\, ]0, \frac{\pi}{6}]$, $K_{-}^{\theta}$ has only simple zeros in the set $]0,2/3[ \,\cup\,  ]2/3,4/3[\,$; for all $\theta \in ]0, \frac{\pi}{6}]\setminus \{\theta_c\}\,$, $K_{+}^{\theta}$ has only simple zeros in the set $]0,2/3[\, \cup\, ]2/3,4/3[\,$; for $\theta = \theta_c$, $K_{+}^{\theta_c}$ has exactly one double zero at $u_{b}$. This corresponds to a critical zero of order $3$ on the edge $[OA]$.

Lemma~\ref{23C-LW2} implies that for $\theta \in ]0, \pi/6]\setminus \{\theta_c\}$, the zeros (if any) of the functions $K_{\pm}^{\theta}$ in the interval $]0,\frac{2}{3}[$ and $]\frac{2}{3},\frac{4}{3}[$ are simple, so that they are smooth in $\theta$. If $u(\theta)$ is such a zero, its derivative with respect to $\theta$ satisfies the relation
$$
1+\tan^2(\theta) = \pm \frac{WCS(u(\theta))}{GS^2(u(\theta))} u'(\theta)\,.
$$

We can now start from the function $\Psi_{2,3}^{\frac{\pi}{6}} = -S_A$, and follow the zeros by continuity, using Corollary~\ref{23C-C1}. This proves Properties~\ref{23CP-P3} on the open edge $[OA]$ when $\theta_c < \theta < \frac{\pi}{6}$, and on the open edges $[OB]$ and $[BA]$ for all $\theta$. When $\theta=\theta_c$, we have a critical zero of order $2$ and a critical zero of order $3$ on $[OA]$. It is easy to see that for $\theta > 0$ very small, there is only one critical zero on the open edge $[OA]$, and we can follow this zero by continuity for $0 < \theta < \theta_c\,$. \hfill $\square$\medskip

\begin{table}[hbt]\scriptsize
\centering
\begin{tabular}{|c|c|c|c|c|c|}
\hline
$u \in [0,4/3]$ & $\approx$ & $GC$ & $GS$ & $K_{+}^{\theta}$ & $K_{-}^{\theta}$ \\ \hline
$0$ & $0$ & $0$ & $\approx 3.7500$ & $\approx 3.7500 \sin\theta$ & $\approx -3.7500 \sin\theta$ \\ \hline
{\clr $\alpha_1(\theta)$} & -- & -- & -- & {\clr $0$} & -- \\ \hline
$u_{1,S}$ & $\approx 0.2413$ & $\approx -0.4167$ & $\approx 0$ & $\approx -0.4167 \cos\theta$ & $\approx -0.4167 \cos\theta$ \\ \hline
$u_{2,S}$ & $\approx 0.3085$ & $\approx -0.2959$ & $\approx 0$ & $\approx -0.2959 \cos\theta$ & $\approx -0.2959 \cos\theta$ \\ \hline
$\frac{1}{3}$ & $\approx 0.3333$ & $\approx -0.2165$ & $\approx 0.1250$ & $
\approx -0.2165 \cos\theta + 0.1250 \sin\theta$ & $\approx -0.2165 \cos\theta - 0.1250 \sin\theta$ \\ \hline
{\clr $\alpha_2(\theta)$} & -- & -- & -- & {\clr $0$ (if $\theta \ge \theta_{c}$)} & -- \\ \hline
$u_b$ & $\approx 0.3913$ & $\approx -0.1161$ & $\approx 0.3745$
& $\approx -0.1161 \cos\theta + 0.3745 \sin\theta$ & $\approx -0.1161 \cos\theta - 0.3745 \sin\theta$ \\ \hline
{\clr $\alpha_3(\theta)$} & -- & -- & -- & {\clr $0$ (if $\theta \ge \theta_{c}$)}& -- \\ \hline
$u_{3,S}$ & $\approx 0.4784$ & $\approx -.6885$ & $\approx 0$ & $\approx -0.6885 \cos\theta$ & $\approx-0.6885 \cos\theta$ \\ \hline
$\frac{2}{3}$ & $\approx 0.6667$ & $\approx -3.2476$ & $\approx -1.8750$
& $\approx -3.2476 \cos\theta - 1.8750 \sin\theta$ & $\approx -3.2476 \cos\theta + 1.8750 \sin\theta$ \\ \hline
$u_{3,SB}$ & $\approx 0.8549$ & $\approx -0.3442$ & $\approx -0.5962$
& $\approx -0.3442 \cos\theta - 0.5962 \sin\theta$ & $\approx -0.3442 \cos\theta + 0.5962 \sin\theta$ \\ \hline
{\clr $\omega_1(\theta)$} & -- & -- & -- & -- & {\clr $0$} \\ \hline
$u_{1,C}$ & $\approx 0.8812$ & $\approx 0$ & $\approx -.4588$ & $\approx - 0.4588 \sin\theta$
& $\approx 0.4588 \sin\theta$ \\ \hline
$u_{2,C}$ & $1$ & $\approx0$ & $\approx -.2500$ & $\approx- 0.2500 \sin\theta$
& $\approx 0.2500 \sin\theta$ \\ \hline
{\clr $\omega_2(\theta)$} & -- & -- & -- & -- & {\clr $0$} \\ \hline
$u_{2,SB}$ & $\approx 1.02480$ & $\approx -0.1479$ & $\approx -0.2562$
& $\approx -0.1479 \cos\theta - 0.2562 \sin\theta$ & $\approx -0.1479 \cos\theta + 0.2562 \sin\theta$ \\ \hline
$u_{1,SB}$ & $\approx 1.0920$ & $\approx -0.2083$ & $\approx -0.3609$
& $\approx -0.2083 \cos\theta - 0.3609 \sin\theta$ & $\approx -0.2083 \cos\theta + 0.3609 \sin\theta$ \\ \hline
{\clr $\omega_3(\theta)$} & -- & -- & -- & -- & {\clr $0$} \\ \hline
$u_{3,C}$ & $\approx 1.1188$ & $\approx 0$ & $\approx -0.4588$ & $\approx - 0.4588 \sin\theta$ & $\approx 0.4588 \sin\theta$ \\ \hline
$\frac{4}{3}$ & $\approx 1.3333$ & $\approx 3.2476$ & $\approx -1.8750$
& $\approx 3.2476 \cos\theta - 1.8750 \sin\theta$ & $\approx 3.2476 \cos\theta + 1.8750 \sin\theta$ \\ \hline
\end{tabular}
\caption{Values of $K_{\pm}^{\theta}\,$}\label{23CP-T1}
\end{table}

\begin{corollary}\label{23CP-C1}
Critical zeros of the function $\Psi_{2,3}^{\theta}$ for $\theta \in ]0,\frac{\pi}{6}]$.\vspace{-3mm}
\begin{enumerate}
  \item The function $\Psi_{2,3}^{\theta}$ always has  unique  critical zero of order $2$ $Z_{6,\theta} \in\, ]O,Z_{6,S_A}]$.\vspace{-3mm}
  \begin{itemize}
  \item If $\theta < \theta_c\,$, the function $\Psi_{2,3}^{\theta}$ has no other critical zero on the open edge $[OA]$.
  \item If $\theta=\theta_c\,$, the function $\Psi_{2,3}^{\theta}$ has  a unique critical zero $Z_b$ of order $3$ in the interval $]Z_{5,S_A},Z_{6,S_A}[$.
  \item
   If $\theta > \theta_c\,$, the function $\Psi_{2,3}^{\theta}$ has   exactly  two critical zeros of order $2$,  $Z_{5,\theta} \in [Z_{5,S_A},Z_b[$ and $Z_{4,\theta} \in \,]Z_b,Z_{4,\theta}]\,$.
   \end{itemize}
  \item The function $\Psi_{2,3}^{\theta}$ has three critical zeros of order $2$ on the open edge $[BA]$:  $Z_{3,\theta} \in [Z_{3,S_A},Z_{1,C}[\,$, $Z_{2,\theta} \in\, ]Z_{2,C},Z_{2,S_A}]\,$, and
      $Z_{1,\theta} \in [Z_{1,S_A},Z_{3,C}[\,$.
  \item The function $\Psi_{2,3}^{\theta}$ has no critical zero on the open edge $[OB]$.
\end{enumerate}
\end{corollary}%

The corollary is illustrated by Figure~\ref{NSP23-F1}: the picture on the left shows the fixed points, the boundary critical zeros of $C_{2,3}$ and $S_A$, as well as the nodal sets $N(C_{2,3})$ and $N(S_A)$ (Maple calculation). The picture on the right shows the fixed points and the boundary critical zeros for $\Psi_{2,3}^{\theta}$.

\begin{figure}[ht!]
\centering
\includegraphics[width=12cm]{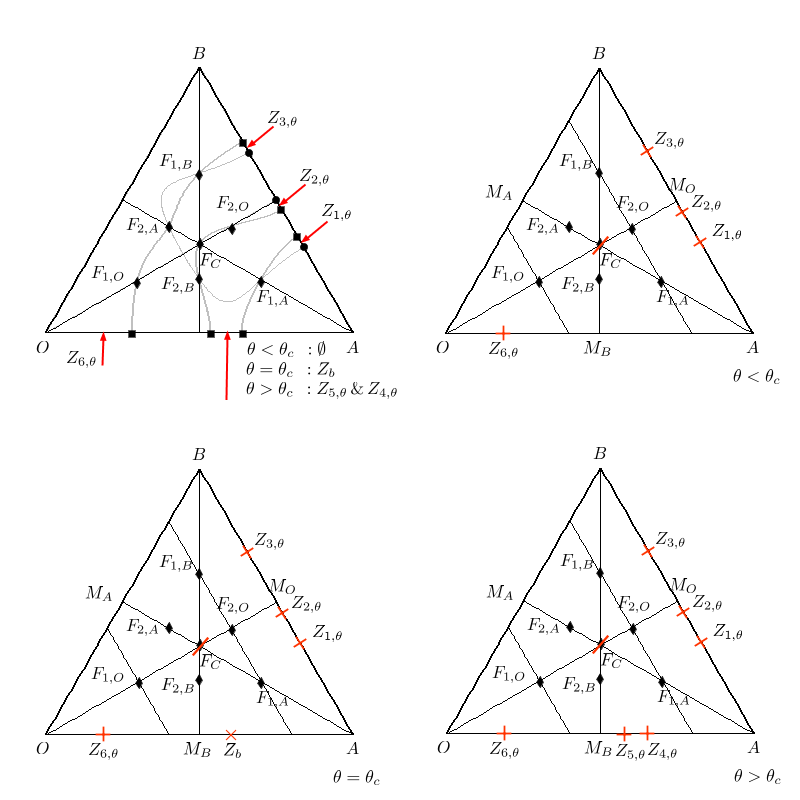}
\caption{Boundary critical zeros for $\Psi_{2,3}^{\theta}$} \label{NSP23-F2}
\end{figure}

\newpage

\subsection{Nodal set of $\Psi_{2,3}^{\theta}$}\label{S23-NSP}

\begin{figure}[htb!]
  \centering
  \includegraphics[width=12cm]{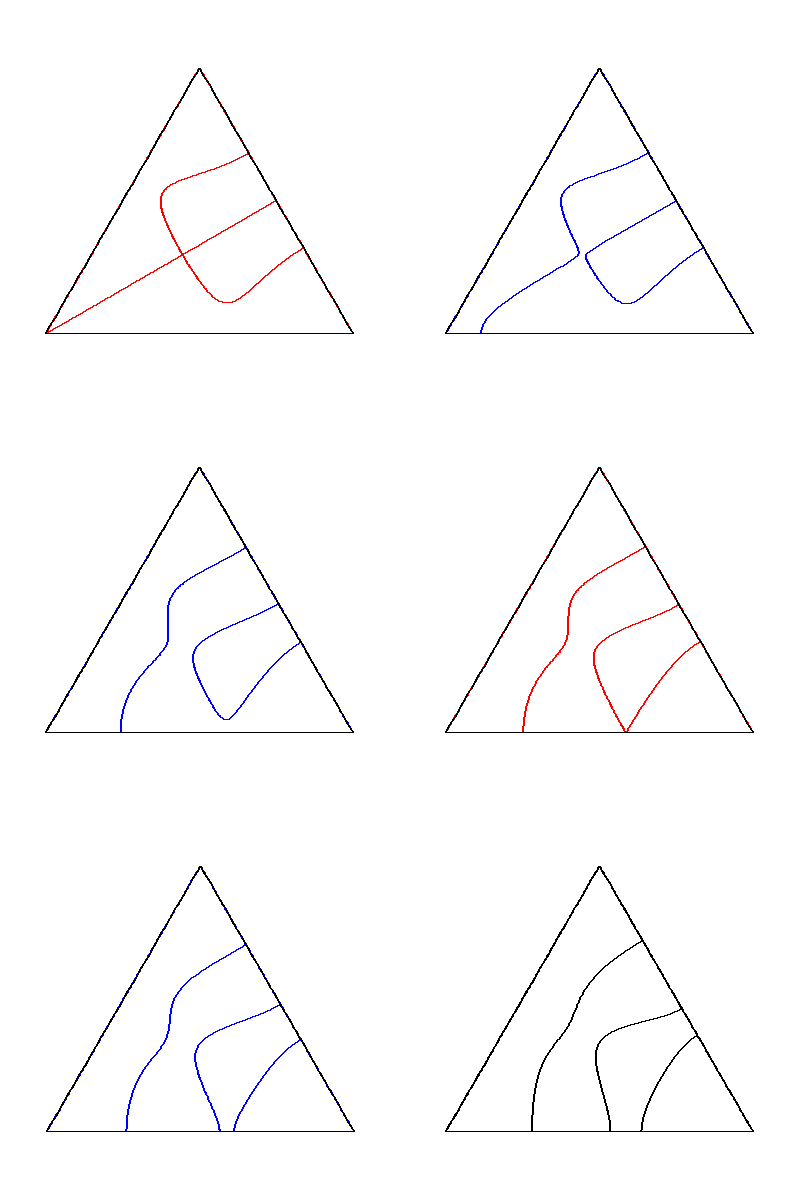}
  \caption{ Bifurcations for $N(\Psi_{2,3}^{\theta}), \theta \in ]0, \frac{\pi}{6}]$}\label{NSP23-F3}
\end{figure}

\begin{proposition}\label{NSP23-P1}[Illustrated by Figure~\ref{NSP23-F3}]\vspace{-3mm}
\begin{enumerate}
\item For $\theta \in ]0,\theta_c[$ the nodal set of $\Psi_{2,3}^{\theta}$ consists of two disjoint injective arcs, one from $Z_{6,\theta}$ to $Z_{3,\theta}$,
through $F_{1,O}$, $F_{2,A}$, and $F_{1,B}$, another from $Z_{1,\theta}$ to $Z_{2,\theta,}$, through the points $\{F_{1,A}$, $F_{2,B}$, $F_{C}$ and $F_{2,O}\}$.
 In particular, the function $\Psi_{2,3}^{\theta}$ has three nodal domains.
\item For $\theta = \theta_c$, the nodal set of $\Psi_{2,3}^{\theta}$ consists of three disjoint injective arcs, one from $Z_{6,\theta}$ to $Z_{3,\theta}$ through the points $F_{1,O}$ $F_{2,A}$, and $F_{1,B}$; one from $Z_b$ to $Z_{2,\theta}$, through the points $F_{2,B}$, $F_C$, and $F_{2,O}$; one from $Z_b$ to $Z_{1,\theta}$, through the point $F_{1,A}$. In particular, the eigenfunction has  four nodal domains.
\item For $\theta_c < \theta \le \frac{\pi}{6}$, the nodal set of $\Psi_{2,3}^{\theta}$ consists of three disjoint injective arcs, one from $Z_{6,\theta}$ to $Z_{3,\theta}$ through the points $F_{1,O}$ $F_{2,A}$, and $F_{1,B}$; one from $Z_{5,\theta}$ to $Z_{2,\theta}$, through the points $F_{2,B}$, $F_C$, and $F_{2,O}$; one from $Z_{4,\theta}$ to $Z_{1,\theta}$, through the point $F_{1,A}$.
 In particular, the eigenfunction has  four nodal domains.
\end{enumerate}
As a consequence,  eigenfunctions associated with the eigenvalue $\lambda_7(\cT)=17$  have at most $4$ nodal domains, so that this eigenvalue  is not Courant-sharp.
\end{proposition}%

\textbf{Proof}. We determine the nodal set of $\Psi_{2,3}^{\theta}$ in each of the six $H$-triangles determined by the medians of the triangle $\cT$.

We first observe that the nodal set has a tangent at the point $F_C$ which makes an angle less than $\frac{\pi}{6}$ with $[OM]$. This implies that the point $F_C$ is an exit point for two $H$-triangles only, namely $T(B,F_C,M_O)$ and $T(O,F_C,M_B)$.

For $a \in [0,1]$, we consider the functions
\begin{equation}\label{23C-20}
BP_{a}^{\theta}(u) = \Psi_{2,3}^{\theta}(u,a-u)\,,
\end{equation}
for $u \in [\frac{a}{3},\frac{2a}{3}]$ \ie, the restrictions of the functions $\Psi_{2,3}^{\theta}$ to the segments $[B_aA_a] = D_a \cap \cT$, whose end points $B_a$ and $A_a$ correspond to the values $\frac{a}{3}$ and $\frac{2a}{3}$ of $u$. We call $M_a$ the mid-point of this segment \ie, its intersection with the median $[OM]$. More precisely, we consider the functions $BP_{\frac{2}{5}}$ and $BP_{\frac{4}{5}}$. According to the proof of Lemma~\ref{MB23-L2}, we have the formulas
\begin{equation}\label{23C-22}
\begin{array}{l}
BP_{\frac{2}{5}} = - 4 \sin(\frac{\pi}{5})\sin(15\pi u)
\sin(\pi u + \frac{3\pi}{10}-\theta) \text{~for~} u \in [\frac{2}{15},\frac{4}{15}]\,,\\[5pt]
BP_{\frac{4}{5}}= - 4 \sin(\frac{\pi}{10})\sin(15\pi u)
\sin(\pi u + \frac{\pi}{10}-\theta) \text{~for~} u \in [\frac{4}{15},\frac{8}{15}]\,.
\end{array}
\end{equation}

Taking into account the intervals for $u$, and the fact that $0 < \theta \le \frac{\pi}{6}$, we can conclude that the corresponding segments are barriers inside the $H$-triangles which they intersect.

Figure~\ref{NSP23-F2} shows the three possible configurations depending on the number of critical zeros of the function $\Psi_{2,3}^{\theta}$ on the open edge $[OA]$.

The nodal set $N(\Psi_{2,3}^{\theta})$ consists of finitely many regular arcs which can only cross at critical zeros (including at the boundary). Because $\theta \in ]0,\frac{\pi}{6}]$, no nodal arc arrives at a vertex. Because the fixed points are regular points, there is only one nodal arc at each fixed points. Only one nodal arc arrives at a boundary critical zero of order $2$ (except $Z_b$). Exactly two nodal arcs arrive at $Z_b$, with equal angles.

We work in each $H$-triangle separately. When working in a given $H$-triangle $\cH$, we call \emph{exit point} a point at which a nodal arc can exit the triangle. Because medians are partial barriers (Lemma~\ref{MB23-L1}), an exit point is either a fixed point, or a boundary critical zero.

Assume that $Z$ is a critical zero of $\Psi_{2,3}^{\theta}$, in the interior of some $H$-triangle $\cH$. There are at least four semi-arcs emanating from $Z$, and we can follow each one of them. Following such an arc, there are only two possibilities: either we arrive at an exit point, or the path we follow is not injective. The latter necessarily occurs if there are at most three exit directions. However, if a nodal path is not injective, it bounds at least a nodal domain $\omega \subset \cH$, and hence its first Dirichlet eigenvalue satisfies $\lambda(\omega) = \lambda_7 = 19$. On the other-hand, $\lambda(\omega) > \lambda(\cH)=21$ (Lemma~\ref{MB13-23-L}). More generally, we have proved the following property.

\begin{lemma}\label{NSP23-LA}
Assume $\cD \subset \cH$ is bounded by partial barriers, with at most $3$ exit directions. Then $\cD$ cannot contain any critical zero in its interior.
\end{lemma}%

We now make a case by case analysis of the six $H$-triangles contained in $\cT$. For the notation, see Figure~\ref{NSP23-F2}.

\noib Triangle $T(B,F_C,M_A)$. There are two exit points $F_{1,B}$ and $F_{2,A}$, each one with one exit direction. Note that $F_C$ is not an exit point. We can apply Lemma~\ref{NSP23-LA}, and conclude that the nodal set $N(\Psi_{2,3}^{\theta})$ inside this $H$-triangle is an arc from $F_{1,B}$ to $F_{2,A}$, without self-intersections.\medskip

\noib Triangle $T(O,F_C,M_A)$. Same arguments as in the preceding case, the nodal set inside this $H$-triangle is an arc from $F_{2,A}$ to $F_{1,O}$, without self-intersections.\medskip

\noib Triangle $T(O,F_C,M_B)$. In this $H$-triangle, there are $4$ exit points, with one direction each: $F_C$, $F_{1,O}$, $F_{2,B}$, and $Z_{6,\theta}$. The barrier $D_{\frac{2}{5}}$ meets $N(\Psi_{2,3}^{\theta})$ at only one point, $F_{1,O}$. In particular, this segment cannot contain any critical zero and divides the $H$-triangle into two sub-domains with two exit points. We can apply Lemma~\ref{NSP23-LA} to each sub-domain, and conclude that the nodal set inside the $H$-triangle consists of two disjoint arcs without self-intersections, one from $F_{1,O}$ to $Z_{6,\theta}$, and one from $F_C$ to $F_{2,B}$.\medskip

\noib Triangle $T(A,F_C,M_B)$. For this $H$-triangle, we have to consider three cases depending on the sign of $\theta-\theta_c$.

(i) Assume $0 < \theta < \theta_c$. In this case, there are only two exit points in this $H$-triangle, and we can reason as in the case of the triangle $T(B,F_C,M_A)$, concluding that the nodal set inside $T(A,F_C,M_B)$ is an arc from $F_{2,B}$ to $F_{1,A}$, without self-intersections.

(ii) Assume $\theta = \theta_c$. In this case, we have three exist points and $4$ exit directions because $Z_b$ is of order $3$. Using the arguments as in the proof of Lemma~\ref{NSP23-LA}, the existence of an interior critical zero would either yield an interior nodal loop, or a nodal loop touching the boundary at $Z_b$. In either case, the energy argument works, and we can conclude that the nodal set inside $T(A,F_C,M_B)$ consists of two disjoint arcs without self-intersections, one from $Z_b$ to $F_{2,B}$, and one from $Z_b$ to $F_{1,A}$.

(iii) Assume $\theta > \theta_c$. Similar to the preceding case, assuming there is an interior critical zero $Z$, we would either get an interior nodal loop, or a curve from $Z_{4,\theta}$ to $Z$, to $Z_{5,\theta}$. We would then get a nodal domain contained in $T(A,F_C,M_B)$ and with boundary intersecting the edge $[OA]$. The energy argument applies, and we can conclude that the nodal set consists of two disjoint curves without self-intersections, one from $Z_{5,\theta}$ to $F_{2,3}$, and one from $Z_{4,\theta}$ to $F_{1,A}$.\medskip

\noib Triangle $T(A,F_C,M_O)$. There are four exit points, with one direction each. Same treatment as for Triangle $T(A,F_C,M_B)$.\medskip

\noib Triangle $T(B,F_C,M_O)$. There are four  exit points, with one exit direction each. Using the barrier $D_{\frac{4}{5}}$, this triangle can be divided into two sub-domains, each with three exit points. It follows that there are no interior critical zero in this $H$-triangle. It follows that the nodal set consists of two disjoint arcs without self-intersections joining pairs of exit points. This excludes (as the barrier $D_{\frac{4}{5}}$ does) the possibility of two arcs $(F_C,Z_{3,\theta})$ and $(F_{2,O},F_{1,B})$. We are left with two possibilities,

(i) an arc $(F_C,F_{1,B})$ and an arc $(Z_{3,\theta},F_{2,O})$, or

(ii) an arc $(F_C,F_{2,O})$ and an arc $(F_{1,B},Z_{3,\theta})$.

For continuity reasons, (ii) holds for $\theta$ small enough.  In both cases, the number of nodal domains is at most four, and we can conclude that the eigenvalue $\lambda_7$ is not Courant-sharp.

In fact both case (i) and (ii) yield at most four nodal domains. \hfill $\square$

\textbf{Remarks}. \\ (i) Let $\Psi$ be an eigenfunction. Once one knows the critical zeros of $\Psi$, together with their orders, and the number of connected components of $N(\Psi) \cap \partial \cT$, one can apply the Euler-type formula of \cite[Proposition~2.8]{HOMN} to obtain the number of nodal domains of $\Psi$. Using the proof of Proposition~\ref{NSP23-P1} and this formula, one can recover the number of nodal domains of $\Psi_{2,3}^{\theta}$\,: three for $\theta \in\, ]0,\theta_c[$, and four for $\theta \in [\theta_c,\frac{\pi}{6}]$.

(ii) It turns out that case (i) at the end of the previous proof can be discarded by yet another barrier. Indeed, we can consider the segment determined in $\cT$ by the line through $F_{1,O}$ and parallel to $[OA]$, whose equation in the parametrization $\cF$ is $s-2t+\frac{3}{5} = 0$. It is parametrized by
\begin{equation}\label{23C-23a}
u \mapsto \left( 2u-\frac{3}{5},u \right) \text{~for~} u \in [\frac{2}{5},\frac{8}{15}]\,.
\end{equation}
As in the proof of Lemma~\ref{MB23-L2}, we define the functions
\begin{equation}\label{23C-23b}
\begin{array}{l}
EC(u) := C_{2,3}(2u-\frac{3}{5},u)\,,\\[5pt]
ES(u) := S_{2,3}(2u-\frac{3}{5},u)\,,\\[5pt]
EP^{\theta}(u) := \Psi_{2,3}^{\theta}(2u-\frac{3}{5},u)\,,
\end{array}
\end{equation}
and we obtain the relations,
\begin{equation}\label{23C-23c}
\begin{array}{l}
EC(u) = 4 \sin(\frac{\pi}{5}) \sin(15\pi u) \cos(\pi u + \frac{\pi}{5})\,,\\
ES(u) = -4 \sin(\frac{\pi}{5}) \sin(15\pi u) \sin(\pi u + \frac{\pi}{5})\,,\\
EP^{\theta}(u) = -4 \sin(\frac{\pi}{5}) \sin(15\pi u) \cos(\pi u + \frac{\pi}{5}+\theta)\,.
\end{array}
\end{equation}

Figure~\ref{NSP23-F3} shows a Maple simulation with the two bifurcations which occur in the interval $[0,\frac{\pi}{6}]$, namely at $0$ (the inner critical zero disappears) and at $\theta_c$, with a critical zero of order $3$ on the edge $[OA]$.

\vfill

\section{Equilateral triangle: epilogue}

Provided one can give a full description of the nodal sets $N(C_{m,n})$ and $N(S_{m,n})$, for $(m,n) = (1,3)$ or $(2,3)$, one can use the checkerboard argument associated with the pair $(C_{m,n},S_{m,n})$ mentioned in Section \ref{sProl}. Since we work with $\theta \in [0,\frac{\pi}{6}]$, it is actually more appropriate to work with the checkerboard associated with the pair $(C_{m,n},\Psi_{m,n}^{\frac{\pi}{6}})$. Indeed,
$$
\sin(\frac{\pi}{6})\Psi_{m,n}^{\theta} = \sin(\frac{\pi}{6} - \theta) C_{m,n} + \sin(\theta) \Psi_{m,n}^{\frac{\pi}{6}}\,.
$$
We have the inclusion
$$
N(\Psi_{m,n}^{\theta}) \subset \{ C_{m,n}\,\Psi_{m,n}^{\frac{\pi}{6}} < 0 \} \cup \left( N(C_{m,n}) \cap N(S_{m,n})\right).
$$

\begin{figure}[h!]
  \centering
  \includegraphics[width=10cm]{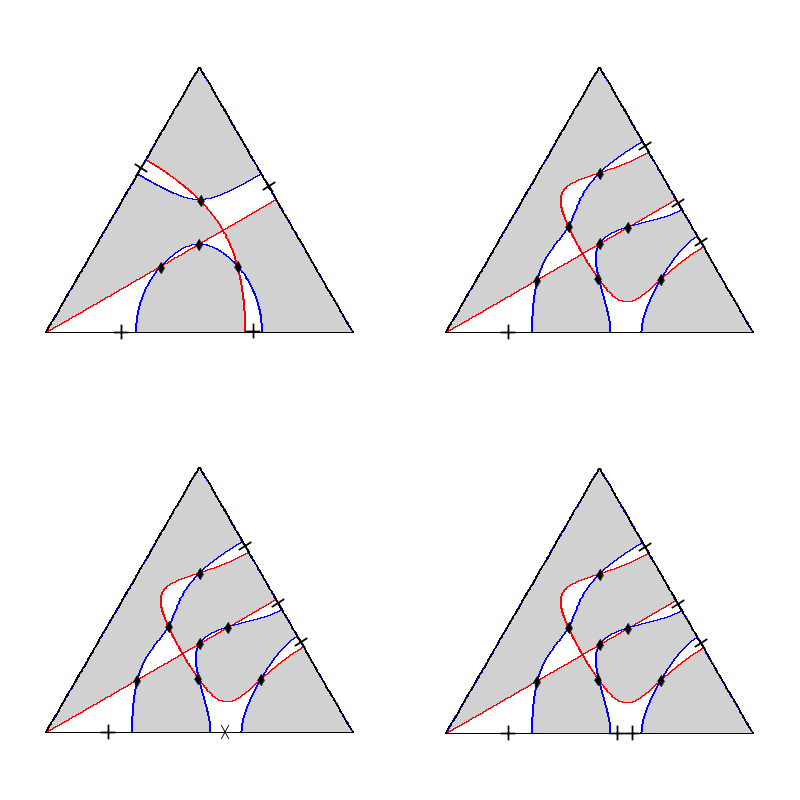}
  \caption{Adapted checkerboards}\label{epi-F1}
\end{figure}

Figure~\ref{epi-F1} displays this checkerboard in the case $(m,n) = (1,3)$ (top row, left), and in the three sub-cases of the case $(m,n) = (2,3)$, depending on the number of critical zeros on the edge $[OA]$ (top row right and bottom row). It is easy to see that the knowledge of the checkerboard and of the critical zeros on the edges determines the nodal set of $\Psi_{m,n}^{\theta}$ when $\theta \in ]0,\frac{\pi}{6}]$. Indeed, recall that the rule is that a nodal line can only leave or enter a white component through points which belong to $N(C_{m,n}) \cap N(S_{m,n})$. But each white component has at most two points of this type in its closure. Note that this was not the case in the checkerboard presented in Figure~\ref{I-F1}, where we need to prove the existence of  additional barriers to determine the nodal picture.

\newpage
\section{ Other triangles}\label{S-OT}

The rectangle, the equilateral triangle, the right-angled isosceles triangle, and the hemi\-equi\-lateral triangle, are the alcoves of the $2$-dimensional root systems $A_1\times A_1$, $A_2$, $B_2$ and $G_2$ respectively, \cite{Bou}. Their Dirichlet and Neumann spectral data can be described completely, \cite{Ber}. They are also the only plane polygonal domains all of whose Dirichlet or Neumann eigenfunctions are trigonometric \cite[Chapt. IV]{McC}.

The Dirichlet  spectral data of the right-angled isosceles triangle and of the hemi\-equi\-lateral triangle are easily deduced from the data for the square and of the equilateral triangle respectively. The investigation of the Courant-sharp Dirichlet eigenvalues of these triangles can also be done, using the above methods. It turns out to be simpler than the investigation of the Courant-sharp eigenvalues of the equilateral triangle. We sketch the proofs in the following sub-sections.

\subsection{The right-angled isosceles triangle}\label{SS-OTB}

Call $\cS_{\pi}$ the open square with side $\pi$. A complete set of Dirichlet eigenfunctions is given by
\begin{equation}\label{OTB-2}
\psi_{m,n}(x,y) = \sin(mx)\sin(ny)\,,
\end{equation}
with associated eigenvalues $m^2+n^2$, where $m$ and $n$ are positive integers.

Call $\cB_{\pi}$ the right-angled isosceles triangle,
\begin{equation}\label{OTB-4}
\cB_{\pi} = \left\lbrace (x,y) \in ]0,\pi[^2 ~|~ y < x \right\rbrace\,.
\end{equation}

A complete set of Dirichlet eigenfunctions for $\cB_{\pi}$ is given by
\begin{equation}\label{OTB-6}
\varphi_{m,n}(x,y) = \sin(mx)\sin(ny) - \sin(nx)\sin(my)\,,
\end{equation}
with associated eigenvalues $m^2+n^2$, for positive integers $m$ and $n$ satisfying $1 \le n \le m-1$.

Denote by $N_{\cS}(\lambda) = \# \{n ~|~ \lambda_n(\cS_{\pi}) < \lambda\}$ the counting function of the Dirichlet eigenvalues of $\cS_{\pi}$. Similarly denote by $N_{\cB}(\lambda)$ the counting function of the Dirichlet eigenvalues of $\cB_{\pi}$. Clearly,
\begin{equation}\label{OTB-8}
N_{\cS}(\lambda) = 2 N_{\cB}(\lambda) + \#\{ m \in \Nb ~|~ 2m^2 < \lambda\}\,.
\end{equation}

Using \cite[Equation (2.1)]{BeHe1}, we get the lower bound
\begin{equation}\label{OTB-10}
N_{\cB}(\lambda) \ge \frac{\pi \lambda}{8} - \frac{(4+\sqrt{2})\sqrt{\lambda}}{4} + \frac{1}{2}\,.
\end{equation}

Using the Faber-Krahn inequality, a Courant-sharp eigenvalue $\lambda_n(\cB_{\pi})$ satisfies
\begin{equation}\label{OTB-12}
\frac{\lambda_n(\cB_{\pi})}{n}  \ge \frac{2}{\pi}j_{0,1}^2 \,  \sim 3.681690532 \,.
\end{equation}

On the other-hand, if $\lambda_n(\cB_{\pi})$ is Courant-sharp, then $N_{\cB}\left( \lambda_n(\cB_{\pi}) \right) = n-1$. Using \eqref{OTB-10} and \eqref{OTB-12}, we obtain the inequality
\begin{equation}\label{OTB-14}
\frac{2}{\pi}j_{0,1}^2 \, n \le \left( \frac{4+\sqrt{2}}{\pi} \right)^2 \left( 1 + \sqrt{1+ \frac{8\pi(n-\frac{3}{2})}{(4+\sqrt{2})^2}}\right)^2 \,.
\end{equation}

Equation \eqref{OTB-14} shows that if $\lambda_n(\cB_{\pi})$ is Courant-sharp, then $ n \le 26$. Calculating the ratios $\frac{\lambda_n(\cB_{\pi})}{n}$ and using \eqref{OTB-12}, we see that the only possible Courant-sharp Dirichlet eigenvalues of $\cB_{\pi}$ are the $\lambda_n(\cB_{\pi})$ for $n \in \{1, \ldots, 7, 9, 10\}$. These eigenvalues have multiplicity $1$, and correspond to the pairs
\begin{equation}\label{OTB-16}
\left\lbrace [2,1], [3,1], [3,2], [4,1], [4,2], [4,3], [5,1], [5,3], [6,1] \right\rbrace \,.
\end{equation}

To conclude, we have to determine the number $\mu(\varphi_{m,n})$ of nodal domains of the corresponding eigenfunctions $\varphi_{m,n}$.

\textbf{Remark}. The eigenvalues $\lambda_5, \lambda_7$ and $\lambda_9$ correspond to the pairs $[4,2], [5,1]$ and $[5,3]$ respectively. The associated eigenfunctions are anti-symmetric with respect to the line $\{x+y=\pi\}$. It follows that they have an even number of nodal domains. These eigenvalues cannot be Courant-sharp.

To deal with the other eigenvalues we use the following result \cite{ABFG}, and the description of some nodal sets of Dirichlet eigenfunctions for the square membrane \cite{BeHe1}.

\begin{properties}\label{OTB-P1}
The Dirichlet eigenfunctions $\varphi_{m,n}$ of the triangle $\cB_{\pi}$ satisfy the following identities.\vspace{-3mm}
\begin{enumerate}
  \item For $m > n$, $$\varphi_{m+n,m-n}(x,y) = \varphi_{m,n}(x+y,x-y)\,.$$
  \item If $d$ is the greatest common divisor of $m$ and $n$, then $$\varphi_{m,n}(x,y) = \varphi_{\frac{m}{d},\frac{n}{d}}(dx,dy)\,.$$
\end{enumerate}
\end{properties}%

We obtain the following nodal count.

\begin{equation}
\begin{array}{|c|c|c|c|}
\hline
\lambda_i(\cB_{\pi} & [m,n] & \mu(\varphi_{m,n}) & \text{~using~} \\ \hline
\lambda_1 & [2,1] & 1 &  \\ \hline
\lambda_2 & [3,1] & 2 &  \\ \hline
\lambda_3 & [3,2] & 2 &  \\ \hline
\lambda_4 & [4,1] & 2 &  \\ \hline
\lambda_5 & [4,2] & 4 & \varphi_{2,1} \\ \hline
\lambda_6 & [4,3] & 3 & \varphi_{7,1} \\ \hline
\lambda_7 & [5,1] & 4 &  \\ \hline
\lambda_9 & [5,3] & 4  & \varphi_{8,2}, \varphi_{4,1} \\ \hline
\lambda_{10} & [6,1] & 3 &  \\
\hline
\end{array}
\end{equation}

\begin{figure}[hbt]
  \centering
  \includegraphics[width=14cm]{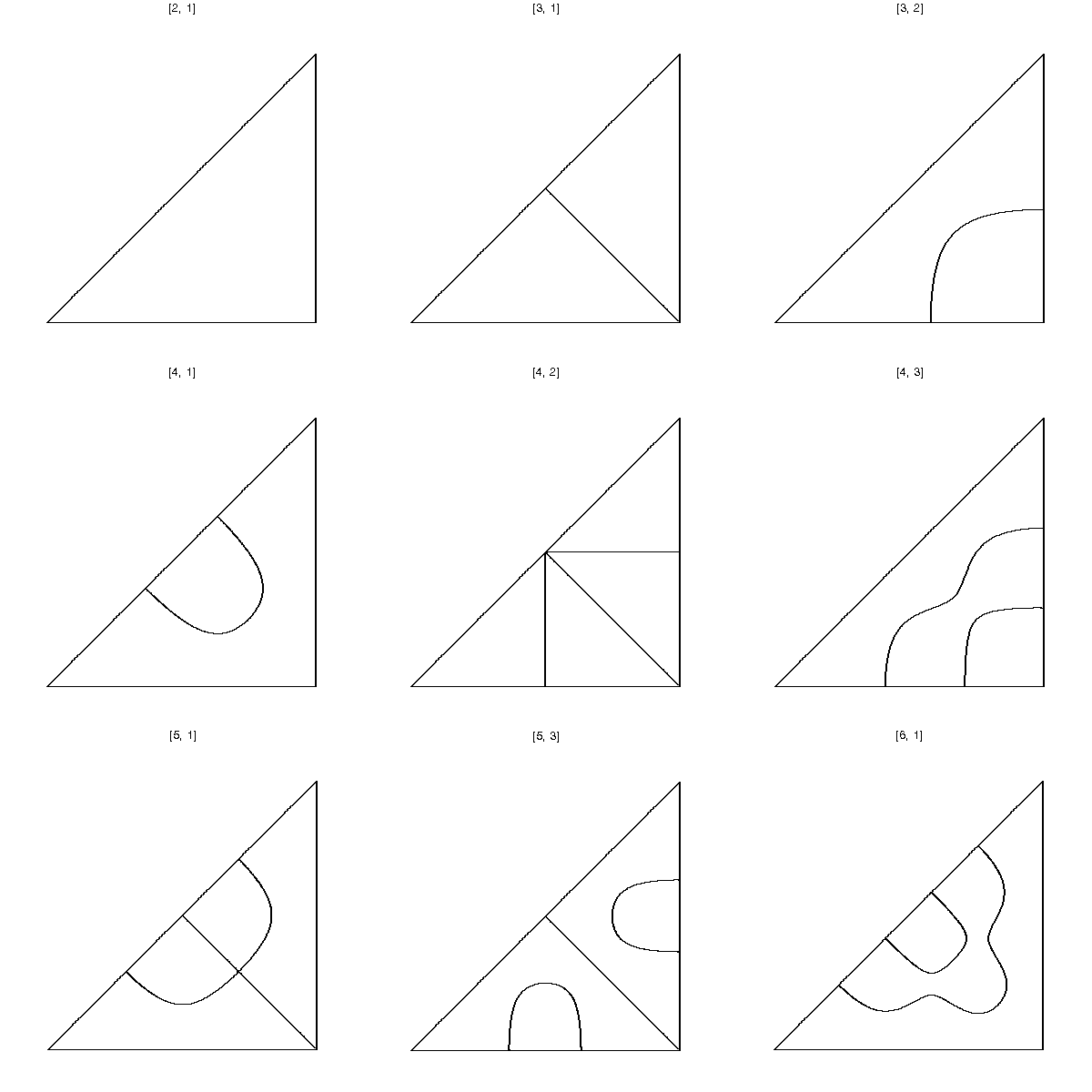}
  \caption{Courant-sharp analysis for the right-angled isosceles triangle}\label{OTB-F1}
\end{figure}

\begin{theorem}\label{OTB-P2}
The only Courant-sharp eigenvalues of the right-angled triangle $\cB_{\pi}$ are $\lambda_1(\cB_{\pi})$ and $\lambda_2(\cB_{\pi})$.
\end{theorem}

\textbf{Remark}. A general algorithm to compute the number of nodal domains of the eigenfunctions $\varphi_{m,n}$ of the triangle $\cB_{\pi}$ is described in \cite{ABFG}.

\subsection{The hemiequilateral triangle}\label{SS-OTH}

A complete set of Dirichlet eigenfunctions of the hemiequilateral triangle $\cH$, with hypothenuse of length $1$, is given by the functions $C_{m,n}$ described in \eqref{ETR-8a}, with $1 \le n < m$. The associated eigenvalues are the numbers $\frac{16\pi^2}{9}(m^2+mn+n^2).$

Using the Faber-Krahn inequality, we obtain that a Courant-sharp eigenvalue $\lambda_n(\cH)$ satisfies the inequality
\begin{equation}\label{OTH-2}
\lambda_n(\cH) \ge \frac{8\pi}{\sqrt{3}} j_{0,1}^2 \, n.
\end{equation}

The lower bound \eqref{ETR-12} yields the following lower bound for the counting function of the Dirichlet eigenvalues of $\cH$,
\begin{equation}\label{OTH-4}
N_{\cH}(\lambda) \ge \frac{\sqrt{3}}{32\pi}\lambda - \frac{6+\sqrt{3}}{8\pi}\sqrt{\lambda} + \frac{1}{2}\,.
\end{equation}

Using \eqref{OTH-2} and \eqref{OTH-4}, we conclude that if $\lambda_n(\cH)$ is Courant-sharp, then $n \le 32$. Using \eqref{OTH-2} again, we see that the only possible Courant-sharp eigenvalues are $\lambda_i(\cH)$, for $i \in \{1, \ldots, 8,10\}$. All these eigenvalues have multiplicity $1$, and correspond to the pairs
$$
\left\lbrace [2,1], [3,1], [3,2], [4,1], [4,2], [5,1], [4,3], [5,2], [5,3] \right\rbrace \,.
$$
Since the eigenfunctions are given explicitly, we can numerically determine the nodal sets, and conclude.

\begin{figure}[bth]
  \centering
  \includegraphics[width=13cm]{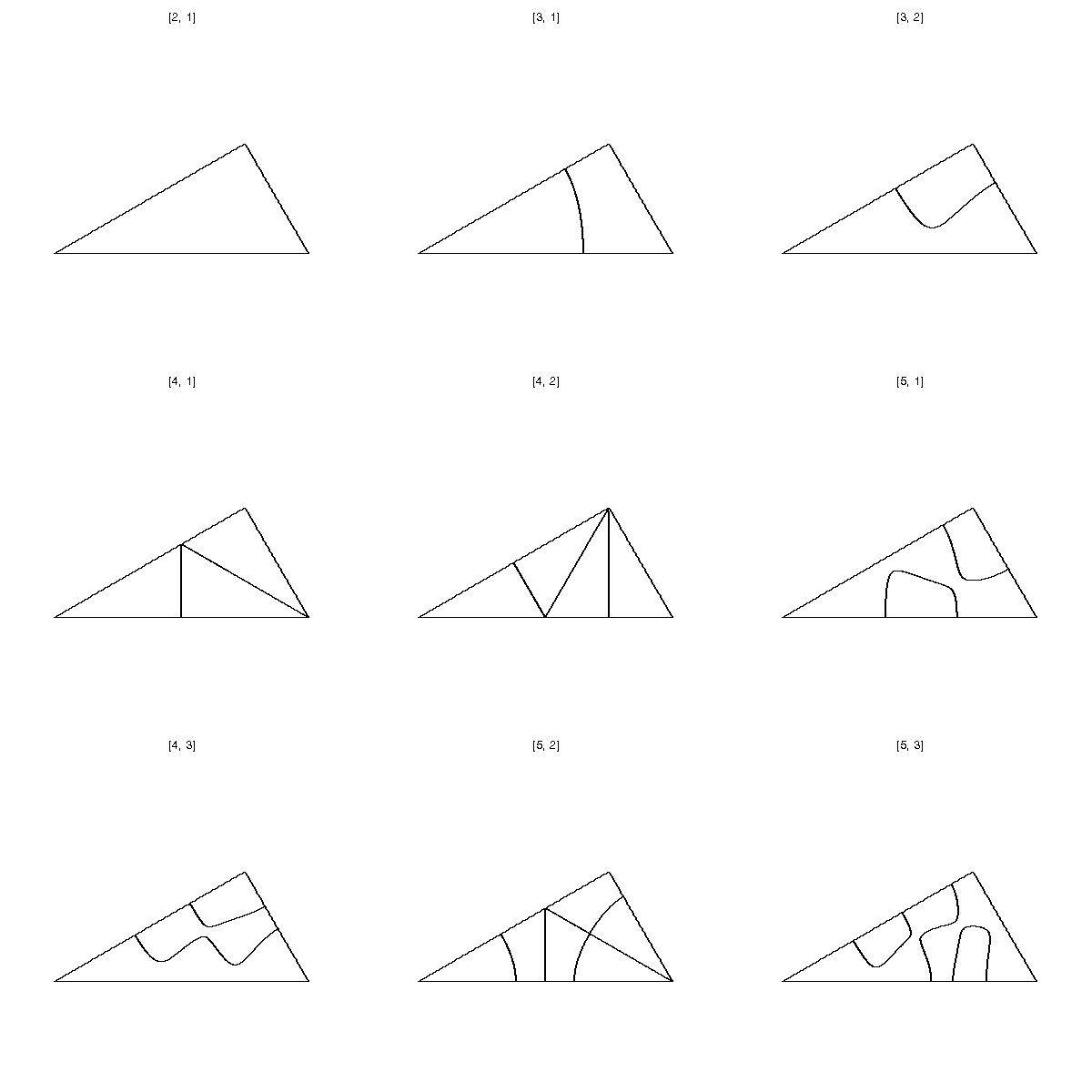}
  \caption{Courant-sharp analysis for the hemiequilateral triangle}\label{OTH-F1}
\end{figure}

\begin{theorem}\label{OTH-P1}
The only Courant-sharp Dirichlet eigenvalues of the hemiequilateral triangle $\cH$ are the first and second.
\end{theorem}%

\bigskip


\end{document}